\documentclass[11pt]{article}

\usepackage{amsmath}
\usepackage{amsfonts}
\usepackage{graphicx}
\usepackage{setspace}
\usepackage{amsmath}
\usepackage{amssymb}
\usepackage{latexsym}
\usepackage{amsmath, amsfonts,amssymb, amsthm, euscript,makeidx,color,mathrsfs}

\usepackage[numbers,sort&compress]{natbib}

\def\sqr#1#2{{\vcenter{\vbox{\hrule height.#2pt
              \hbox{\vrule width.#2pt height#1pt \kern#1pt \vrule width.#2pt}
              \hrule height.#2pt}}}}
\usepackage{hyperref}
\hypersetup{
    colorlinks=true,
    linkcolor=blue,
    urlcolor=cyan,
    citecolor=cyan
}
\usepackage{caption}

\RequirePackage[capitalize,nameinlink]{cleveref}

\crefname{section}{section}{sections}
\crefname{subsection}{subsection}{subsections}
\Crefname{section}{Section}{Sections}
\Crefname{subsection}{Subsection}{Subsections}

\crefname{condition}{Condition}{Conditions}

\Crefname{figure}{Figure}{Figures}

\crefformat{equation}{\textup{#2(#1)#3}}
\crefrangeformat{equation}{\textup{#3(#1)#4--#5(#2)#6}}
\crefmultiformat{equation}{\textup{#2(#1)#3}}{ and \textup{#2(#1)#3}}
{, \textup{#2(#1)#3}}{ and \textup{#2(#1)#3}}
\crefrangemultiformat{equation}{\textup{#3(#1)#4--#5(#2)#6}}%
{ and \textup{#3(#1)#4--#5(#2)#6}}{, \textup{#3(#1)#4--#5(#2)#6}}{ and \textup{#3(#1)#4--#5(#2)#6}}

\Crefformat{equation}{#2Equation~\textup{(#1)}#3}
\Crefrangeformat{equation}{Equations~\textup{#3(#1)#4--#5(#2)#6}}
\Crefmultiformat{equation}{Equations~\textup{#2(#1)#3}}{ and \textup{#2(#1)#3}}
{, \textup{#2(#1)#3}}{ and \textup{#2(#1)#3}}
\Crefrangemultiformat{equation}{Equations~\textup{#3(#1)#4--#5(#2)#6}}%
{ and \textup{#3(#1)#4--#5(#2)#6}}{, \textup{#3(#1)#4--#5(#2)#6}}{ and \textup{#3(#1)#4--#5(#2)#6}}

\crefdefaultlabelformat{#2\textup{#1}#3}

\newtheorem {theorem}{Theorem}[section]
\newtheorem {lemma}[theorem]{{\bf Lemma}}
\newtheorem {corollary}[theorem]{{\bf Corollary}}
\newtheorem {proposition}[theorem]{{\bf Proposition}}
\theoremstyle{remark}
\newtheorem {remark}{{\bf Remark}}[section]

\theoremstyle{definition}

\theoremstyle{plain} \numberwithin {equation}{section}

\numberwithin{assumption}{section}

\allowdisplaybreaks

\def\no{\noindent}

\usepackage{graphicx}
\usepackage{ifpdf}
\ifpdf
  \usepackage{epstopdf}
\fi

\usepackage{enumerate}

\def\deq{\mathop{\buildrel\Delta\over=}}

\begin{document}

\title{Null controllability for stochastic
 fourth order semi-discrete parabolic equations\footnote{This work is supported by the Fundamental Research Funds for the Central Universities 2682024CX013.}}
\author{
	Yu Wang\footnote{School of Mathematics, Southwest Jiaotong University, Chengdu, P. R. China.
	Email: yuwangmath@163.com.}
	~~~ and ~~~
	Qingmei Zhao\footnote{Corresponding author. School of Mathematics Sciences, Sichuan Normal University, Chengdu, P. R. China. E-mail:
	qmmath@163.com.}
} 

\date{}

\maketitle

\begin{abstract}

	This paper is devoted to studying  null controllability  for a class of  stochastic  fourth order  semi-discrete parabolic equations, where the spatial variable is discretized with finite difference scheme and the time is kept as a continuous variable.
	For this purpose, we establish a new global Carleman estimate for a  backward stochastic fourth order semi-discrete parabolic operators, in which the large parameter is connected to the mesh size. 
	A relaxed observability estimate is established for backward stochastic  fourth order semi-discrete parabolic equations by this new Carleman estimate,  with an explicit observability constant that depends on the discretization parameter and coefficients of lower order terms. 
	Then, the $\phi$-null  controllability of the stochastic fourth order semi-discrete parabolic equations is proved using the standard duality technique.
				
\end{abstract}

\bigskip

\no{\bf Mathematics Subject Classification}.  Primary 93B05; Secondary 93B07, 93C20
\bigskip

\no{\bf Key Words}.  Stochastic fourth order semi-discrete parabolic equation,  null controllability, observability,  global Carleman
estimate.

\section{Introduction}

Let $T>0$ and $(\Omega, \mathcal{F}, \mathbf{F},
\mathbb{P})$ with $\mathbf F=\{\mathcal{F}_{t}\}_{t \geq 0}$ be a complete filtered probability space on which a one-dimensional standard Brownian motion $\{W(t)\}_{t \geq 0}$ is defined and $\mathbf{F}$ is the natural filtration generated by $W(\cdot)$,  augmented by all the $\mathbb{P}$ null sets in $\mathcal{F}$. 
Write $ \mathbb{F} $ for the progressive $\sigma$-field with respect to $\mathbf{F}$. 
 Let $\mathcal{H}$ be a Banach space. Denote by $L^2_{\mathcal{F}_\tau}(\Omega; \mathcal{H})$ the space of all $\mathcal{F}_\tau$-measurable random variables $\xi$ such that $\mathbb{E} | \xi |^2_{\mathcal{H}}<\infty$; by $L^2_{\mathbb{F}}(0,T;\mathcal{H})$ the
Banach space consisting of all $\mathcal{H}$-valued
$\{\mathcal{F}_t\}_{t\geq 0}$-adapted processes $X(\cdot)$ such that
$\mathbb{E}(|X(\cdot)|^2_{L^2(0,T;\mathcal{H})})<\infty$, with the
canonical norm; by $L^\infty_{\mathbb{F}}(0,T;\mathcal{H})$ the
Banach space consisting of all $\mathcal{H}$-valued
$\{\mathcal{F}_t\}_{t\geq 0}$-adapted essentially bounded processes;
and by $L^2_{\mathbb{F}}(\Omega; C([0,T];\mathcal{H}))$ the Banach
space consisting of all $\mathcal{H}$-valued
$\{\mathcal{F}_t\}_{t\geq 0}$-adapted continuous processes
$X(\cdot)$ such that
$\mathbb{E}(|X(\cdot)|^2_{C([0,T];\mathcal{H})})<\infty$.

Let $ G = (0,1) $ and $ G_{0} $ be a nonempty open subset of $ G $.
As usual, $ \chi_{G_{0}} $   denotes the characteristic function of $ G_{0} $. Firstly, we introduce the following controlled fourth order  stochastic  parabolic equation in a continuous framework:
\begin{align}
	\label{eqContinuousEq}
	\begin{cases}
		d y + \partial_{x}^{4} y d t = (a_{1} y + \chi_{G_{0}} u ) d t + (a_{2} y + v) d W(t) & \text { in } (0,T) \times G, \\ 
		y(t,0) = y(t,1) = 0 & \text { on } (0,T), \\ 
		\partial_{x} y(t,0) = \partial_{x} y(t,1) = 0 & \text { on } (0,T), \\ 
		y(0) = y_{0} & \text { in } G 
		.
	\end{cases}
\end{align}
Here, the initial datum $ y_{0} \in  L^2(G) $, the coefficients $ a_{1}, a_{2} $ in $ L^{\infty}_{\mathbb{F}}( 0,T;L^{\infty}(G) ) $ and the control $(u, v) \in L_{\mathbb{F}}^{2}(0, T ; L^{2}(G_{0})) \times$ $L_{\mathbb{F}}^{2}(0, T ; L^{2}(G))$. 
By the classical well-posedness result for stochastic evolution equations, we know that
\cref{eqContinuousEq}  admits a unique weak solution $y\in L^2_{\mathbb{F}}(\Omega;C([0,T];L^2(G)))\times L^2_{\mathbb{F}}(0,T;H^2(G)\cap H_0^1(G))$ (e.g., \cite[Theorem 3.24]{Lue2021a}).
The system
(\ref{eqContinuousEq}) is said to be null controllable at time $T$ if for any given $y_0$, there exist controls $(u,v)\in   L_{\mathbb{F}}^{2}(0, T ; L^{2}(G_{0})) \times$ $L_{\mathbb{F}}^{2}(0, T ; L^{2}(G))$ such
that  the corresponding solution satisfies $y(x, T) = 0$, $\mathbb{P}$-a.s.

Stochastic fourth order parabolic equation in a continuous framework is widely used to describe some fast diffusion phenomenon
under stochastic perturbations (e.g.,\cite{Cook1970}). Further, it is the linearized version of the stochastic Cahn-Hilliard equation, which is used to model phase separation
phenomena in a melted alloy that is quenched to a temperature  (e.g.,\cite{Cardon2001}).
Compared to second order parabolic equation, to derive controllability for higher
order parabolic equation is not an easy task even in the continuous framework. 
To the best of our knowledge,
the controllability of deterministic and stochastic fourth order parabolic equations has been extensively studied in
many work (see \cite{Cerpa2011,Gao2015,Guerrero2019,Guzman2020,Lv2022,Zhou2012} and the references cited therein). 


As the numerical approximation of this control system, it is natural to consider the null controllability of discrete systems.

For $ N \in \mathbb{N} $, put the space discretization parameter $ h \deq 1/(N+1) $.
We consider the pairs $ (t, x_{i}) $ with $ t \in (0,T) $ and $ x_{i} = i h $ for $ i = 0, \cdots, N+1 $. 
Utilizing the centered finite difference method to the space variable for the system \cref{eqContinuousEq}, and considering two additional points $ x_{-1} = - h $ and $ x_{N+2} = (N+2) h$ to discretize the boundary conditions, we consider the following stochastic semi-discrete system:
\begin{align}
	\label{eqDiscreteEq}
	\begin{cases}
		\begin{aligned}
			&d y^{i} + \frac{1}{4} ( y^{i+2} - 4 y^{i+1} + 6 y^{i} - 4 y^{i-1} + y^{i-2}) d t 
		,\\ 
		&\quad \quad 
		= (a^{i}_{1} y^{i} + \chi^{i}_{G_{0}} u^{i} ) d t + (a^{i}_{2} y^{i} + v^{i}) d W(t) & & t \in  (0,T), \quad i=1,\cdots,N   ,\\ 
		&y^{0}(t) = y^{N+1}(t) = 0 &&  t \in  (0,T) ,\\ 
		&y^{-1}(t) = y^{N+2}(t) = 0 && t \in (0,T) ,\\ 
		&y^{i}(0) = y^{i}_{0} & &i=1,\cdots,N .
		\end{aligned}
	\end{cases}
\end{align}

It is well known that the  uniform null controllability property is  hard to address. 
For instance, as demonstrated by the counterexample provided in \cite{Zua2006}, the uniform null controllability property may not hold even for the deterministic second order semi-discrete equation in the two dimensional case. 
To the best of our knowledge, there are also only a few results on the uniform null controllability of deterministic second order semi-discrete parabolic equations for some special cases; see \cite{Lopez1998, Labbe2006, Nguyen2015, Allonsius2018}.



Therefore, there are papers consider other controllability rather than uniform null controllability.
In \cite{BoyerJMPA,BoyerSICON}, a uniform null controllability result for the lower part of the spectrum of deterministic   semi-discrete 
second order
parabolic problem is established.
Later,  \cite{Boyer2014} proved that the deterministic second order semi-discrete  parabolic equations are $\phi$-null controllable by establishing a discrete global Carleman estimate directly for semi-discrete  parabolic operator.
This $\phi$-null  controllability is defined as the norm of the semi-discrete solution at time $ T $ can be dominated by $ \sqrt{\phi(h)}$, where $\phi$ is a real-valued function that tends to zero as the space discretization parameter $h$ tends to zero.



Since then, there have been many results in the study of controllability problems and inverse problems of deterministic semi-discrete  parabolic systems by using discrete Carleman estimates. 
We refer readers to \cite{Nguyen2014,Cerpa2022,Lecaros2021}.
However, as far as we know, \cite{Zhao2024} is only one paper which are concerned with controllability of the stochastic second order semi-discrete parabolic equations.



To state our main results, we first introduce some notations. 
We define the following regular partition of the interval $[0,1]$ as
$ \mathcal{K} \deq  \{x_i \mid x_i \mid i=0,1, \ldots, N, N+1 \} $.
Considering any set of points $\mathcal{W} \subset \mathcal{K}$, we define the dual meshes $\mathcal{W}^{\prime}$ and $\mathcal{W}^*$, respectively, by
\begin{align*}
	\mathcal{W}^{\prime}\deq\tau_{+}(\mathcal{W}) \cap \tau_{-}(\mathcal{W}), \quad \mathcal{W}^*\deq\tau_{+}(\mathcal{W}) \cup \tau_{-}(\mathcal{W})
\end{align*}
where
\begin{align*}
	\tau_{ \pm}(\mathcal{W})\deq\bigg \{ x \pm \frac{h}{2}  \bigg \rvert\, x \in \mathcal{W}\bigg\}
\end{align*}
We denote $\overline{\mathcal{W}}\deq\left(\mathcal{W}^*\right)^*$ and $\mathring{\mathcal{W}}\deq\left(\mathcal{W}^{\prime}\right)^{\prime}$. 
Note that for two consecutive points $x_i, x_{i+1} \in \mathcal{W}$, we have $x_{i+1} - x_i = h$ provided $\mathring{\overline{\mathcal{W}}} = \mathcal{W}$.
We call such a subset $\mathcal{W} \subset \mathcal{K}$ a regular mesh.
Define the boundary of a regular mesh $\mathcal{W}$ as $\partial \mathcal{W}\deq\overline{\mathcal{W}} \backslash \mathcal{W}$. 

Then, we introduce the semi-discrete sets  $\mathcal{Q} \deq (0, T)  \times \mathcal{W}$. Moreover, we say that the semi-discrete set is regular if the space variable is a regular mesh.
We also define the dual semi-discrete sets by $\mathcal{Q}^{\prime}\deq  (0, T) \times \mathcal{W}^{\prime} $ and $\mathcal{Q}^* \deq (0, T) \times \mathcal{W}^* $. 
Similarly, the semi-discrete boundary is given by $\partial \mathcal{Q}= (0, T) \times \partial \mathcal{W} $.

Next, we define the average operator $A_h$ and the difference operator $D_h$ by
\begin{align*}
	& A_h(u)(t, x)\deq\frac{\tau_{+} u(t, x)+\tau_{-} u(t, x)}{2}, \\
& D_h(u)(t, x)\deq\frac{\tau_{+} u(t, x)-\tau_{-} u(t, x)}{h},
\end{align*}
where $\tau_{ \pm} u(t, x) \deq u\big(t, x \pm \frac{h}{2}\big)$.

Here and in the sequel, $L(\mathcal{W})$ denotes the set of real-valued functions defined in $\mathcal{W}$. 
We also denote by $L_h^2(\mathcal{W})$ the Hilbert space with the inner product 
\begin{align*}
	\langle u, v\rangle_{L_h^2(\mathcal{W})}\deq  \int_{\mathcal{W}} u v\deq h \sum_{x \in \mathcal{W}} u(x) v(x)
\end{align*}
For $ u \in L(\mathcal{W}) $, we define its $ L^{\infty}_{h}(\mathcal{W}) $-norm as 
\begin{align*}
	|u|_{L^{\infty}_{h}(\mathcal{W})} \deq \max_{x \in \mathcal{W}} \{ |u(x) | \}.
\end{align*}

With above notations, defining $ Q \deq (0,T) \times \mathcal{M} $ where $ \mathcal{M} \deq \mathring{\mathcal{K}} $, the controlled semi-discrete system \cref{eqDiscreteEq} can be written as 
\begin{align}
	\label{eqDiscreteEqOpe}
	\begin{cases}
		d y + D_{h}^{4} y  d t  = (a_{1} y + \chi_{G_{0}} u ) d t + (a_{2} y + v) d W(t) & \text { in } Q,\\ 
		y(t,0) = y(t,1) = 0 & \text { on } (0,T), \\ 
		D_{h} y\big(t, - \frac{h}{2}\big) = D_{h} y\big(t, 1+ \frac{h}{2}\big) = 0 & \text { on } (0,T), \\ 
		y(0, x) = y_{0} & \text { in } \mathcal{M},
	\end{cases}
\end{align}
where the initial datum $ y_{0} \in  L^2_{h}(\mathcal{M}) $, the coefficients $ a_{1}, a_{2} $ in $ L^{\infty}_{\mathbb{F}}(0,T;L^{\infty}_{h}(\mathcal{M})) $ and the control $(u, v) \in  L^{2}_{\mathbb{F}}(0,T;   L^{2}_{h}(G_{0} \cap \mathcal{M})) \times L^{2}_{\mathbb{F}}(0,T;L^{2}_{h}( \mathcal{M}))$. 
According to the classical theory of stochastic differential equations (see \cite[Theorem 3.2]{Lue2021a}),  we know that there is a unique solution  $ y \in L_{\mathbb{F}}^{2}(\Omega ; C([0, T] ; L^{2}_{h}(\mathcal{M}))) $ to \cref{eqDiscreteEqOpe}.

The main null controllability result of this paper is the following.

\begin{theorem}
	\label{thmCon}
	There exist positive constants $ h_{0} $ and $ C  $ such that for all $ T > 0 $ and  $ h \leq \min\{h_{0}, h_{1}\} $ with
	\begin{align*}
		h_{1} = C (1+T^{-1} + \mathcal{H}^{2/7})^{-1},
	\end{align*}
	where $ \mathcal{H} = |a_{1}|_{L^{\infty}_{\mathbb{F}}(0,T;L^{\infty}_{h}(\mathcal{M}))} +  |a_{2}|_{L^{\infty}_{\mathbb{F}}(0,T;L^{\infty}_{h}(\mathcal{M}))} $,
	there exist $ (u,v) \in L^{2}_{\mathbb{F}}(0,T;L^{2}_{h}(G_{0} \cap \mathcal{M})) \times L^{2}_{\mathbb{F}}(0,T;L^{2}_{h}( \mathcal{M}))$ such that the solution $ y $ to \cref{eqDiscreteEqOpe} satisfies 
	\begin{align*}
		\mathbb{E} \int_Q|v|^2 d t
		+\mathbb{E} \int_0^T \int_{G_{0} \cap \mathcal{M}}|u|^2 d t \leq C_{obs}    \int_{\mathcal{M}}\left|y_0\right|^2,
	\end{align*}
	and 
	\begin{align}
		\label{eqNull}
		\mathbb{E} \int_{\mathcal{M}}|y(T)|^2 \leq C_{obs} e^{-\frac{C}{h}}  \int_{\mathcal{M}}\left|y_0\right|^2 .
	\end{align}
	Here, $ C_{obs} = e^{C(1+T^{-1}+\mathcal{H}^{2/7}+T+\mathcal{H}^{2}T)} $.
\end{theorem}

Here and in what follows, unless otherwise stated, $ C $ stands for a generic positive constant depending only
on $ G $ and $ G_{0} $ whose value may vary from line to line.

\begin{remark}
	Choosing $ \phi(h) $ equals the right hand side of \cref{eqNull}, we have
	$|y(T)|_{L^{2}_{\mathcal{F}_{T}}(\Omega;L^{2}_{h}(\mathcal{M}))}\le \sqrt{\phi(h)}$ with $\lim\limits_{h\rightarrow 0} \phi(h)=0$, which means that   system \cref{eqDiscreteEqOpe} is $ \phi $-null controllable at time $ T $. 
	It would be quite interesting to study the uniform null controllability; i.e., for all $ h > 0 $, there exist $ (u, v) $ which have uniform bounds, such that  the corresponding solution to \cref{eqDiscreteEqOpe}  satisfies $y(T,x) = 0$, $\mathbb{P}$-a.s for any $x\in \mathcal{M}$.
	It has been done for the deterministic  semi-discrete heat equation on a uniform mesh in \cite{Allonsius2018}.
	However, we currently do not have a method for achieving this for the stochastic case.
\end{remark}

\begin{remark}
	Similar to proving \cref{thmCon}, employing the techniques presented in this paper allows us to obtain $ \phi $-null controllability results for the backward stochastic fourth order semi-discrete  parabolic equation.
\end{remark}

\begin{remark}
	Due to randomness, extra control is needed in the diffusion term, as demonstrated in \cite{Tang2009}. 
	Therefore, the control of \cref{eqDiscreteEqOpe}  is  a pair of function $ (u, v) $.
	It is an interesting question to study the $ \phi $-null controllability for system \cref{eqDiscreteEqOpe} with only one control in the drift. 
\end{remark}

To prove the controllability result of \cref{eqDiscreteEqOpe}, we introduce the adjoint equation:
\begin{align}
	\label{eqDiscreteEqAdjoint}
	\begin{cases}
		d z - D_{h}^{4} z  d t  = - (a_{1} z + a_{2} Z ) d t + Z d W(t) & \text { in } Q,\\ 
		z(t,0) = z(t,1) = 0 & \text { on } (0,T), \\ 
		D_{h} z\big(t, - \frac{h}{2}\big) = D_{h} z\big(t, 1+ \frac{h}{2}\big) = 0 & \text { on } (0,T), \\ 
		z(T, x) = z_{T} & \text { in } \mathcal{M}.
	\end{cases}
\end{align}
Due to the classical theory of backward stochastic differential equations (see \cite[Theorem 4.2]{Lue2021a}),  we know that there is a unique solution  $(z, Z) \in L_{\mathbb{F}}^{2}(\Omega ; C([0, T] ; L^{2}_{h}(\mathcal{M}))) \times L_{\mathbb{F}}^{2}(0, T ; L^{2}_{h}(\mathcal{M}))$ to \cref{eqDiscreteEqAdjoint}.

We have the following observability estimate
\begin{theorem}
	\label{thmObservabilityEstiames}
	There exist positive constants  $ h_{0} $ and $ C   $ such that for all $ T > 0 $ and  $ h \leq \min\{h_{0}, h_{1}\} $,
	the solution $ (z,Z) $ to system \cref{eqDiscreteEqAdjoint} satisfies 
	\begin{align*}
		\mathbb{E} \int_{\mathcal{M}} |z(0)|^{2} \leq 
		C_{obs} \Big(
			\mathbb{E} \int_{Q} |Z|^{2} d t 
			+ \mathbb{E} \int_{0}^{T} \int_{G_{0}\cap \mathcal{M}} |z|^{2} d t 
			+ e^{- \frac{C}{h}} \int_{\mathcal{M}} |z|^{2} \Big |_{t=T}
		\Big)
		,
	\end{align*}
	where $ h_{1} $ and $ C_{obs} $ is given in \cref{thmCon}.
\end{theorem}

\begin{remark}
Note that the last term in the right hand side of observability estimate in \cref{thmObservabilityEstiames} is novel
compared to the stochastic continuous case.	And so far, we don't know how to remove it.
\end{remark}

We will use discrete Carleman estimates to prove \cref{thmObservabilityEstiames}.
To state our Carleman estimate, we first introduce the weight functions.
Let $ \widetilde{G} $ be an open interval contains $ [0,1] $, $ G_{1} \subset  \overline{G_{1}} \subset G_{0} $ and $ G_{2} \subset  \overline{G_{2}} \subset G_{1} $  be two open sets.
By \cite[Lemma 1.1]{Fursikov1996}, we know that there exists a function $ \psi \in C^{4}(\widetilde{G})  $ satisfies
\begin{equation}
	\label{eqPsi}
	\psi >0 \text{ in } \widetilde{G}, 
	\quad| \partial_{x} \psi| \geq C_{0} >0 \text { in } \overline{ G \backslash G_{2}},
	\quad \partial_{x} \psi(0) > 0,
	  \text{ and }
	\partial_{x}\psi(1) < 0
	.
\end{equation}
Furthermore, we assume $ \overline{\mathcal{M} \cap G_{2}} \subset \mathcal{M} \cap G_{1} $.
For any positive parameters $ \lambda $ and $ \mu $, choosing $ 0 < \delta < \frac{1}{2} $, put 
\begin{equation}
	\label{eqAlpha} 
	\begin{cases}
		 \phi(x) = e^{\mu(|2\psi|_{C(\overline{G})}+ \psi(x))} - e^{4 \mu |\psi|_{C(\overline{G})}},
		\quad \varphi(x) = e^{\mu(|2\psi|_{C(\overline{G})}+ \psi(x))},
		\\  
		 s (t) = \lambda \theta(t), 
		\quad \theta(t) = \frac{1}{(t+\delta T)(T+\delta T-t)},
		\\
		 r(t,x) = e^{s(t) \phi(x)},
		\quad \rho (t,x) = (r(t,x))^{-1}.
	\end{cases}
\end{equation}

We have the following Carleman estimate.

\begin{theorem}
	\label{thmCarleman}
	Let the function $s$ and $ \phi $ be defined according to \cref{eqAlpha}. 
	There exists a positive constant $ \mu_{0} $ such that for any $ \mu \geq \mu_{0} $, one can find positive constants $ C = C(\mu) $, $  h_0  $ and $ \varepsilon_{0} $ with $ 0< \varepsilon_{0} \leq 1 $ and a sufficiently large $ \lambda_{0} \geq 1 $
	such that for all $\lambda \geq \lambda_0 \left(T+T^2\right), 0< h \leq h_0$ and $\lambda h\left(\delta T^2\right)^{-1} \leq \varepsilon_0$, it holds that 
	\begin{align}
		\notag
		\label{eqFinalCarleman}
		& 
		\mathbb{E}\int_Q s^7 e^{2 s \phi}|w|^2 d t
		+ \mathbb{E} \int_{Q^*} s^5 e^{2 s \phi}\left|D_h w\right|^2 d t
		+ \mathbb{E} \int_{\overline{Q}} s^3 e^{2 s \phi}\left|D_h^2 w\right|^2 d t
		+ \mathbb{E} \int_{Q^*} s e^{2 s \phi}\left|D_h^3 w\right|^2 d t 
		\\ \notag
		& 
		\leq C \Big(
			\mathbb{E} \int_{0}^{T}\int_{G_{0} \cap \mathcal{M}} s^{7} e^{2 s \phi} |w|^{2} d t
			+ \mathbb{E} \int_{Q} e^{2 s \phi} |f|^{2} d t
			+ \mathbb{E} \int_{Q} s^{4} e^{2 s \phi} |g|^{2} d t
		\\
		& \quad \quad ~ 
		+ h^{-4} \mathbb{E} \int_{\mathcal{M}} e^{2 s \phi}|w|^2 \Big|_{t=0}
		+ h^{-4} \mathbb{E} \int_{\mathcal{M}} e^{2 s \phi}|w|^2 \Big |_{t=T} \Big),
	\end{align}
	for all $ w \in L^{2}_{\mathbb{F}}(\Omega; C([0,T]; L^{2}_{h}(\mathcal{M}))) $  and $ f,  g \in L^{2}_{\mathbb{F}}(0,T;L^{2}_{h}(\mathcal{M})) $ satisfying 
	\begin{align}
		\label{eqW}
		- d w + D_{h}^{4} w d t = f d t + g d W(t) 
	\end{align}
	with  $w=0$ on $\partial \mathcal{M}$ and $D_h w=0$ on $\partial \mathcal{M}^*$.
\end{theorem}

The rest of this paper is organized as follows.  In \cref{sec2}, we give some preliminaries.  The Carleman estimate \cref{eqFinalCarleman}  is proved in \cref{sec3}.  In \cref{sec4}, we prove \cref{thmObservabilityEstiames}.  Finally, we prove the main  controllability result \cref{thmCon} based on the Hilbert uniqueness method in \cref{sec5}.

\section{Some preliminary results}
\label{sec2}

This section is devoted to presenting  some preliminary results we needed.

\begin{lemma}\cite[Lemma 2.1]{Cerpa2022}
	\label{lemmaLibF}
	For any $u, v \in C(\overline{\mathcal{W}})$, we have for the difference operator
	\begin{align}
		\label{eqDuv}
		D_h(u v)=D_h u A_h v+A_h u D_h v \text {, on } \mathcal{W}^* .
	\end{align}	
	Similarly, the average of the product gives
	\begin{align}
		\label{eqAuv}
		A_h(u v)=A_h u A_h v+\frac{h^2}{4} D_h u D_h v \text {, on } \mathcal{W}^* .
	\end{align}
	Finally, on $\mathring{\mathcal{W}} $ we have
	\begin{align}
		\label{eqA2u}
	u=A_h^2 u-\frac{h^2}{4} D_h^2 u.
	\end{align}
\end{lemma}

\begin{corollary}\cite[Corollary 2.2]{Cerpa2022}
	Let $\overline{\mathcal{W}} \subseteq \mathcal{M}$ be a regular mesh.
    For $u \in L(\overline{\mathcal{W}})$,
	\begin{align}
		\label{eqAu2}
		A_h  (u^2) = (A_h u) ^2+\frac{h^2}{4} (D_h u )^2, \text { on } \mathcal{W}^*.
	\end{align}
	In particular, for all $u \in L(\overline{\mathcal{W}})$,
	\begin{align}
		\label{eqAu2leq}
		A_h(u^2) \geq(A_h u)^2 \text {, on } \mathcal{W}^* .
	\end{align}
	For $u \in L(\overline{\mathcal{W}})$
	\begin{align}
		\label{eqDu2}
		D_h(u^2)=2 D_h u A_h u .
	\end{align}
\end{corollary}

To deal with the boundary conditions, we define the outward normal for $ x \in \partial \mathcal{W} $ as 
\begin{align}
	\label{eqNu}
	\nu(x) \deq \begin{cases}
		1, & \text {if }\tau_{-}(x) \in \mathcal{W}^{*} \text{ and } \tau_{+}(x) \notin \mathcal{W}^{*}, \\
		- 1, & \text {if } \tau_{-}(x) \notin \mathcal{W}^{*} \text{ and } \tau_{+}(x) \in \mathcal{W}^{*}, \\
		0, &  \text{ otherwise }.
	\end{cases}
\end{align}
We also introduce the trace operator for $ u \in L(\mathcal{W}^{*}) $ as 
\begin{align*}
	\forall \; x \in \partial \mathcal{W},  \quad  t_{r}(u) \deq 
	\begin{cases}
		\tau_{-} u(x), &  \text {if }  \nu(x)=1, \\
		\tau_{+} u(x), &  \text {if }  \nu(x)=-1, \\
		0, &   \text {if } \nu(x)=0.
	\end{cases}
\end{align*}
Then by the definition of  trace operator, for $ x \in \partial \mathcal{W} $ and $ u \in L(\overline{\mathcal{W}}) $, it is easy to check that
\begin{align}\label{tra}
t_r(|A_hu|^2)\nu-\frac{h^2}{4}t_r(|D_hu|^2)\nu=|u|^2\nu-hut_r(D_hu) .
\end{align}
Finally, let us  introduce the discrete integration on the boundary for $u \in L(\partial \mathcal{W})$ as
\begin{align*}
\int_{\partial \mathcal{W}} u \deq \sum_{x \in \partial \mathcal{W}} u(x) .
\end{align*}

\begin{proposition}\cite[Proposition 2.3]{Cerpa2022}
	Let $\mathcal{W}$ be a semi-discrete regular mesh. For $u \in L(\overline{\mathcal{W}})$ and $v \in L (\mathcal{W}^*)$ we have
	\begin{align}
		\label{eqDIgbp}
		\int_\mathcal{W} u D_h v=-\int_{\mathcal{W}^*} D_h u v+\int_{\partial \mathcal{W}} u t_{r}(v) \nu
	\end{align}
	and
	\begin{align}
		\label{eqAIgbp}
		\int_\mathcal{W} u A_h v=\int_{\mathcal{W}^*} A_h u v-\frac{h}{2} \int_{\partial \mathcal{W}} u t_{r}(v) .
	\end{align}
\end{proposition}

\begin{corollary}\cite[Corollary 2.4]{Cerpa2022}
	Let $\mathcal{W}$ be a semi-discrete regular mesh. For $u, v \in L(\overline{\mathcal{W}})$ we have
	\begin{align}
		\label{eqD2Ibp}
		& \int_\mathcal{W} u D_h^2 v=\int_{\overline{\mathcal{W}}} v D_h^2 u-\int_{\partial \mathcal{W}^*} D_h u t_{r}(v) \nu+\int_{\partial \mathcal{W}} u t_{r}\left(D_h v\right) \nu, \\
		\label{eqA2Ibp}
		& \int_\mathcal{W} u A_h^2 v=\int_{\overline{\mathcal{W}}} v A_h^2 u-\frac{h}{2} \int_{\partial \mathcal{W}^*} A_h u t_{r}(v)-\frac{h}{2} \int_{\partial \mathcal{W}} u t_{r}\left(A_h v\right),
	\end{align}
	and
	\begin{align}
		\label{eqADInt}
		\notag
		\int_\mathcal{W} u A_h D_h  v & =-\int_{\overline{\mathcal{W}}} v A_h D_h  u+\frac{h}{2} \int_{\partial \mathcal{W}^*} D_h u t_{r}(v)+\int_{\partial \mathcal{W}} u t_{r}\left(A_h v\right) \nu \\
		& =-\int_{\overline{\mathcal{W}}} v A_h D_h  u-\frac{h}{2} \int_{\partial \mathcal{W}} u t_{r}\left(D_h v\right)+\int_{\partial \mathcal{W}^*} A_h u t_{r}(v) \nu .
	\end{align}
\end{corollary}

\begin{proposition}\cite[Corollary 2.7]{Cerpa2022}
	\label{propADf}
	Let $ f $ be a $ (n+4) $-times differentiable function defined on $ \mathbb{R} $ and $ m, n \in \mathbb{N} $. Then 
	\begin{align*}
		A_h^m D_h^n f=f^{(n)}+R_{A_h^m} (f^{(n)} )+R_{D_h^n}(f)+R_{A_h^m D_h^n}(f),
	\end{align*}
	where 
	\begin{align*}
		\begin{gathered}
			R_{D_h^n}(f):=h^2 \sum_{k=0}^n\begin{pmatrix}
				n \\
				k
			\end{pmatrix}(-1)^k\bigg(\frac{(n-2 k)}{2}\bigg)^{n+2} \int_0^1 \frac{(1-\sigma)^{n+1}}{(n+1) !} f^{(n+2)}\bigg(\cdot+\frac{(n-2 k) h}{2} \sigma\bigg) d \sigma \\
			R_{A_h^n}(f):=\frac{h^2}{2^{n+2}} \sum_{k=0}^n\begin{pmatrix}
				n \\
				k
			\end{pmatrix}(n-2 k)^2 \int_0^1(1-\sigma) f^{(2)}\bigg(\cdot+\frac{(n-2 k) h}{2} \sigma\bigg) d \sigma
			,
			\end{gathered}
	\end{align*}
	and $ R_{A_h^m D_h^n}:=R_{A_h^m} \circ R_{D_h^n} $. Here, $ f^{(n)} $ denotes the $ n $-th derivative of $ f $ with respect to $ x $.
\end{proposition}

In the sequel, for any $n\in \mathbb{N}$, we denote by $\mathcal{O}(s^n)$ a function of
order $s^n$ for sufficiently large $s$ (which is independent of $\mu$ and $h$), by $\mathcal{O}_\mu(s^n)$ a function of order $s^n$
for fixed $\mu$ and sufficiently
large $s$ (which is independent of $h$), and
by $\mathcal{O}_\mu((sh)^n)$ a function of
order $(sh)^n$ for fixed $\mu$ and sufficiently
large $s$, sufficiently
small $h$.

\begin{lemma}\cite[Lemma 4.3]{Lecaros2021}
	For $ \alpha, \beta \in \mathbb{N} $, we have 
	\begin{align}
		\label{eqprpr} \notag
		\partial_{x}^{\beta} (r \partial_{x}^{\alpha} \rho) 
		& = 
		\alpha^{\beta} (- s \varphi)^{\alpha} \mu^{\alpha+\beta} (\partial_{x} \psi)^{\alpha + \beta}
		+ \alpha \beta(s \varphi)^{\alpha} \mu^{\alpha + \beta - 1} \mathcal{O}(1)
		+ s^{\alpha - 1} \alpha (\alpha - 1) \mathcal{O}_{\mu}(1)
		\\
		& = s^{\alpha} \mathcal{O}_{\mu}(1)
		.
	\end{align}
\end{lemma}

\begin{corollary}\cite[Corollary 4.4]{Lecaros2021}
	For $ \alpha, \beta, \gamma \in \mathbb{N} $, we have 
	\begin{align}
		\label{eqpr2pr} \notag
		\partial_{x}^{\gamma} (r^{2} (\partial_{x}^{\alpha} \rho) \partial_{x}^{\beta} \rho)
		& = 
		(\alpha + \beta)^{\gamma}(- s \varphi)^{\alpha+\beta} \mu^{\alpha+\beta+\gamma} (\partial_{x} \psi)^{\alpha + \beta + \gamma}
		+ \gamma (\alpha + \beta) (s \varphi)^{\alpha + \beta} \mu^{\alpha + \beta + \gamma -1} \mathcal{O}(1)
		\\ \notag
		& \quad 
		+ s^{\alpha + \beta - 1} [\alpha(\alpha-1)+ \beta(\beta-1)] \mathcal{O}_{\mu}(1)
		\\
		& =
		s^{\alpha+\beta}\mathcal{O}_{\mu}(1)
		.
	\end{align}
\end{corollary}
By  \cref{eqAlpha,eqprpr}, it easy to see that the following result holds.
\begin{lemma}
	For $ \alpha, \beta, \gamma, \sigma\in \mathbb{N} $, we have 
	\begin{align}
		\label{eqpr2pra} \notag
		\partial_{x}^ \sigma(r \partial_x^\alpha \rho \partial_{x}^\beta(r \partial_x^\gamma \rho))
		& = \gamma^\beta
		(\alpha + \gamma)^{\sigma}(- s \varphi)^{\alpha+\gamma} \mu^{\alpha+\beta+\gamma+\sigma} (\partial_{x} \psi)^{\alpha + \beta + \gamma+\sigma}
		+  (s \varphi)^{\alpha + \gamma} \mu^{\alpha + \beta + \gamma+\sigma -1} \mathcal{O}(1)
		\\ \notag
		& \quad 
		+ s^{\alpha + \gamma - 1}  \mathcal{O}_{\mu}(1)
		\\
		& =
		s^{\alpha+\gamma}\mathcal{O}_{\mu}(1)
		.
	\end{align}
\end{lemma}

\begin{lemma}\cite[Lemma 4.8]{Lecaros2021}
	Let $\alpha, j, k,  m, n \in \mathbb{N}$. Provided  $ \lambda h (\delta T^{2})^{-1}\leq 1$, we have
	\begin{align}
		\label{eqdprpr}  
		A_h^j D_h^k \partial_{x}^\alpha(r A_h^m D_h^n \rho)=\partial_x^k \partial_{x}^\alpha(r \partial_x^n \rho)+s^n \mathcal{O}_\mu((s h)^2)=s^n \mathcal{O}_\mu(1)
		.
	\end{align}
\end{lemma}

\begin{theorem}\cite[Theorem 4.9]{Lecaros2021}
	Let $\alpha, \beta, j, k, l, m, n, p \in \mathbb{N}$. Provided  $ \lambda h (\delta T^{2})^{-1}\leq 1$, we have
	\begin{align}
		\label{eqdpr2pr} \notag
		A_h^p D_h^l \partial_{x}^\beta(r^2 A_h^j D_h^k(\partial_{x}^\alpha \rho) A_h^m D_h^n(\rho)) & =\partial_x^l \partial_{x}^\beta(r^2 \partial_x^k \partial_{x}^\alpha \rho \partial_x^n \rho)+s^{n+k+\alpha} \mathcal{O}_\mu((s h)^2) \\
		& =s^{n+k+\alpha} \mathcal{O}_\mu(1) .
	\end{align}
\end{theorem}

\begin{theorem}\cite[Theorem 2.9]{Cerpa2022}
	Let $\alpha, j, k,  m, n \in \mathbb{N}$. Provided  $ \lambda h (\delta T^{2})^{-1}\leq 1$, we have
	\begin{align}
		\label{eqtCareleman} 
		\partial_t A_{h}^{j} D_{h}^{k}  (r A_{h}^{m} D_{h}^{n} (\partial_x^\alpha \rho) )=T \theta s^{\alpha+n} \mathcal{O}_{\mu}(1) .
	\end{align}
\end{theorem}
By \cref{propADf} and \cref{eqdprpr}, it is easy to check that the following estimate holds.
\begin{lemma}
	Let $i, j, k, l, m, n, p, q \in \mathbb{N}$. Provided  $ \lambda h (\delta T^{2})^{-1}\leq 1$, we have
\begin{align}
	\label{eqdprpra} \notag
	A_h^i D_h^j (r A_h^k D_h^l \rho A_h^m D_h^n(r A_{h}^{p} D_{h}^{q}\rho)) & = \partial_{x}^j(r \partial_x^l \rho \partial_{x}^n(r \partial_x^q \rho))+s^{l+q} \mathcal{O}_\mu((s h)^2) \\
	& =s^{l+q} \mathcal{O}_\mu(1).
\end{align}
\end{lemma}

\section{Proof of the Carleman estimate}
\label{sec3}

\begin{proof}[Proof of \cref{thmCarleman}]

\emph{Step 1.}
At first, we define 
\begin{align}
	\label{eqC1}
	\mathcal{P}_{h} w   = - d w + D_{h}^{4} w
	.
\end{align}
Letting $ v = r w $ and recalling that $ \rho = r^{-1} $, utilizing It\^o's formula, we have 
\begin{align}
	\label{eqC2}
	-r d w = - r d (\rho v) = - d v +\lambda \phi \partial_{t} \theta v d t 
	.
\end{align}
Thanks to \cref{lemmaLibF} and $w=\rho v$, it holds that 
\begin{align}
	\label{eqC3}
	D_{h}^{4} w  
	& =
	D_h^4 \rho A_h^4v
	+4 A_h D_h^3 \rho A_h^3 D_h v
	+6 A_h^2 D_h^2 \rho A_h^2 D_h^2 v 
	+4 A_h^3 D_h \rho A_h D_h^3 v
	+ A_h^4 \rho D_h^4 v .
\end{align}
Combining \cref{eqC1,eqC2,eqC3}, we deduce that
\begin{align}
	\label{eqC4}
	r\mathcal{P}_{h} w  = I_{1} d t + I_{2} + I_{3} d t,
\end{align}
where 
\begin{align}
	\notag
	\label{eqI}
	& I_{1} = 
		6 r A_{h}^{2} D_{h}^{2}  \rho A_{h}^{2} D_{h}^{2}  v 
		+ r D_{h}^{4} \rho A_{h}^{4} v 
		+ r A_{h}^{4} \rho D_{h}^{4} v 
		+ 6 A_{h} D_{h}(r A_{h}^{2}D_{h}^{2} \rho) A_{h} D_{h} v 
		+ \Phi(v)
		,
	\\[2mm] \notag
	& I_{2} =
		- d v 
		+ 4 r A_{h} D_{h}^{3}  \rho  A_{h}^{3} D_{h} v d t
		+ 4 r A_{h}^{3} D_{h}  \rho  A_{h} D_{h}^{3}  v d t 
		+ 2 A_{h} D_{h} (r A_{h} D_{h}^{3} \rho) v d t 
		,
	\\[2mm] \notag
	& I_{3} = 
		 \lambda \phi \partial_{t} \theta  v   
		- 6 A_{h} D_{h}(r A_{h}^{2}D_{h}^{2} \rho) A_{h} D_{h} v 
		- 2 A_{h} D_{h} (r A_{h} D_{h}^{3} \rho) v  
		- \Phi(v)  
	\\
	& \begin{aligned}
		\Phi(v) & = 
		3 h^{2} D_{h}^{2} [A_{h}  D_{h} (r A_{h}^{2} D_{h}^{2} \rho) A_{h} D_{h} v ] 
		+ \frac{h^{2}}{2} A_{h}^{2} [ A_{h} D_{h}  (r D_{h}^{4} \rho) A_{h}   D_{h} v ]
		+ \frac{h^{2}}{4} D_{h}^{2} (r D_{h}^{4} \rho) A_{h}^{2} v 
		\\
		& \quad 
		- D_{h}^{2} (r A_{h}^{4} \rho) D_{h}^{2} v 
		+ 2  A_{h} D_{h}  [ A_{h} D_{h}  (r A_{h}^{4} \rho) D_{h}^{2} v ] 
		.
	\end{aligned}
\end{align}
 Using \cref{eqC4},  it follows that
\begin{align*}
	2  I_{1}r \mathcal{P}_{h} w= 2 I_{1}^{2} d t + 2 I_{1} I_{2} +2 I_{1} I_{3} d t, 
\end{align*}
which, combined with Cauchy-Schwarz inequality, \cref{eqC1,eqW}, implies
\begin{align}
	\label{eqC5}
	\mathbb{E} \int_{Q} |r f|^{2} d t 
	\geq 
	2 \mathbb{E} \int_{Q} I_{1} I_{2} 
	- \mathbb{E} \int_{Q} |I_{3}|^{2} d t .
\end{align}
The next steps provide an estimate for the right hand side of \cref{eqC5}.  For convenience,
we denote 
\begin{align}
	\label{eqDefIij}
	\mathbb{E} \int_{Q} I_{1} I_{2} = \sum_{i=1}^{5} \sum_{j=1}^{4} I_{i j},
\end{align}
where $ I_{i j} $ stands for the inner product in $ L^{2}_{\mathbb{F}}(0,T;L_{h}^{2}(\mathcal{M})) $ between the $ i $-th term of $ I_{1} $ and $ j $-th term of $ I_{2} $.

In order to shorten the formulas used in the sequel, we define the following lower order terms 
\begin{align*}
	\mathcal{A}_{1} &  =
	\mathbb{E} \int_{Q} [s^{6} \mathcal{O}_{\mu}(1) + s^{7} \mathcal{O}_{\mu}((sh)^{2}) +  s^{7} \mu^{7} \varphi^{7} \mathcal{O}(1)] |v|^{2} d t
	\\
	& \quad 
	+ \mathbb{E} \int_{Q^{*}} [s^{4} \mathcal{O}_{\mu}(1) + s^{5} \mathcal{O}_{\mu}((sh)^{2}) + s^{5} \mu^{5} \varphi^{5} \mathcal{O}(1)] |D_{h} v|^{2} d t
	\\
	& \quad 
	+ \mathbb{E} \int_{\overline{Q}} [s^{2} \mathcal{O}_{\mu}(1) + s^{3} \mathcal{O}_{\mu}((sh)^{2}) + s^{3} \mu^{3} \varphi^{3} \mathcal{O}(1) ] |D^{2}_{h} v|^{2} d t
	\\
	& \quad 
	+ \mathbb{E} \int_{Q^{*}} [  \mathcal{O}_{\mu}(1) + s  \mathcal{O}_{\mu}((sh)^{2})+ s \mu \varphi \mathcal{O}(1)] |D_{h}^{3} v |^{2} d t
	,
	\\
	\mathcal{A}_{2} & = 
	\mathbb{E} \int_{Q} s^{4} \mathcal{O}_{\mu}(1) |d v|^{2}
	+ \mathbb{E} \int_{\overline{Q}} h^{2} \mathcal{O}_{\mu}(1) |D_{h}^{2} (d v)|^{2}
	,
	\\
	\mathcal{A}_{3} & = 
	\mathbb{E} \int_{\mathcal{M}} s^{4} \mathcal{O}_{\mu}(1) |v|^{2} \Big|_{t=0}
	+ \mathbb{E} \int_{\mathcal{M}^{*}} s^{2} \mathcal{O}_{\mu}(1) |D_{h} v|^{2} \Big|_{t=0}
	+ \mathbb{E} \int_{\overline{\mathcal{M}}} \mathcal{O}_{\mu}(1) |D_{h}^{2} v|^{2} \Big|_{t=0}
	\\
	& \quad 
	+ \mathbb{E} \int_{\mathcal{M}} s^{4} \mathcal{O}_{\mu}(1) |v|^{2} \Big|_{t=T}
	+ \mathbb{E} \int_{\mathcal{M}^{*}} s^{2} \mathcal{O}_{\mu}(1) |D_{h} v|^{2} \Big|_{t=T}
	+ \mathbb{E} \int_{\overline{\mathcal{M}}} \mathcal{O}_{\mu}(1) |D_{h}^{2} v|^{2} \Big|_{t=T}
	,
	\\
	\mathcal{B} & = 
	\mathbb{E} \int_{\partial Q} [s^{4} \mathcal{O}_{\mu}(1) + s^{5} \mathcal{O}_{\mu}((sh)^{2}) ] t_{r} (|D_{h} v|^{2}) d t
	+ \mathbb{E} \int_{\partial Q} [s^{2} \mathcal{O}_{\mu}(1) + s^{3} \mathcal{O}_{\mu}(sh) ]  |D^{2}_{h} v|^{2} d t
	\\
	& \quad 
	+ \mathbb{E} \int_{\partial Q} [\mathcal{O}_{\mu}(1) + s \mathcal{O}_{\mu}(sh) ] t_{r}( |D^{3}_{h} v|^{2}) d t
	.
\end{align*}

\emph{Step 2.}
Let us compute $ I_{ i 1 } $ for $ i =1, \cdots ,5 $.

Since $w=0$ on $\partial \mathcal{M}$ and $D_h w=0$ on $\partial \mathcal{M}^*$, noting that $ v= r w $, we obtain 
\begin{align}
	\label{eqvBoundary}
	v = 0 \text{ on } \partial Q, \ 
	D_{h} v = A_{h} v = 0 \text{ on } \partial Q^{*}
	,
\end{align}
and
\begin{align}
	\label{traBoundary}
	t_r(v) = 0 \text{ on } \partial Q^* \text{ and } \partial \overline{Q}^*, \
		t_r(D_{h} v) =	t_r( A_{h} v) = 0 \text{ on } \partial \overline{Q}
	.
\end{align}
From \cref{eqI,eqADInt,eqvBoundary,eqDuv,eqAuv,eqA2u}, letting $ q_{11} = r A_{h}^{2} D_{h}^{2}  \rho $, we have 
\begin{align}
	\notag
	\label{eqI11e1}
	I_{11} 
	& = - 6 \mathbb{E} \int_{Q} q_{11} A_{h}^{2} D_{h}^{2}  v d v 
	=
	6 \mathbb{E} \int_{\overline{Q}} A_{h} D_{h} v A_{h} D_{h} (q_{11} d v)
	\\\notag
	& = 
	6 \mathbb{E} \int_{\overline{Q}} A_{h} D_{h} q_{11} A_{h}^{2} (d v) A_{h} D_{h} v 
	+ 6 \mathbb{E} \int_{\overline{Q}} A_{h}^{2}   q_{11} A_{h} D_{h} (d v) A_{h} D_{h} v 
	\\\notag
	& \quad 
	+ \frac{3}{2} h^{2} \mathbb{E} \int_{\overline{Q}} D_{h}^{2}   q_{11} A_{h} D_{h} (d v) A_{h} D_{h} v 
	+ \frac{3}{2} h^{2} \mathbb{E} \int_{\overline{Q}} A_{h} D_{h}    q_{11}  D_{h}^{2} (d v) A_{h} D_{h} v 
	\\\notag
	& = 
	6 \mathbb{E} \int_{Q} A_{h} D_{h} q_{11}  A_{h} D_{h} v d v
	+ 6 \mathbb{E} \int_{\overline{Q}}   \bigg( q_{11} + \frac{h^{2}}{2} D_{h}^{2} q_{11} \bigg) A_{h} D_{h} v  A_{h} D_{h}  (dv) 
	\\
	& \quad 
	+ 3 h^{2} \mathbb{E} \int_{\overline{Q}} A_{h} D_{h}    q_{11}  A_{h} D_{h} v  D_{h}^{2} (d v). 
\end{align}
By It\^o's formula, it follows that
\begin{align}
	\notag
	\label{eqI11e2}
	& 6 \mathbb{E} \int_{\overline{Q}}   \bigg( q_{11} + \frac{h^{2}}{2} D_{h}^{2} q_{11} \bigg) A_{h} D_{h} v  A_{h} D_{h}  (dv) 
	\\\notag
	&
	=
	3 \mathbb{E} \int_{\overline{M}} \bigg[  \bigg( q_{11} + \frac{h^{2}}{2} D_{h}^{2} q_{11} \bigg) | A_{h} D_{h} v |^{2} \bigg] \Bigg|_{0}^{T}
	- 3 \mathbb{E} \int_{\overline{Q}} \partial_{t} \bigg( q_{11} + \frac{h^{2}}{2} D_{h}^{2} q_{11} \bigg) | A_{h} D_{h} v |^{2}  d t
	\\
	& \quad 
	- 3 \mathbb{E} \int_{\overline{Q}} \bigg( q_{11} + \frac{h^{2}}{2} D_{h}^{2} q_{11} \bigg) | A_{h} D_{h} (dv) |^{2}. 
\end{align}
By \cref{eqAIgbp,eqAu2leq,eqdprpr,eqvBoundary},   for $ \lambda h (\delta T^{2})^{-1}\leq 1$, we conclude that 
	\begin{align}
		\notag
		\label{eqI11e2a}
		&3\mathbb{E} \int_{\overline{M}}  \bigg( q_{11} + \frac{h^{2}}{2} D_{h}^{2} q_{11} \bigg) | A_{h} D_{h} v |^{2}
		 \geq 
	-\mathbb{E} \int_{\overline{M}} [s^{2} \mathcal{O}_{\mu}(1) + \mathcal{O}_{\mu}((sh)^{2})]|A_{h}D_{h} v|^{2} 
	\\\notag
	& \geq 
	- \mathbb{E} \int_{\overline{M}} [ s^{2} \mathcal{O}_{\mu}(1) + \mathcal{O}_{\mu}((sh)^{2}) ] A_{h}  (|D_{h}v|^2)\\
	& \geq 
	- \mathbb{E} \int_{M^{*}} [ s^{2} \mathcal{O}_{\mu}(1) + \mathcal{O}_{\mu}((sh)^{2}) ] |D_{h}v|^2
	.
\end{align}
Therefore, similar to \cref{eqI11e2a}, by \cref{eqtCareleman}, it holds that
	\begin{align}
	\notag
	\label{eqI11e2c}
	&- 3 \mathbb{E} \int_{\overline{Q}} \partial_{t} \bigg( q_{11} + \frac{h^{2}}{2} D_{h}^{2} q_{11} \bigg) | A_{h} D_{h} v |^{2}  d t
	- 3 \mathbb{E} \int_{\overline{Q}} \bigg( q_{11} + \frac{h^{2}}{2} D_{h}^{2} q_{11} \bigg) | A_{h} D_{h} (dv) |^{2}
		\\
	& \geq  
	-\mathbb{E} \int_{Q^*} T\theta[  s^{2} \mathcal{O}_{\mu}(1) + \mathcal{O}_{\mu}((sh)^{2})]|D_{h} v|^{2}  d t
	- \mathbb{E} \int_{Q^*} [ s^{2} \mathcal{O}_{\mu}(1) + \mathcal{O}_{\mu}((sh)^{2}) ] |D_{h}(dv)|^2
	.
\end{align}
By Cauchy-Schwarz inequality and combining \cref{eqI11e2,eqI11e2a,eqI11e2c}, for $ \lambda h (\delta T^{2})^{-1}\leq 1$ and $ \lambda \geq T + T^{2} $, we obtain 
	\begin{align}
	\notag
	\label{eqI11e2b}
	& 6 \mathbb{E} \int_{\overline{Q}}   \bigg( q_{11} + \frac{h^{2}}{2} D_{h}^{2} q_{11} \bigg) A_{h} D_{h} v  d (A_{h} D_{h}  v) 
\\\notag
	& \geq 
- \mathbb{E} \int_{Q^{*}} [ s^{2} \mathcal{O}_{\mu}(1) + \mathcal{O}_{\mu}((sh)^{2}) ] |D_{h} (d v)|^{2}-  \mathcal{A}_{1} - \mathcal{A}_{3}
\\\notag
& = 
\mathbb{E} \int_{Q} [ s^{2} \mathcal{O}_{\mu}(1) + \mathcal{O}_{\mu}((sh)^{2}) ]  D_{h}^{2} (d v) d v 
-  \mathcal{A}_{1} - \mathcal{A}_{3}
\\
& \geq 
-\frac{1}{4} \mathbb{E} \int_{Q} |D_{h}^{2} (d v)|^{2} 
-  \mathcal{A}_{1} - \mathcal{A}_{2} - \mathcal{A}_{3}
.
\end{align}
Combining \cref{eqI11e1,eqI11e2b,eqD2Ibp,traBoundary,eqvBoundary}, for $ \lambda h (\delta T^{2})^{-1}\leq 1$ and $ \lambda \geq T + T^{2} $, we have 
\begin{align}
	\notag
	\label{eqI11}
	I_{11} 
	&\geq 
	6 \mathbb{E} \int_{Q} A_{h} D_{h} (r A_{h}^{2} D_{h}^{2} \rho)  A_{h} D_{h} v d v
	+ 3 h^{2} \mathbb{E} \int_{\overline{Q}} D_{h}^{2} [A_{h} D_{h}    (r A_{h}^{2} D_{h}^{2} \rho)  A_{h} D_{h} v]   d v 
	\\
	& \quad 
	-\frac{1}{4} \mathbb{E} \int_{Q} |D_{h}^{2} (d v)|^{2} 
	-  \mathcal{A}_{1} - \mathcal{A}_{2} - \mathcal{A}_{3}
	.
\end{align}

From \cref{eqA2Ibp,eqvBoundary,eqAuv,eqA2u}, letting $ q_{21} = r D_{h}^{4} \rho$, we obtain 
\begin{align}
	\notag
	\label{eqI21e1}
	I_{21} 
	& =
	- \mathbb{E} \int_{Q} q_{21} A_{h}^{4} v d v  
	=
	- \mathbb{E} \int_{\overline{Q}} A_{h}^{2} (q_{21} d v ) A_{h}^{2} v 
	\\\notag
	& = 
	- \mathbb{E} \int_{\overline{Q}} A_{h}^{2} q_{21}   A_{h}^{2} (d v )A_{h}^{2} v 
	- \frac{h^{2}}{2} \mathbb{E} \int_{\overline{Q}} A_{h} D_{h} q_{21}    A_{h} D_{h} (d v ) A_{h}^{2} v
	- \frac{h^{2}}{4} \mathbb{E} \int_{\overline{Q}} D_{h}^{2} q_{21}    A_{h}^{2} (d v ) A_{h}^{2} v
	\\
	& \quad 
	+ \frac{h^{2}}{4} \mathbb{E} \int_{\overline{Q}} D_{h}^{2} q_{21}      A_{h}^{2} v dv
	.
\end{align}
Thanks to It\^o's formula and similar to \cref{eqI11e2a}, from \cref{eqdprpr,eqtCareleman}, for $ \lambda h (\delta T^{2})^{-1}\leq 1$ and $ \lambda \geq T + T^{2} $, we have 
\begin{align}
	\notag
	\label{eqI21e2}
	& 
	- \mathbb{E} \int_{\overline{Q}} A_{h}^{2} q_{21}   A_{h}^{2} (d v )A_{h}^{2} v 
	\\ \notag
	& = 
	- \frac{1}{2} \mathbb{E} \int_{\overline{\mathcal{M}}} ( A_{h}^{2} q_{21} |A_{h}^{2} v|^{2} ) \Big |_{0}^{T}
	+ \frac{1}{2} \mathbb{E} \int_{\overline{Q}} \partial_{t}(A_{h}^{2} q_{21}) |A_{h}^{2} v|^{2} d t 
	+ \frac{1}{2} \mathbb{E} \int_{\overline{Q}} A_{h}^{2} q_{21} |A_{h}^{2}(d v)|^{2}
	\\
	& \geq 
	-  \mathcal{A}_{1} - \mathcal{A}_{2} - \mathcal{A}_{3}
	.
\end{align}
Using the identity $ A_{h} D_{h} (d v ) A_{h}^{2} v=  \frac{1}{2} d [D_{h} (|A_{h} v|^{2})]-A_{h} D_{h} v  A_{h}^{2} (dv)-A_{h} D_{h} (d v) A_{h}^{2}(dv)$ and Cauchy-Schwarz inequality, from \cref{eqAu2leq,traBoundary,eqDIgbp,eqAIgbp,eqdprpr,eqtCareleman}, for $ \lambda h (\delta T^{2})^{-1}\leq 1$, we deduce that 
\begin{align}
	\notag
	\label{eqI21e3}
	& - \frac{h^{2}}{2} \mathbb{E} \int_{\overline{Q}} A_{h} D_{h} q_{21}    A_{h} D_{h} (d v ) A_{h}^{2} v
	\\ \notag
	& =
	- \frac{h^{2}}{4} \mathbb{E} \int_{\overline{Q}} d [ A_{h} D_{h} q_{21} D_{h} (|A_{h} v|^{2}) ]
		+\frac{h^{2}}{4} \mathbb{E} \int_{\overline{Q}} \partial_{t} (A_{h} D_{h} q_{21}) A_{h} D_{h} v A_{h}^{2} v d t 
	\\ \notag
	& \quad 
		+ \frac{h^{2}}{2} \mathbb{E} \int_{\overline{Q}} A_{h} D_{h} q_{21} A_{h} D_{h} v A_{h}^{2} (d v ) 
	+ \frac{h^{2}}{2} \mathbb{E} \int_{\overline{Q}} A_{h} D_{h} q_{21} A_{h} D_{h} (d v) A_{h}^{2} (d v )
	\\ \notag
	& =	- \frac{h^{2}}{4}	\mathbb{E} \int_{\overline{\mathcal{M}}}  A_{h} D_{h} q_{21} D_{h} (|A_{h} v|^{2})  \Big|_{0}^{T}
		+ \mathbb{E} \int_{\overline{Q}} T \theta s^{2} \mathcal{O}_{\mu}((sh)^{2}) A_{h} D_{h} v A_{h}^{2} v d t 
	\\ \notag
	& \quad 
	+ \frac{h^{2}}{2} \mathbb{E} \int_{\overline{Q}} A_{h} D_{h} q_{21} A_{h} D_{h} v A_{h}^{2} (d v )
 -\frac{h^{2}}{4} \mathbb{E} \int_{\overline{Q}} A_{h} D_{h}^2 q_{21} | A_{h} (d v )|^2
	\\ \notag
	& \geq
	 \frac{h^{2}}{2} \mathbb{E} \int_{\overline{Q}} A_{h} D_{h} q_{21} A_{h} D_{h} v A_{h}^{2} (d v ) 
	- \mathcal{A}_{1} - \mathcal{A}_{2} - \mathcal{A}_{3}
	\\
	& =
	 \frac{h^{2}}{2} \mathbb{E} \int_{\overline{Q}} A_{h}^{2} (A_{h} D_{h} q_{21} A_{h} D_{h} v)  d v 
	- \mathcal{A}_{1} - \mathcal{A}_{2} - \mathcal{A}_{3}.
\end{align}
Thanks to It\^o's formula, from \cref{eqAu2leq,eqdprpr,eqAIgbp,eqtCareleman}, for $ \lambda h (\delta T^{2})^{-1}\leq 1$ and $ \lambda \geq T + T^{2} $, it holds that 
\begin{align}
	\notag
	\label{eqI21e4}
	& 
	- \frac{h^{2}}{4} \mathbb{E} \int_{\overline{Q}} D_{h}^{2} q_{21}    A_{h}^{2} (d v ) A_{h}^{2} v
	\\ \notag
	& = 
	-\frac{h^{2}}{8} \mathbb{E} \int_{\overline{\mathcal{M}}} D_{h}^{2} q_{21} |A_{h}^{2} v|^{2}\Big|_{0}^{T}
	+ \frac{h^{2}}{8} \mathbb{E} \int_{\overline{Q}} \partial_{t} (D_{h}^{2} q_{21} ) |A_{h}^{2} v|^{2} d t 
	+ \frac{h^{2}}{8} \mathbb{E} \int_{\overline{Q}}    D_{h}^{2} q_{21}   |A_{h}^{2} (d v)|^{2} 
	\\
	& \geq 
	- \mathcal{A}_{1} - \mathcal{A}_{2} - \mathcal{A}_{3}
	.
\end{align}
Combining \cref{eqI21e1,eqI21e2,eqI21e3,eqI21e4}, for $ \lambda h (\delta T^{2})^{-1}\leq 1$ and $ \lambda \geq T + T^{2} $, it follows that 
\begin{align}
	\label{eqI21}
	I_{21} 
	\geq 
	 \frac{h^{2}}{2} \mathbb{E} \int_{\overline{Q}} A_{h}^{2} [A_{h} D_{h} (r D_{h}^{4} \rho) A_{h} D_{h} v]  d v 
	+ \frac{h^{2}}{4} \mathbb{E} \int_{Q} D_{h}^{2} (r D_{h}^{4} \rho)   A_{h}^{2} v d v 
	- \mathcal{A}_{1} - \mathcal{A}_{2} - \mathcal{A}_{3}.
\end{align}

From \cref{eqD2Ibp,eqvBoundary,eqDuv,eqA2u,eqADInt,traBoundary}, letting $ q_{31} = r A_{h}^{4} \rho  $, we have 
\begin{align}
	\notag
	\label{eqI31e1}
	I_{31} 
	& = 
	- \mathbb{E} \int_{Q} q_{31} D_{h}^{4} v d v 
	=
	- \mathbb{E} \int_{\overline{Q}} D_{h}^{2} (q_{31} d v) D_{h}^{2} v 
	\\ \notag
	& =
	- \mathbb{E} \int_{\overline{Q}} D_{h}^{2} q_{31} A_{h}^{2} (d v) D_{h}^{2} v 
	- 2 \mathbb{E} \int_{\overline{Q}} A_{h} D_{h}  q_{31} A_{h}D_{h} (d v) D_{h}^{2} v 
	- \mathbb{E} \int_{\overline{Q}} A_{h}^{2} q_{31} D_{h}^{2} (d v) D_{h}^{2} v 
	\\ \notag
	& = 
	- \mathbb{E} \int_{\overline{Q}} D_{h}^{2} q_{31}  D_{h}^{2} v  d v
	+ 2 \mathbb{E} \int_{\overline{Q}} A_{h} D_{h}  ( A_{h} D_{h}  q_{31} D_{h}^{2} v ) d v
	\\
	& \quad 
	- \mathbb{E} \int_{\overline{Q}} \bigg(A_{h}^{2} q_{31} + \frac{h^{2}}{4} D_{h}^{2} q_{31} \bigg) D_{h}^{2} (d v) D_{h}^{2} v 
	.
\end{align}
Thanks to It\^o's formula, from \cref{eqdprpr,eqtCareleman}, for $ \lambda h (\delta T^{2})^{-1}\leq 1$, $ h\leq 1 $  and $ \lambda \geq T + T^{2} $, it holds that 
\begin{align}
	\notag
	\label{eqI31e2}
	& - \mathbb{E} \int_{\overline{Q}} \bigg(A_{h}^{2} q_{31} + \frac{h^{2}}{4} D_{h}^{2} q_{31} \bigg) D_{h}^{2} (d v) D_{h}^{2} v
	\\ \notag
	& =
	- \frac{1}{2} \mathbb{E} \int_{\overline{Q}}  \bigg(A_{h}^{2} q_{31} + \frac{h^{2}}{4} D_{h}^{2} q_{31}\bigg) |D_{h}^{2} v |^{2}\bigg|_{0}^{T}
	+ \frac{1}{2} \mathbb{E} \int_{\overline{Q}} \partial_{t}\bigg(A_{h}^{2} q_{31} + \frac{h^{2}}{4} D_{h}^{2} q_{31} \bigg) |D_{h}^{2}  v|^{2} d t 
	\\ \notag
	& \quad 
	+ \frac{1}{2} \mathbb{E} \int_{\overline{Q}}  \bigg(A_{h}^{2} q_{31} + \frac{h^{2}}{4} D_{h}^{2} q_{31} \bigg) |D_{h}^{2}  (d v)|^{2} 
	\\
	& \geq 
	 \frac{1}{2}\mathbb{E} \int_{\overline{Q}}  |D_{h}^{2} (d v)|^{2}  
	- \mathcal{A}_{1} -  \mathcal{A}_{2} - \mathcal{A}_{3} 
	.
\end{align}
Combining \cref{eqI31e2,eqI31e1}, for $ \lambda h (\delta T^{2})^{-1}\leq 1$, $ h\leq 1 $  and $ \lambda \geq T + T^{2} $, we have 
\begin{align}
	\label{eqI31}
	\notag
	&I_{31}
	\geq
	- \mathbb{E} \int_{\overline{Q}} D_{h}^{2} (r A_{h}^{4} \rho)  D_{h}^{2} v  d v
	+ 2 \mathbb{E} \int_{\overline{Q}} A_{h} D_{h}  [ A_{h} D_{h}  (r A_{h}^{4} \rho) D_{h}^{2} v ]  d v
		\\ 
	& \quad \quad \quad 
	+ \frac{1}{2}\mathbb{E} \int_{\overline{Q}}  |D_{h}^{2} (d v)|^{2}  
	- \mathcal{A}_{1} -  \mathcal{A}_{2} - \mathcal{A}_{3} 
	.
\end{align}

Hence, combining \cref{eqI11,eqI21,eqI31,eqI}, for $ \lambda h (\delta T^{2})^{-1}\leq 1$, $ h\leq 1 $  and $ \lambda \geq T + T^{2} $, we obtain 
\begin{align}
	\label{eqIi1e1}
	I_{11} + I_{21} + I_{31} + I_{41} + I_{51} 
	\geq 
	\frac{1}{4} \mathbb{E} \int_{\overline{Q}}  |D_{h}^{2} (d v)|^{2}  
	- \mathcal{A}_{1} -  \mathcal{A}_{2} - \mathcal{A}_{3} 
	.
\end{align}

Noting that \cref{eqW} and $ v = r w $, by \cref{propADf,eqAlpha,eqDuv}, for $ \lambda h (\delta T^{2})^{-1}\leq 1$, we have 
\begin{align}
	\notag
	\label{eqIi1e2}
	\mathbb{E} \int_{\overline{Q}}  |D_{h}^{2} (d v)|^{2}  
	& = 
	\mathbb{E} \int_{\overline{Q}} |D_{h}^{2} (r g)|^{2}  d t
	\\ \notag
	& \geq 
	\frac{1}{3}  \mathbb{E} \int_{\overline{Q}} r^{2} |D_{h}^{2} g |^{2} d t 
	- \mathbb{E} \int_{\overline{Q}} \mathcal{O}_{\mu}((sh)^{4}) r^{2} |D_{h}^{2} g|^{2} d t
	- \mathbb{E} \int_{Q^{*}} s^{2} \mathcal{O}_{\mu}(1) r^{2} |D_{h} g|^{2} d t
	\\
	& \quad 
	- \mathbb{E} \int_{Q} s^{4} \mathcal{O}_{\mu}(1) r^{2} | g|^{2} d t
	.
\end{align}
From \cref{eqDIgbp,eqDuv,propADf,eqAlpha}, thanks to Cauchy-Schwarz inequality, for $ \lambda h (\delta T^{2})^{-1}\leq 1$, we obtain 
\begin{align}
	\label{eqIi1e3}
	-\mathbb{E} \int_{Q^{*}} s^{2} \mathcal{O}_{\mu}(1) r^{2} |D_{h} g|^{2} d t
	\geq 
	-\frac{1}{4} \mathbb{E} \int_{\overline{Q}} r^{2} |D_{h}^{2} g |^{2} d t 
	- \mathbb{E} \int_{Q} s^{4} \mathcal{O}_{\mu}(1) r^{2} | g|^{2}  d t 
	.
\end{align} 
Hence, combining  \cref{eqIi1e1,eqIi1e2,eqIi1e3}, for all $ \mu \geq 1 $, one can find  constant $ h_{1} \leq 1 $, such that for all $ \lambda \geq T + T^{2} $, $ h\leq h_{1} $ and $ \lambda h (\delta T^{2})^{-1}\leq 1 $, there exists a constant $ C > 0 $ such that 
\begin{align}
	\label{eqIi1}
	I_{11} + I_{21} + I_{31} + I_{41} + I_{51} 
	\geq 
	\frac{1}{48} \mathbb{E} \int_{\overline{Q}} r^{2} |D_{h}^{2} g |^{2} d t 
	- C(\mu) \mathbb{E} \int_{Q} s^{4} r^{2} |g|^{2} d t 
	- \mathcal{A}_{1} - \mathcal{A}_{3}
	.
\end{align}

\emph{Step 3.} In this step, we compute $ I_{ij} $ for $ i = 1, \cdots, 5 $ and $ j = 2, 3, 4 $.

Noting that \cref{eqDu2,eqDIgbp,eqA2u,eqdpr2pr,eqpr2pr,eqPsi,eqNu},  letting $ q_{12} = r^{2} A_{h}^{2} D_{h}^{2} \rho A_{h} D_{h}^{3} \rho $, for $ \lambda h (\delta T^{2})^{-1}\leq 1$, we have 
\begin{align}
	\notag
	\label{eqI12}
	I_{12}
	& = 
	24 \mathbb{E} \int_{Q} q_{12} A_{h}^{2} D_{h}^{2} v A_{h}^{3} D_{h} v  d t 
	= 
	12 \mathbb{E} \int_{Q} q_{12} D_{h} [ (A_{h}^{2} D_{h}  v )^{2} ]  d t 
	\\ \notag
	& =
	- 12 \mathbb{E} \int_{Q^{*}} D_{h} q_{12} |A_{h}^{2} D_{h} v|^{2}  d t 
	+ 12 \mathbb{E} \int_{\partial Q} q_{12} t_{r} ( |A_{h}^{2} D_{h} v|^{2}  ) \nu d t 
	\\ \notag
	& =
	- 12 \mathbb{E} \int_{Q^{*}} D_{h} q_{12} \bigg|D_{h} v + \frac{h^{2}}{4} D_{h}^{3} v \bigg|^{2}  d t 
	+ 12 \mathbb{E} \int_{\partial Q} q_{12} t_{r} \bigg( \bigg|D_{h} v + \frac{h^{2}}{4} D_{h}^{3} v \bigg|^{2}  \bigg) \nu d t 
	\\ \notag
	& \geq 
	60 \mathbb{E} \int_{Q^{*}} s^{5} \mu^{6} \varphi^{5} |\partial_{x} \psi|^{6} |D_{h} v |^{2} d t 
	- \mathbb{E} \int_{Q^{*}} [s  \mathcal{O}_{\mu}((sh)^{2}) + s  \mathcal{O}_{\mu}((sh)^{4})] |D_{h}^{3} v|^{2} d t
	\\ \notag
	& \quad 
	- \mathbb{E} \int_{Q^{*}} [s^{4} \mathcal{O}_{\mu}(1) + s^{5} \mathcal{O}_{\mu}((sh)^{2}) + s^{5} \mu^{5} \varphi^{5} \mathcal{O}(1)] |D_{h} v|^{2} d t
	\\ \notag
	& \quad 
	- 12 \mathbb{E} \int_{\partial Q} s^{5} \mu^{5} \varphi^{5} (\partial_{x} \psi)^{5} t_{r} (|D_{h} v|^{2}) \nu d t 
	- \mathbb{E} \int_{\partial Q} [s  \mathcal{O}_{\mu}((sh)^{2}) + s  \mathcal{O}_{\mu}((sh)^{4})] t_{r} (|D_{h}^3 v|^{2})  d t 
\\ \notag
& \quad
		- \mathbb{E} \int_{\partial Q} [s^{4} \mathcal{O}_{\mu}(1) + s^{5} \mathcal{O}_{\mu}((sh)^{2}) ] t_{r} ( |D_{h} v|^{2}) d t\\
	& \geq 
	60 \mathbb{E} \int_{Q^{*}} s^{5} \mu^{6} \varphi^{5} |\partial_{x} \psi|^{6} |D_{h} v |^{2} d t 
	+ 12 \mathbb{E} \int_{\partial Q} s^{5} \mu^{5} \varphi^{5} |\partial_{x} \psi|^{5} t_{r} (|D_{h} v|^{2})  d t 
	- \mathcal{A}_{1} - \mathcal{B}
	.
\end{align}

From \cref{eqDu2,eqDIgbp,eqAu2,eqAIgbp}, letting $ q_{13} = r^{2} A_{h}^{2} D_{h}^{2} \rho A_{h}^{3} D_{h} \rho $, for, we obtain 
\begin{align}
	\notag
	\label{eqI13e1}
	I_{13} 
	& = 
	24 \mathbb{E} \int_{Q} q_{13} A_{h}^{2} D_{h}^{2} v A_{h} D_{h}^{3} v d t 
	= 
	12 \mathbb{E} \int_{Q} q_{13} D_{h} [ |A_{h} D_{h}^{2} v|^{2} ]d t 
	\\\notag
	& = 
	-12 \mathbb{E} \int_{Q^{*}} D_{h} q_{13}  |A_{h} D_{h}^{2} v|^{2} d t 
	+ 12  \mathbb{E} \int_{\partial Q}  q_{13} t_{r} (|A_{h} D_{h}^{2} v|^{2}) \nu d t 
	\\\notag
	& = 
	-12 \mathbb{E} \int_{Q^{*}} D_{h} q_{13}  A_{h} (|  D_{h}^{2} v|^{2}) d t 
	+ 3 h^{2} \mathbb{E} \int_{Q^{*}} D_{h} q_{13}     |  D_{h}^{3} v|^{2} d t 
	+ 12  \mathbb{E} \int_{\partial Q}  q_{13} t_{r} (|A_{h} D_{h}^{2} v|^{2}) \nu d t 
	\\\notag
	& = 
	-12 \mathbb{E} \int_{\overline{Q}} A_{h} D_{h} q_{13}   |  D_{h}^{2} v|^{2} d t 
	+ 6 h \mathbb{E} \int_{\partial Q^{*}}   D_{h} q_{13}  t_{r} (|  D_{h}^{2} v|^{2}) d t 
	+ 3 h^{2} \mathbb{E} \int_{Q^{*}} D_{h} q_{13}     |  D_{h}^{3} v|^{2} d t 
	\\
	& \quad 
	+ 12  \mathbb{E} \int_{\partial Q}  q_{13} t_{r} (|A_{h} D_{h}^{2} v|^{2}) \nu d t 
	.
\end{align}
From \cref{eqdpr2pr,eqpr2pr}, for $ \lambda h (\delta T^{2})^{-1}\leq 1$, we have 
\begin{align}
	\notag
	\label{eqI13e2}
	-12 \mathbb{E} \int_{\overline{Q}} A_{h} D_{h} q_{13}   |  D_{h}^{2} v|^{2} d t 
	& \geq 
	36 \mathbb{E} \int_{\overline{Q}} s^{3} \mu^{4} \varphi^{3} |\partial_{x} \psi|^{4} |D_{h}^{2} v |^{2} d t 
	\\
	& \quad 
	- \mathbb{E} \int_{\overline{Q}} [s^{2} \mathcal{O}_{\mu}(1) + s^{3} \mathcal{O}_{\mu}((sh)^{2}) + s^{3} \mu^{3} \varphi^{3} \mathcal{O}(1) ] |D^{2}_{h} v|^{2} d t
	,
\end{align}
and
\begin{align}
	\notag
	\label{eqI13e3}
	& 6 h \mathbb{E} \int_{\partial Q^{*}}   D_{h} q_{13}  t_{r} (|  D_{h}^{2} v|^{2}) d t 
	+ 3 h^{2} \mathbb{E} \int_{Q^{*}} D_{h} q_{13}     |  D_{h}^{3} v|^{2} d t 
	\\
	& 
	\geq 
	-\mathbb{E} \int_{\partial Q}   s^{2} \mathcal{O}_{\mu}(1)   |D^{2}_{h} v|^{2} d t
	-\mathbb{E} \int_{Q^{*}}   s  \mathcal{O}_{\mu}((sh)^{2}) |D_{h}^{3} v |^{2} d t
	.
\end{align}
From \cref{eqdpr2pr,eqpr2pr,tra} and Cauchy-Schwarz inequality, for $ \lambda h (\delta T^{2})^{-1}\leq 1$, it holds that 
\begin{align}
	\notag
	\label{eqI13e4}
	& 12  \mathbb{E} \int_{\partial Q}  q_{13} t_{r} (|A_{h} D_{h}^{2} v|^{2}) \nu d t 
	\\ \notag
	& = 
	12  \mathbb{E} \int_{\partial Q}  q_{13} |D_{h}^{2} v|^{2} \nu d t
	+  3 h^{2} \mathbb{E} \int_{\partial Q}  q_{13}  t_{r} (|D_{h}^{3} v|^{2}) \nu d t
	- 12 h  \mathbb{E} \int_{\partial Q}  q_{13} D^{2}_{h} v   t_{r} (D_{h}^{3} v)  d t
	\\ \notag
	& \geq 
	- 12  \mathbb{E} \int_{\partial Q}  s^{3} \mu^{3} \varphi^{3} (\partial_{x} \psi)^{3} |D_{h}^{2} v|^{2} \nu d t
	- \mathbb{E} \int_{\partial Q} [s^{2} \mathcal{O}_{\mu}(1) + s^{3} \mathcal{O}_{\mu}(sh) ]  |D^{2}_{h} v|^{2} d t
	\\
	& \quad 
	- \mathbb{E} \int_{\partial Q} s \mathcal{O}_{\mu}(s h) t_{r}(|D_{h}^{3} v|^{2}) d t.
\end{align}
Combining \cref{eqI13e1,eqI13e2,eqI13e3,eqI13e4}, for $ \lambda h (\delta T^{2})^{-1}\leq 1$ and $ \lambda \geq T^{2} $, we have 
\begin{align}
	\label{eqI13}
	I_{13} 
	\geq 
	36 \mathbb{E} \int_{\overline{Q}} s^{3} \mu^{4} \varphi^{3} |\partial_{x} \psi|^{4} |D_{h}^{2} v |^{2} d t 
	+ 12   \mathbb{E} \int_{\partial Q}  s^{3} \mu^{3} \varphi^{3} |\partial_{x} \psi|^{3} |D_{h}^{2} v|^{2}   d t
	- \mathcal{A}_{1} - \mathcal{B}.
\end{align}
 
From \cref{eqA2u,eqvBoundary,traBoundary,eqDIgbp,eqDu2,eqDuv,eqdprpra,eqpr2pra}, noting that Cauchy-Schwarz inequality, letting $ q_{14} = r A_{h}^{2} D_{h}^{2} \rho A_{h}D_{h} (r A_{h} D_{h}^{3} \rho)$, for $ \lambda h (\delta T^{2})^{-1}\leq 1$ and $ \lambda \geq T^{2} $, we have 
\begin{align}
	\notag
	\label{eqI14}
	I_{14}
	& = 
	12 \mathbb{E} \int_{Q} q_{14} A_{h}^{2} D_{h}^{2} v v d t 
	= 
	12 \mathbb{E} \int_{Q} q_{14}  D_{h}^{2} v v d t 
	+ 3 h^{2} \mathbb{E} \int_{Q} q_{14}  D_{h}^{4} v v d t 
	\\ \notag
	& = 
	-12 \mathbb{E} \int_{Q^{*}} D_{h}q_{14} A_{h} v D_{h} v  d t
	-12 \mathbb{E} \int_{Q^{*}} A_{h}q_{14}  |D_{h} v|^{2}  d t
	- 3 h^{2} \mathbb{E} \int_{Q^{*}} D_{h}  q_{14} A_{h} v D_{h}^{3} v  d t 
	\\ \notag
	& \quad 
	- 3 h^{2} \mathbb{E} \int_{Q^{*}} A_{h}  q_{14} D_{h} v D_{h}^{3} v  d t 
	\\ \notag
	& = 
	6 \mathbb{E} \int_{\overline{Q}} D_{h}^{2} q_{14} | v |^{2}  d t
	-12 \mathbb{E} \int_{Q^{*}} A_{h}q_{14}  |D_{h} v|^{2}  d t
	- 3 h^{2} \mathbb{E} \int_{Q^{*}} D_{h}  q_{14} A_{h} v D_{h}^{3} v  d t 
	\\ \notag
	& \quad 
	- 3 h^{2} \mathbb{E} \int_{Q^{*}} A_{h}  q_{14} D_{h} v D_{h}^{3} v  d t 
	\\ \notag
	& \geq 
	36 \mathbb{E} \int_{Q^{*}} s^{5} \mu^{6} \varphi^{5} |\partial_{x} \psi|^{6} |D_{h} v|^{2} d t 
	- \mathbb{E} \int_{Q} s^{5} \mathcal{O}_{\mu}(1)  |v|^{2} d t 
	- \mathbb{E} \int_{Q^{*}} s \mathcal{O}_{\mu}((sh)^{2}) |D_{h}^{3} v|^{2} d t 
	\\ \notag 
	& \quad 
	- \mathbb{E} \int_{Q^{*}} [
		s^{5} \mu^{5} \varphi^{5} \mathcal{O}(1)
		+ s^{5} \mathcal{O}_{\mu}((s h)^{2}) 
		+ s^{4}\mathcal{O}_{\mu}(1)
	] |D_{h} v|^{2} d t
	\\
	& \geq 
	36 \mathbb{E} \int_{Q^{*}} s^{5} \mu^{6} \varphi^{5} |\partial_{x} \psi|^{6} |D_{h} v|^{2} d t 
	- \mathcal{A}_{1}
	.
\end{align}

From \cref{eqAuv,eqDu2,eqA2u,eqADInt,eqDIgbp}, letting $ q_{22} = r^{2} D_{h}^{4} \rho A_{h} D_{h}^{3} \rho $, we have 
\begin{align}
	\notag
	\label{eqI22e1}
	I_{22}
	& =
	4 \mathbb{E} \int_{Q} q_{22} A_{h}^{4} v A_{h}^{3} D_{h} v d t
	=  
	4 \mathbb{E} \int_{Q} q_{22} A_{h} ( A_{h}^{3} v A_{h}^{2} D_{h} v ) d t
	- h^{2} \mathbb{E} \int_{Q} q_{22}   A_{h}^{3} D_{h} v A_{h}^{2} D_{h}^{2} v  d t
	\\ \notag
	& =
	2 \mathbb{E} \int_{Q} q_{22} A_{h} D_{h}   (|A_{h}^{2}  v |^{2})  d t
	- \frac{h^{2}}{2} \mathbb{E} \int_{Q} q_{22}  D_{h}   (|A_{h}^{2} D_{h}  v |^{2})  d t
	\\ \notag
	& =
	2 \mathbb{E} \int_{Q} q_{22} A_{h} D_{h}   \bigg(\bigg| v + \frac{h^{2}}{4} D_{h}^{2} v  \bigg|^{2}\bigg)  d t
	- \frac{h^{2}}{2} \mathbb{E} \int_{Q} q_{22}  D_{h}   \bigg(\bigg| D_{h}  v  + \frac{h^{2}}{4} D_{h}^{3} v \bigg|^{2}\bigg)  d t
	\\ \notag
	& =
	- 2 \mathbb{E} \int_{\overline{Q}} A_{h} D_{h} q_{22}     \bigg| v + \frac{h^{2}}{4} D_{h}^{2} v  \bigg|^{2}   d t
	+ h \mathbb{E} \int_{\partial Q^{*}}   D_{h} q_{22} t_{r}   \bigg( \bigg| v + \frac{h^{2}}{4} D_{h}^{2} v  \bigg|^{2}  \bigg) d t
	\\ \notag
	& \quad 
	+  2 \mathbb{E} \int_{\partial Q }   q_{22} t_{r}   \bigg( A_{h} \bigg| v + \frac{h^{2}}{4} D_{h}^{2} v  \bigg|^{2}  \bigg) \nu  d t
	+ \frac{h^{2}}{2} \mathbb{E} \int_{Q^{*}} D_{h} q_{22}     \bigg| D_{h}  v  + \frac{h^{2}}{4} D_{h}^{3} v \bigg|^{2}   d t
	\\
	& \quad 
	- \frac{h^{2}}{2} \mathbb{E} \int_{\partial Q} q_{22} t_{r}  \bigg(\bigg| D_{h}  v  + \frac{h^{2}}{4} D_{h}^{3} v \bigg|^{2}\bigg) \nu d t
	.
\end{align}
From \cref{eqdpr2pr,eqpr2pr}, for $ \lambda h (\delta T^{2})^{-1}\leq 1$, we obtain 
\begin{align}
	\notag
	\label{eqI22e2}
	& - 2 \mathbb{E} \int_{\overline{Q}} A_{h} D_{h} q_{22}     \bigg| v + \frac{h^{2}}{4} D_{h}^{2} v  \bigg|^{2}   d t
	\\ \notag
	& \geq 
	14 \mathbb{E} \int_{Q} s^{7} \mu^{8} \varphi^{7} |\partial_{x} \psi|^{8} |v|^{2} d t 
	- \mathbb{E} \int_{Q} [s^{6} \mathcal{O}_{\mu}(1) + s^{7} \mathcal{O}_{\mu}((sh)^{2}) +  s^{7} \mu^{7} \varphi^{7} \mathcal{O}(1)] |v|^{2} d t
	\\
	& \quad 
	- \mathbb{E} \int_{\overline{Q}} s^{3}  \mathcal{O}_{\mu}((sh)^{2}) |D_{h}^{2} v |^{2} d t 
	.
\end{align}
Thanks to \cref{eqdpr2pr,traBoundary}, for $ \lambda h (\delta T^{2})^{-1}\leq 1$, we have 
\begin{align}
		\notag
	\label{eqI22e3}
	&h \mathbb{E} \int_{\partial Q^{*}}   D_{h} q_{22} t_{r}   \bigg( \bigg| v + \frac{h^{2}}{4} D_{h}^{2} v  \bigg|^{2}  \bigg) d t
	- \frac{h^{2}}{2} \mathbb{E} \int_{\partial Q} q_{22} t_{r}  \bigg(\bigg| D_{h}  v  + \frac{h^{2}}{4} D_{h}^{3} v \bigg|^{2}\bigg) \nu d t
		\\ \notag
	&  \geq 
		- \mathbb{E} \int_{\partial Q} s^{2} \mathcal{O}_{\mu}((sh)^{5})  |D_{h}^{2} v| d t -\mathbb{E} \int_{\partial Q} s^{5} \mathcal{O}_{\mu}((sh)^{2}) t_{r} (|D_{h} v|^{2}) d t  
	\\
	& \quad 
	- \mathbb{E} \int_{\partial Q} s  \mathcal{O}_{\mu}((sh)^{6}) t_{r} (|D_{h}^{3} v|^{2}) d t 
	.
\end{align}
and 
\begin{align}
	\notag
	\label{eqI22e4}
	&	2\mathbb{E} \int_{\partial Q }   q_{22} t_{r}   \bigg( A_{h} \bigg| v + \frac{h^{2}}{4} D_{h}^{2} v  \bigg|^{2}  \bigg) \nu  d t
	\\  
	& \geq 
	- \mathbb{E} \int_{\partial Q} s^{3} \mathcal{O}_{\mu}((sh)^{2})  |D_{h}^{2} v| d t -\mathbb{E} \int_{\partial Q} s^{5} \mathcal{O}_{\mu}((sh)^{2}) t_{r} (|D_{h} v|^{2}) d t  
	- \mathbb{E} \int_{\partial Q} s  \mathcal{O}_{\mu}((sh)^{2}) t_{r} (|D_{h}^{3} v|^{2}) d t 
	.
\end{align}
From \cref{eqdpr2pr}, for $ \lambda h (\delta T^{2})^{-1}\leq 1$, we obtain 
\begin{align}
	\label{eqI22e5}
	\frac{h^{2}}{2} \mathbb{E} \int_{Q^{*}} D_{h} q_{22}     \bigg| D_{h}  v  + \frac{h^{2}}{4} D_{h}^{3} v \bigg|^{2}   d t
	\geq 
	-\mathbb{E} \int_{Q^{*}}  s^{5} \mathcal{O}_{\mu}((sh)^{2}) |D_{h} v |^{2} d t 
	- \mathbb{E} \int_{Q^{*}}  s  \mathcal{O}_{\mu}((sh)^{6}) |D_{h}^{3} v |^{2} d t 
	,
\end{align}
Hence, combining \cref{eqI22e1,eqI22e2,eqI22e3,eqI22e4,eqI22e5}, for all $ \mu \geq 1 $, one can find  constant $ h'_{1} \leq h_1 $, such that for all $ \lambda \geq T + T^{2} $, $ h\leq h_{2} $ and $ \lambda h (\delta T^{2})^{-1}\leq 1 $, it holds that
\begin{align}
	\label{eqI22}
	I_{22}
	\geq 
	14 \mathbb{E} \int_{Q} s^{7} \mu^{8} \varphi^{7} |\partial_{x} \psi|^{8} |v|^{2} d t 
	- \mathcal{A}_{1} - \mathcal{B}
	.
\end{align}

From \cref{eqA2u}, letting $ q_{23} = r^{2} D_{h}^{4} \rho A_{h}^{3} D_{h} \rho $, for  $ \lambda h (\delta T^{2})^{-1}\leq 1$, we have 
\begin{align}
	\label{eqI23e1}
	I_{23}	
	& =
	4 \mathbb{E} \int_{Q}  q_{23} A_{h}^{4} v  A_{h} D_{h}^{3} v d t 
	=
	4 \mathbb{E} \int_{Q}  q_{23} \bigg(  v + \frac{h^{2}}{4} D_{h}^{2} v + \frac{h^{2}}{4} A_{h}^{2} D_{h}^{2} v \bigg) A_{h} D_{h}^{3} v d t 
	.
\end{align}
From \cref{eqDIgbp,eqvBoundary,traBoundary,eqDuv,eqAIgbp,eqDu2,eqdpr2pr,eqpr2pr,eqAu2,eqD2Ibp}, for  $ \lambda h (\delta T^{2})^{-1}\leq 1$, we obtain 
\begin{align}
	\notag
	\label{eqI23e2}
	& 4 \mathbb{E} \int_{Q} q_{23} v A_{h} D_{h}^{3} v d t =4 \mathbb{E} \int_{\overline{Q}} D_h^2(q_{23} v) A_{h} D_{h} v d t
	\\ \notag
	& =
	4 \mathbb{E} \int_{\overline{Q}} D^{2}_{h} q_{23} A_{h}^{2} v A_{h} D_{h} v d t 
	+ 8 \mathbb{E} \int_{\overline{Q}} A_{h} D_{h} q_{23} |A_{h} D_{h} v|^{2} d t 
	+4 \mathbb{E} \int_{\overline{Q}} A_{h}^{2} q_{23}  A_{h} D_{h} v D_{h}^{2} v d t 
	\\ \notag
	& = 
	-2 \mathbb{E} \int_{\overline{Q^{*}}} D^{3}_{h} q_{23} |A_{h} v|^{2}  d t 
	+ 6 \mathbb{E} \int_{\overline{Q^{*}}} A_{h}^{2} D_{h} q_{23} |D_{h} v|^{2}  d t 
	- 2 h^{2} \mathbb{E} \int_{\overline{Q}} A_{h} D_{h} q_{23} |D_{h}^{2} v|^{2}  d t
	\\ \notag
	& \geq 
	- 30 \mathbb{E} \int_{Q^{*}} s^{5} \mu^{6} \varphi^{5} |\partial_{x} \psi|^{6} |D_{h}v|^{2} d t 
	- \mathbb{E} \int_{\overline{Q}} s^{3} \mathcal{O}_{\mu}((sh)^{2}) |D_{h}^{2} v|^{2} d t 
	- \mathbb{E} \int_{Q} s^{5} \mathcal{O}_{\mu}(1) | v|^{2} d t 
	\\
	& \quad 
	- \mathbb{E} \int_{Q^{*}} [ 
	s^{5} \mu^{5} \varphi^{5}  \mathcal{O}(1) 
	+ s^{4} \mathcal{O}_{\mu}(1)
	+ s^{5} \mathcal{O}_{\mu}((sh)^{2})
	]|D_{h} v|^{2} d t 
	.
\end{align}
From \cref{eqAIgbp,eqAuv,eqDu2,eqDIgbp,eqdpr2pr}, for  $ \lambda h (\delta T^{2})^{-1}\leq 1$, we get 
\begin{align}
	\notag
	\label{eqI23e3}
	&
	h^{2} \mathbb{E} \int_{Q} q_{23} D_{h}^{2} v A_{h} D_{h}^{3} v d t 
	\\ \notag
	& = 
	h^{2} \mathbb{E} \int_{Q^{*}} A_{h} q_{23} A_{h} D_{h}^{2} v D_{h}^{3} v d t 
	+ \frac{h^{2}}{4} \mathbb{E} \int_{Q^{*}} D_{h} q_{23} |D_{h}^{3} v|^{2} d t 
	- \frac{h^{3}}{2} \mathbb{E} \int_{\partial Q} q_{23} D_{h}^{2} v t_{r}(D_{h}^{3} v)d t 
	\\ \notag
	& = 
	- \frac{h^{2}}{2} \mathbb{E} \int_{\overline{Q}} A_{h} D_{h} q_{23}  |D_{h}^{2} v|^{2} d t 
	+ \frac{h^{2}}{2} \mathbb{E} \int_{\partial Q^{*}} A_{h} q_{23} t_{r} (|D_{h}^{2} v |^{2}) \nu d t 
	+ \frac{h^{4}}{4} \mathbb{E} \int_{Q^{*}} D_{h} q_{23} |D_{h}^{3} v|^{2} d t
	\\ \notag
	& \quad 
	- \frac{h^{3}}{2} \mathbb{E} \int_{\partial Q} q_{23} D_{h}^{2} v t_{r}(D_{h}^{3} v)d t 
	\\ \notag
	& \geq 
	-\mathbb{E} \int_{\overline{Q}} s^{3} \mathcal{O}_{\mu}((sh)^{2}) |D_{h}^{2} v |^{2} d t 
	- \mathbb{E} \int_{Q^{*}} s  \mathcal{O}_{\mu}((sh)^{4}) |D_{h}^{3} v |^{2} d t 
	- \mathbb{E} \int_{\partial Q} s^{3} \mathcal{O}_{\mu}((sh)^{2})  t_{r} (|D_{h}^{2} v |^{2}) d t 
	\\
	& \quad 
	- \mathbb{E} \int_{\partial Q} s \mathcal{O}_{\mu}((sh)^{4}) t_{r} (|D_{h}^{3} v |^{2}) d t 
	.
\end{align}
From \cref{eqDu2,eqDIgbp,eqdpr2pr,eqAu2leq,eqAIgbp,tra}, for $ \lambda h (\delta T^{2})^{-1}\leq 1$, it holds that 
\begin{align}
	\notag
	\label{eqI23e4}
	& 
	h^{2} \mathbb{E} \int_{Q} q_{23} A_{h}^{2} D_{h}^{2} v A_{h} D_{h}^{3} v d t 
	\\ \notag
	& = 
	- \frac{h^{2}}{2} \mathbb{E} \int_{Q^{*}} D_{h} q_{23} |A_{h}  D_{h}^{2} v|^{2}  d t 
	+\frac{h^{2}}{2} \mathbb{E} \int_{\partial Q} q_{23}  t_{r}(|A_{h} D_{h}^{2} v|^{2}) \nu d t 
	\\ \notag
	& =
	-\frac{h^{2}}{2} \mathbb{E} \int_{Q^{*}} D_{h} q_{23} |A_{h}  D_{h}^{2} v|^{2}  d t 
	+ \frac{h^{2}}{2} \mathbb{E} \int_{\partial Q} q_{23}  \bigg(
		|D_{h}^{2} v |^{2}\nu
		- h D_{h}^{2} v t_{r} (D_{h}^{3} v) 
		+ \frac{h^{2}}{4} t_{r}(|D_{h}^{3} v|^{2}) \nu
	\bigg) d t 
	\\
	& \geq 
	-\mathbb{E} \int_{\overline{Q}} s^{3} \mathcal{O}_{\mu}((sh)^{2}) |D_{h}^{2} v|^{2} d t 
	- \mathbb{E} \int_{\partial Q} s^{3} \mathcal{O}_{\mu}((sh)^{2}) |D_{h}^{2} v|^{2} d t 
	-\mathbb{E} \int_{\partial Q} s \mathcal{O}_{\mu}((sh)^{4}) t_{r}(|D_{h}^{3} v|^{2}) d t 
	.
\end{align}
Combining \cref{eqI23e1,eqI23e2,eqI23e3,eqI23e4}, for $ \lambda h (\delta T^{2})^{-1}\leq 1$ and $ \lambda \geq T^{2} $, we have 
\begin{align}
	\label{eqI23}
	I_{23} \geq 
	- 30 \mathbb{E} \int_{Q^{*}} s^{5} \mu^{6} \varphi^{5} |\partial_{x} \psi|^{6} |D_{h}v|^{2} d t 
	- \mathcal{A}_{1} - \mathcal{B}
	.
\end{align}

From \cref{eqA2u}, letting $ q_{24} = r D_{h}^{4} \rho  A_{h} D_{h} (r A_{h} D_{h}^{3} \rho) $,  we have 
\begin{align}
	\label{eqI24e1}
	I_{24} 
	& = 
	2 \mathbb{E} \int_{Q} q_{24} A_{h}^{4} v v d t 
	=
	2 \mathbb{E} \int_{Q} q_{24} \bigg( 
		v 
		+ \frac{h^{2}}{4}  D_{h}^{2} v   
		+ \frac{h^{2}}{4} A_{h}^{2} D_{h}^{2} v  
	\bigg) v d t 
	.
\end{align}
From \cref{eqdprpra,eqpr2pra}, for $ \lambda h (\delta T^{2})^{-1}\leq 1$, we have 
\begin{align}
	\label{eqI24e2}
	2 \mathbb{E} \int_{Q} q_{24} v^{2} d  t 
	\geq 
	- 6 \mathbb{E} \int_{Q} s^{7} \mu^{8} \varphi^{7} |\partial_{x} \psi|^{8} v^{2} d t 
	- \mathbb{E} \int_{Q} [s^{6} \mathcal{O}_{\mu}(1) + s^{7} \mathcal{O}_{\mu}((sh)^{2}) +  s^{7} \mu^{7} \varphi^{7} \mathcal{O}(1)] |v|^{2} d t
	.
\end{align}
From \cref{eqvBoundary,traBoundary,eqDIgbp,eqDuv,eqDu2,eqdprpr}, for $ \lambda h (\delta T^{2})^{-1}\leq 1$, we obtain 
\begin{align}
	\label{eqI24e3}
	\notag
	\frac{h^{2}}{2} \mathbb{E} \int_{Q} q_{24} D_{h}^{2} v v d t 
	& = 
	- \frac{h^{2}}{2} \mathbb{E} \int_{Q^{*}} D_{h} q_{24}  A_{h }v D_{h} v d t 
	- \frac{h^{2}}{2} \mathbb{E} \int_{Q^{*}} A_{h} q_{24}  |D_{h }v|^{2}  d t 
	\\ \notag
	& = 
	\frac{h^{2}}{4} \mathbb{E}\int_{\overline{Q}} D_{h}^{2} q_{24} v^{2} d t 
	- \frac{h^{2}}{2} \mathbb{E} \int_{Q^{*}} A_{h} q_{24}  |D_{h }v|^{2}  d t 
	\\
	& \geq 
	-\mathbb{E}\int_{Q}  s^{5} \mathcal{O}_{\mu}(1)  v^{2} d t 
	- \mathbb{E}\int_{Q^{*}}  s^{5} \mathcal{O}_{\mu}((sh)^{2})  |D_{h} v|^{2} d t 
	.
\end{align}
From \cref{eqDIgbp,eqvBoundary,eqDuv,eqA2u,eqDu2,eqdprpr,eqAIgbp,eqAu2leq},  thanks to Cauchy-Schwarz inequality, for $ \lambda h (\delta T^{2})^{-1}\leq 1$,  we obtain 
\begin{align}
	\label{eqI24e4}
	\notag
	& 
	\frac{h^{2}}{2} \mathbb{E} \int_{Q} q_{24} A_{h}^{2} D_{h}^{2} v v d t 
	\\ \notag
	& = 
	- \frac{h^{2}}{2} \mathbb{E} \int_{Q^{*}} D_{h} q_{24} A_{h} v A_{h}^{2} D_{h} v  d t 
	- \frac{h^{2}}{2} \mathbb{E} \int_{Q^{*}} A_{h} q_{24} D_{h} v A_{h}^{2} D_{h} v  d t 
	\\ \notag
	& = 
	- \frac{h^{2}}{2} \mathbb{E} \int_{Q^{*}} D_{h} q_{24} A_{h} v  D_{h} v  d t 
	- \frac{h^{4}}{8} \mathbb{E} \int_{Q^{*}} D_{h} q_{24} A_{h} v  D_{h}^{3} v  d t 
	- \frac{h^{2}}{2} \mathbb{E} \int_{Q^{*}} A_{h} q_{24} |D_{h} v|^{2}  d t 
	\\ \notag
	& \quad 
	- \frac{h^{4}}{8} \mathbb{E} \int_{Q^{*}} A_{h} q_{24}  D_{h} v D_{h}^{3} v   d t \\
	& \geq
	-\mathbb{E} \int_{Q} s^{5} \mathcal{O}_{\mu}((sh)^{2}) |v|^{2} d t
	- \mathbb{E} \int_{Q^{*}} s^{5} \mathcal{O}_{\mu}((sh)^{2}) |D_{h} v|^{2} d t 
	- \mathbb{E} \int_{Q^{*}} s  \mathcal{O}_{\mu}((sh)^{4}) |D_{h}^{3} v|^{2} d t 
	.
\end{align}
Combining  \cref{eqI24e1,eqI24e2,eqI24e3,eqI24e4}, for $ \lambda h (\delta T^{2})^{-1}\leq 1$ and  $ \lambda \geq T^{2} $, we have 
\begin{align}
	\label{eqI24}
	I_{24}
	\geq 
	- 6 \mathbb{E} \int_{Q} s^{7} \mu^{8} \varphi^{7} |\partial_{x} \psi|^{8} v^{2} d t 
	- \mathcal{A}_{1}
	.
\end{align}

From \cref{eqDIgbp,eqDuv,eqA2u}, letting $ q_{32} = r^{2} A_{h}^{4} \rho A_{h} D_{h}^{3} \rho $, we have 
\begin{align}
	\notag
	\label{eqI32e1}
	I_{32} 
	& = 
	4 \mathbb{E} \int_{Q} q_{32} D_{h}^{4} v A_{h}^{3} D_{h} v d t 
	=
	4 \mathbb{E} \int_{Q} q_{32} D_{h}^{4} v A_{h} D_{h} v d t 
	+ h^{2} \mathbb{E} \int_{Q} q_{32} D_{h}^{4} v A_{h} D_{h}^{3} v d t 
	\\ \notag
	& = 
	- 4 \mathbb{E} \int_{Q^{*}} D_{h} q_{32} A_{h}^{2} D_{h} v D_{h}^{3} v  d t 
	- 4 \mathbb{E} \int_{Q^{*}} A_{h} q_{32} A_{h} D_{h}^{2} v D_{h}^{3} v  d t 
	+ 4 \mathbb{E} \int_{\partial Q} q_{32} A_{h} D_{h} v t_{r} (D_{h}^{3} v) \nu  d t 
	\\
	& \quad 
	+ h^{2} \mathbb{E} \int_{Q} q_{32} D_{h}^{4} v A_{h} D_{h}^{3} v d t 
	.
\end{align}
From \cref{eqA2u,eqDIgbp,eqDuv,eqvBoundary,traBoundary,eqDu2,eqdpr2pr,eqpr2pr}, for $ \lambda h (\delta T^{2})^{-1}\leq 1$, we obtain 
\begin{align}
	\notag
	\label{eqI32e2}
	&
	- 4 \mathbb{E} \int_{Q^{*}} D_{h} q_{32} A_{h}^{2} D_{h} v D_{h}^{3} v  d t 
	\\ \notag
	& = 
	- 4 \mathbb{E} \int_{Q^{*}} D_{h} q_{32}  D_{h} v D_{h}^{3} v  d t 
	- h^{2} \mathbb{E} \int_{Q^{*}} D_{h} q_{32}  |D_{h}^{3} v|^{2}  d t 
	\\ \notag
	& =
	4 \mathbb{E} \int_{\overline{Q}} D_{h}^{2} q_{32} A_{h} D_{h} v D_{h}^{2} v  d t 
	+ 4 \mathbb{E} \int_{\overline{Q}} A_{h} D_{h} q_{32}  |D_{h}^{2} v|^{2}  d t 
	- h^{2} \mathbb{E} \int_{Q^{*}} D_{h} q_{32}  |D_{h}^{3} v|^{2}  d t 
	\\ \notag
	& =
	- 2 \mathbb{E} \int_{\overline{Q^{*}}} D_{h}^{3} q_{32} |D_{h} v|^{2}   d t 
	+ 4 \mathbb{E} \int_{\overline{Q}} A_{h} D_{h} q_{32}  |D_{h}^{2} v|^{2}  d t 
	- h^{2} \mathbb{E} \int_{Q^{*}} D_{h} q_{32}  |D_{h}^{3} v|^{2}  d t 
	\\ \notag
	& \geq 
	- 12 \mathbb{E} \int_{\overline{Q}} s^{3} \mu^{4} \varphi^{3} |\partial_{x} \psi|^{4} |D_{h}^{2} v |^{2} d t 
	- \mathbb{E} \int_{Q^{*}} s^{3} \mathcal{O}_{\mu}(1) |D_{h} v|^{2} d t 	
	- \mathbb{E} \int_{Q^{*}} s \mathcal{O}_{\mu}((sh)^{2}) |D_{h}^{3} v|^{2} d t 	
	\\
	& \quad 
	- \mathbb{E} \int_{\overline{Q}} [s^{2} \mathcal{O}_{\mu}(1) + s^{3} \mathcal{O}_{\mu}((sh)^{2}) + s^{3} \mu^{3} \varphi^{3} \mathcal{O}(1) ] |D^{2}_{h} v|^{2} d t
	.
\end{align}
From \cref{eqDu2,eqDIgbp,eqdpr2pr,eqpr2pr}, for $ \lambda h (\delta T^{2})^{-1}\leq 1$, it holds that 
\begin{align}
	\notag
	\label{eqI32e3}
	& 
	- 4 \mathbb{E} \int_{Q^{*}} A_{h} q_{32} A_{h} D_{h}^{2} v D_{h}^{3} v  d t 
	\\ \notag
	& =
	2 \mathbb{E} \int_{\overline{Q}} A_{h}  D_{h} q_{32} |D_{h}^{2} v |^{2}   d t 
	- 2 \mathbb{E} \int_{ \partial Q^{*}} A_{h}  q_{32} t_{r} (|D_{h}^{2} v |^{2}) \nu   d t 
	\\ \notag
	& \geq 
	- 6 \mathbb{E} \int_{\overline{Q}} s^{3} \mu^{4} \varphi^{3} |\partial_{x} \psi|^{4} |D_{h}^{2} v |^{2} d t 
	+ 2 \mathbb{E} \int_{\partial Q^{*}} s^{3} \mu^{3} \varphi^{3} (\partial_{x} \psi)^{3} t_{r} (|D_{h}^{2} v |^{2}) \nu   d t 
	\\ \notag
	& \quad 
	- \mathbb{E} \int_{\overline{Q}} [s^{2} \mathcal{O}_{\mu}(1) + s^{3} \mathcal{O}_{\mu}((sh)^{2}) + s^{3} \mu^{3} \varphi^{3} \mathcal{O}(1) ] |D^{2}_{h} v|^{2} d t
	\\
	& \quad 
	-  \mathbb{E} \int_{\partial Q^{*}} [s^{2} \mathcal{O}_{\mu}(1) + s^{3} \mathcal{O}_{\mu}((sh)^{2})+s^{3} \mu^{2} \varphi^{3} \mathcal{O}(1)] t_{r} (|D_{h}^{2} v |^{2})   d t 
	.
\end{align}
From \cref{eqvBoundary,eqdpr2pr,eqpr2pr}, thanks to Cauchy-Schwarz inequality, for $ \lambda h (\delta T^{2})^{-1}\leq 1$, it holds that 
\begin{align}
	\notag
	\label{eqI32e4}
	&
	4 \mathbb{E} \int_{\partial Q} q_{32} A_{h} D_{h} v t_{r} (D_{h}^{3} v) \nu  d t 
	\geq 
	-\mathbb{E} \int_{\partial Q} s^{3} \mathcal{O}_{\mu}(1) |A_{h} D_{h} v| t_{r} (|D_{h}^{3} v|)  d t 
	\\ \notag
	& = 
	- \frac{h}{2}\mathbb{E} \int_{\partial Q} s^{3} \mathcal{O}_{\mu}(1) |D_{h}^{2} v| t_{r} (|D_{h}^{3} v|)  d t 
	\\
	& \geq 
	-\mathbb{E} \int_{\partial Q} s^{3} \mathcal{O}_{\mu}(s h) |D_{h}^{2} v|^{2} d t 
	-\mathbb{E} \int_{\partial Q} s \mathcal{O}_{\mu}(s h) t_{r} (|D_{h}^{3} v|^{2}) d t 
	.
\end{align}
From \cref{eqDu2,eqDIgbp,eqdpr2pr}, for $ \lambda h (\delta T^{2})^{-1}\leq 1$, we obtain 
\begin{align}
	\notag
	\label{eqI32e5}
	h^{2} \mathbb{E} \int_{Q} q_{32} D_{h}^{4} v A_{h} D_{h}^{3} v d t 
	& =
	- \frac{h^{2}}{2} \mathbb{E} \int_{Q^{*}} D_{h} q_{32}  |D_{h}^{3} v|^{2} d t 
	+ \frac{h^{2}}{2} \mathbb{E} \int_{\partial Q} q_{32}  t_{r}(|D_{h}^{3} v|^{2}) \nu d t 
	\\
	& \geq 
	-\mathbb{E} \int_{Q^{*}} s \mathcal{O}_{\mu}((sh)^{2})   |D_{h}^{3} v|^{2} d t 
	- \mathbb{E} \int_{\partial Q} s \mathcal{O}_{\mu}((sh)^{2})  t_{r}(|D_{h}^{3} v|^{2}) d t 
	.
\end{align}
Combining \cref{eqI32e1,eqI32e2,eqI32e3,eqI32e4,eqI32e5}, for $ \lambda h (\delta T^{2})^{-1}\leq 1$ and  $ \lambda \geq T^{2} $, we have 
\begin{align}
	\label{eqI32}
	I_{32} 
	\geq 
	- 18 \mathbb{E} \int_{\overline{Q}} s^{3} \mu^{4} \varphi^{3} |\partial_{x} \psi|^{4} |D_{h}^{2} v |^{2} d t 
	- 2 \mathbb{E} \int_{\partial Q} s^{3} \mu^{3} \varphi^{3} |\partial_{x} \psi|^{3} |D_{h}^{2} v |^{2}     d t 
	- \mathcal{A}_{1} - \mathcal{B}
	.
\end{align}

From \cref{eqDu2,eqDIgbp,eqdpr2pr,eqprpr}, letting $ q_{33} = r^{2} A_{h}^{4} \rho A_{h}^{3} D_{h} \rho $,  for $ \lambda h (\delta T^{2})^{-1}\leq 1$, we have 
\begin{align}
	\notag
	\label{eqI33}
	I_{33} 
	& = 
	4 \mathbb{E} \int_{Q} q_{33} D_{h}^{4} v A_{h} D_{h}^{3} v d t 
	=
	- 2 \mathbb{E} \int_{Q^{*}}D_{h} q_{33} | D_{h}^{3} v |^{2}  d t 
	+ 2 \mathbb{E} \int_{\partial Q} q_{33} t_{r} (|D_{h}^{3} v|) \nu d t 
	\\\notag
	& \geq 
	2 \mathbb{E} \int_{Q^{*}}s \mu^{2}  \varphi |\partial_{x} \psi|^{2} | D_{h}^{3} v |^{2}  d t 
	+ 2 \mathbb{E} \int_{\partial Q} s \mu  \varphi |\partial_{x} \psi| t_{r} (|D_{h}^{3} v|)  d t 
	\\\notag
	& \quad 
	- \mathbb{E} \int_{Q^{*}} [  \mathcal{O}_{\mu}(1) + s  \mathcal{O}_{\mu}((sh)^{2})+ s \mu \varphi \mathcal{O}(1)] |D_{h}^{3} v |^{2} d t\\
	& \quad
	- \mathbb{E} \int_{\partial Q} [\mathcal{O}_{\mu}(1) + s \mathcal{O}_{\mu}(sh)+ s  \varphi \mathcal{O}(1) ] t_{r}( |D^{3}_{h} v|^{2}) d t
	.
\end{align}

Letting $ q_{34} = r A_{h}^{4} \rho A_{h} D_{h}(r A_{h} D_{h}^{3} \rho) $, from \cref{eqD2Ibp,eqvBoundary,traBoundary,eqDuv,eqA2u,eqDu2,eqDIgbp,eqdprpra,eqpr2pra},   for $ \lambda h (\delta T^{2})^{-1}\leq 1$ and $ \lambda \geq T^{2} $, we have 
\begin{align}
	\notag
	\label{eqI34}
	I_{34}
	&  = 
	2 \mathbb{E} \int_{Q} q_{34} D_{h}^{4} v v d t 
	= 
	2 \mathbb{E} \int_{\overline{Q}} D_{h}^{2} (q_{34} v) D_{h}^{2} v d t 
	\\ \notag
	& =
	2 \mathbb{E} \int_{\overline{Q}} D_{h}^{2} q_{34} A_{h}^{2} v D_{h}^{2} v d t 
	+ 4 \mathbb{E} \int_{\overline{Q}} A_{h} D_{h} q_{34} A_{h} D_{h} v D_{h}^{2} v d t 
	+ 2 \mathbb{E} \int_{\overline{Q}} A_{h}^{2} q_{34} | D_{h}^{2} v |^{2} d t 
	\\ \notag
	& = 
	2 \mathbb{E} \int_{\overline{Q}} D_{h}^{2} q_{34} v D_{h}^{2} v d t 
	+ \frac{h^{2}}{2} \mathbb{E} \int_{\overline{Q}} D_{h}^{2} q_{34} | D_{h}^{2} v |^{2} d t - 2 \mathbb{E} \int_{\overline{Q^{*}}} A_{h} D_{h}^{2} q_{34} | D_{h} v |^{2} d t
	\\ \notag
	& \quad 
	+ 2 \mathbb{E} \int_{\overline{Q}} A_{h}^{2} q_{34} | D_{h}^{2} v |^{2} d t
	\\ \notag
	& \geq 
	- 6 \mathbb{E} \int_{\overline{Q}} s^{3} \mu^{4} \varphi^{3} |\partial_{x} \psi|^{4} |D_{h}^{2} v|^{2} d t 
	- \mathbb{E} \int_{Q^{*}} s^{3} \mathcal{O}_{\mu}(1) |D_{h} v|^{2} d t 
	- \mathbb{E} \int_{Q} s^4 \mathcal{O}_{\mu}(1) | v|^{2} d t 
	\\ \notag
	& \quad 
	- \mathbb{E} \int_{\overline{Q}} [s^{2} \mathcal{O}_{\mu}(1) + (s + s^{3}) \mathcal{O}_{\mu}((sh)^{2}) + s^{3} \mu^{3} \varphi^{3} \mathcal{O}(1) ] |D^{2}_{h} v|^{2} d t
	\\
	& \geq 
	- 6 \mathbb{E} \int_{\overline{Q}} s^{3} \mu^{4} \varphi^{3} |\partial_{x} \psi|^{4} |D_{h}^{2} v|^{2} d t 
	- \mathcal{A}_{1}
	.
\end{align}

Letting $ q_{42} = A_{h}D_{h} (r A_{h}^{2} D_{h}^{2} \rho) r A_{h} D_{h}^{3} \rho $ from \cref{eqA2u,eqAu2,eqAIgbp,eqAu2leq,eqAIgbp,eqdprpr,eqprpr}, thanks to Cauchy-Schwarz inequality, for $ \lambda h (\delta T^{2})^{-1}\leq 1$, we obtain 
\begin{align}
	\notag
	\label{eqI42}
	I_{42}
	& = 
	24 \mathbb{E}  \int_{Q} q_{42} A_{h} D_{h} v A_{h}^{3} D_{h} v d t 
	=
	24 \mathbb{E}  \int_{Q} q_{42} | A_{h} D_{h} v |^{2} d t 
	+ 6  h^{2} \mathbb{E}  \int_{Q} q_{42}  A_{h} D_{h} v A_{h} D_{h}^{3} v  d t 
	\\ \notag
	& =
	24 \mathbb{E}  \int_{Q} q_{42} A_{h} (|  D_{h} v |^{2}) d t 
	- 6  h^{2} \mathbb{E}  \int_{Q} q_{42}  | D_{h}^{2} v |^{2}  d t 
	+ 6  h^{2} \mathbb{E}  \int_{Q} q_{42}  A_{h} D_{h} v A_{h} D_{h}^{3} v  d t 
	\\ \notag
	& =
	24 \mathbb{E}  \int_{Q^{*}} A_{h} q_{42} |  D_{h} v |^{2} d t 
	- 12 h \mathbb{E}  \int_{\partial Q} q_{42} t_{r}(|  D_{h} v |^{2}) d t 
	- 6  h^{2} \mathbb{E}  \int_{Q} q_{42}  | D_{h}^{2} v |^{2}  d t 
	\\ \notag
	& \quad 
	+ 6  h^{2} \mathbb{E}  \int_{Q} q_{42}  A_{h} D_{h} v A_{h} D_{h}^{3} v  d t 
	\\ \notag
	& \geq 
	-48 \mathbb{E} \int_{Q^{*}} s^{5} \mu^{6} \varphi^{5} |\partial_{x} \psi|^{6} |D_{h} v|^{2} d t 
	- \mathbb{E} \int_{Q^{*}} [s^{4} \mathcal{O}_{\mu}(1) + s^{5} \mathcal{O}_{\mu}((sh)^{2}) + s^{5} \mu^{5} \varphi^{5} \mathcal{O}(1)] |D_{h} v|^{2} d t
	\\ \notag
	& \quad 
	- \mathbb{E} \int_{Q} s^{3} \mathcal{O}_{\mu}((sh)^{2}) |D_{h}^{2} v|^{2} d t 
	- \mathbb{E} \int_{Q^{*}} s \mathcal{O}_{\mu}((sh)^{2}) |D_{h}^{3} v|^{2} d t 
	- \mathbb{E} \int_{\partial Q} s^{4} \mathcal{O}_{\mu}((sh)) t_{r} (|D_{h}v|^{2}) d t 
	\\
	& \geq 
	-48 \mathbb{E} \int_{Q^{*}} s^{5} \mu^{6} \varphi^{5} |\partial_{x} \psi|^{6} |D_{h} v|^{2} d t 
	- \mathcal{A}_{1} - \mathcal{B}
	.
\end{align}

Letting  $ q_{43} = A_{h} D_{h} (r A_{h}^{2} D_{h}^{2} \rho) r A_{h}^{3} D_{h} \rho $, from \cref{eqDIgbp,eqDuv}, we have 
\begin{align}
	\notag
	\label{eqI43e1}
	I_{43} 
	& = 
	24 \mathbb{E} \int_{Q} q_{43} A_{h} D_{h} v A_{h} D_{h}^{3} v d t 
	\\\notag
	& =
	-24 \mathbb{E} \int_{Q^{*}} D_{h} q_{43} A_{h}^{2} D_{h} v A_{h} D_{h}^{2} v d t 
	-24 \mathbb{E} \int_{Q^{*}} A_{h} q_{43}  | A_{h} D_{h}^{2} v |^{2} d t 
	\\
	& \quad 
	+ 24 \mathbb{E} \int_{\partial Q} q_{43} A_{h} D_{h} v t_{r} (A_{h} D_{h}^{2} v) \nu d t 
	.
\end{align}
From \cref{eqDu2,eqDIgbp,eqAu2,eqAIgbp,eqvBoundary,traBoundary,eqdprpr,tra}, for $ \lambda h (\delta T^{2})^{-1}\leq 1$, we obtain 
\begin{align}
	\notag
	\label{eqI43e2}
	& 
	-24 \mathbb{E} \int_{Q^{*}} D_{h} q_{43} A_{h}^{2} D_{h} v A_{h} D_{h}^{2} v d t 
	\\\notag
	& =
	12 \mathbb{E} \int_{\overline{Q}} D_{h}^{2} q_{43} | A_{h} D_{h} v |^{2} d t 
	- 12 \mathbb{E} \int_{\partial Q^{*}} D_{h} q_{43} t_{r} (| A_{h} D_{h} v |^{2} ) \nu d t
	\\\notag
	& =  
	12 \mathbb{E} \int_{\overline{Q^{*}}} A_{h} D_{h}^{2} q_{43} | D_{h} v |^{2} d t 
	- 3 h^{2} \mathbb{E} \int_{\overline{Q}} D_{h}^{2} q_{43} |D_{h}^{2} v |^{2} d t 
	- 3 h^{2} \mathbb{E} \int_{\partial Q^*} D_{h} q_{43}  t_{r} (|D_{h}^{2} v |^{2} )\nu d t 
	\\ 
	& \geq 
	-\mathbb{E} \int_{Q^{*}} s^{3} \mathcal{O}_{\mu} (1) |D_{h} v|^{2} d t 
	- \mathbb{E} \int_{\overline{Q}} s  \mathcal{O}_{\mu} ((sh)^{2}) |D_{h}^{2} v|^{2} d t 
	- \mathbb{E} \int_{\partial Q} s \mathcal{O}_{\mu} ((sh)^{2})  |D_{h} ^{2} v|^{2} d t 
	.
\end{align}
From \cref{eqAu2,eqAIgbp,eqdprpr,eqprpr}, for $ \lambda h (\delta T^{2})^{-1}\leq 1$, it holds that 
\begin{align}
	\notag
	\label{eqI43e3}
	& -24 \mathbb{E} \int_{Q^{*}} A_{h} q_{43}  | A_{h} D_{h}^{2} v |^{2} d t 
	\\\notag
	& =
	-24 \mathbb{E} \int_{\overline{Q}} A_{h}^{2} q_{43} |D_{h}^{2} v |^{2} d t 
	+ 12 h \mathbb{E} \int_{\partial Q^{*}} A_{h} q_{43} t_{r}(|D_{h}^{2}v|^{2}) d t 
	+ 6 h^{2} \mathbb{E} \int_{Q^{*}} A_{h} q_{43} |D_{h}^{3} v|^{2} d t
	\\\notag
	& \geq 
	48 \mathbb{E} \int_{\overline{Q}} s^{3} \mu^{4} \varphi^{3} |\partial_{x} \psi|^{4} |D_{h}^{2} v|^{2} d t 
	- \mathbb{E} \int_{\overline{Q}} [s^{2} \mathcal{O}_{\mu}(1) + s^{3} \mathcal{O}_{\mu}((sh)^{2}) + s^{3} \mu^{3} \varphi^{3} \mathcal{O}(1) ] |D^{2}_{h} v|^{2} d t
	\\
	& \quad 
	- \mathbb{E} \int_{Q^{*}} s \mathcal{O}_{\mu}((sh)^{2}) |D_{h}^{3} v|^{2} d t 
	- \mathbb{E} \int_{\partial Q} s^{2} \mathcal{O}_{\mu}(sh) |D_{h}^{2} v|^{2} d t 
	.
\end{align}
From \cref{eqvBoundary} and Cauchy-Schwarz inequality, by $A_{h} D_{h} v=\frac{h}{2}D_h^2v$, $t_{r} (A_{h} D_{h}^{2} v)=D^2_hv-\frac{h}{2}t_{r} (D_{h}^{3})\nu$ on $\partial Q$, then for $ \lambda h (\delta T^{2})^{-1}\leq 1$, we have 
\begin{align}
	\notag
	\label{eqI43e4}
	24 \mathbb{E} \int_{\partial Q} q_{43} A_{h} D_{h} v t_{r} (A_{h} D_{h}^{2} v) \nu d t 
	& \geq 
	-\mathbb{E} \int_{\partial Q} s^{2} \mathcal{O}_{\mu}(sh)  |D_{h}^{2} v|^{2} d t 
	\\
	& \quad 
	- \mathbb{E} \int_{\partial Q}  \mathcal{O}_{\mu}((sh)^{2}) t_{r} (|D_{h}^{3} v|^{2} )d t.
\end{align}
Combining \cref{eqI43e1,eqI43e2,eqI43e3,eqI43e4}, for $ \lambda h (\delta T^{2})^{-1}\leq 1$ and $ \lambda \geq T^{2} $, we obtain 
\begin{align}
	\label{eqI43}
	I_{43} 
	\geq 
	48 \mathbb{E} \int_{\overline{Q}} s^{3} \mu^{4} \varphi^{3} |\partial_{x} \psi|^{4} |D_{h}^{2} v|^{2} d t 
	- \mathcal{A}_{1} - \mathcal{B}
	.
\end{align}

From \cref{eqdprpr,eqprpr,eqAu2leq,eqAIgbp}, thanks to Cauchy-Schwarz inequality, for $ \lambda h (\delta T^{2})^{-1}\leq 1$ and $ \lambda \geq T^{2} $, we have 
\begin{align}
	\notag
	\label{eqI44}
	I_{44} 
	& = 12 \mathbb{E} \int_{Q} A_{h} D_{h} (r A_{h}^{2} D_{h}^{2} \rho) A_{h} D_{h} (r A_{h}  D_{h}^{3} \rho) A_{h} D_{h} v v d t 
	\\ \notag
	&
	\geq
	-\mathbb{E} \int_{Q^{*}} s^{4} \mathcal{O}_{\mu}(1) |D_{h} v|^{2} d t 
	- \mathbb{E} \int_{Q} s^{6} \mathcal{O}_{\mu}(1) |v|^{2} d t 
	\\
	& \geq 
	- \mathcal{A}_{1}
	.
\end{align}

From \cref{eqI,eqDuv,eqAuv,eqA2u,eqdprpr,eqprpr}, letting 
\begin{equation}
	\label{eqQ}
	\begin{cases}
		\begin{aligned}
			& q_{1} = A_{h}  D_{h} (r A_{h}^{2} D_{h}^{2} \rho),
			\quad 
			q_{2} = A_{h}D_{h} (rD_{h}^{4} \rho),
			\quad 
			q_{3} = D_{h}^{2} (r D_{h}^{4} \rho),
			\\
			\quad 
			& q_{4} = D_{h}^{2} (r A_{h}^{4} \rho),
			\quad 
			q_{5} = A_{h} D_{h} (r A_{h}^{4} \rho),
		\end{aligned}
	\end{cases}
\end{equation}
for $ \lambda h (\delta T^{2})^{-1}\leq 1$, $h\le 1$, we obtain 
\begin{align}
	\notag
	\label{eqPhiE}
	\Phi(v) 
	&  =
	3 h^{2} D_{h}^{2} (q_{1} A_{h} D_{h} v ) 
	+ \frac{h^{2}}{2} A_{h}^{2} (  q_{2} A_{h}   D_{h} v )
	+ \frac{h^{2}}{4} q_{3} A_{h}^{2} v 
	- q_{4} D_{h}^{2} v 
	+ 2  A_{h} D_{h}  (  q_{5} D_{h}^{2} v ) 
	\\ \notag
	& =
	3 h^{2}  (
		 D_{h}^{2} q_{1} A_{h}^{3} D_{h} v   
		+   2 A_{h} D_{h}  q_{1} A_{h}^{2} D_{h}^{2} v   
		+    A_{h}^{2} q_{1} A_{h} D_{h}^{3} v   
	)
	+ \frac{h^{2}}{2} q_{2} A_{h}D_{h} v 
	\\ \notag
	& \quad 
	+ \frac{h^{4}}{8}  (
		D_{h}^{2} q_{2} A_{h}^{3} D_{h} v   
	   +  2  A_{h} D_{h}  q_{2} A_{h}^{2} D_{h}^{2} v   
	   +    A_{h}^{2} q_{2} A_{h} D_{h}^{3} v   
   	)
	+ \frac{h^{2}}{4} q_{3} A_{h}^{2} v 
	- q_{4} D_{h}^{2} v 
	\\ \notag
	& \quad
	+ 2 A_{h} D_{h} q_{5} A_{h}^{2} D_{h}^{2} v 
	+ \frac{h^{2}}{2} D_{h}^{2} q_{5} A_{h} D_{h}^{3} v 
	+ 2 A_{h} (A_{h} q_{5}   D_{h}^{3} v )
	\\ \notag
	& \geq 
	-\mathcal{O}_{\mu}((sh)^{2}) (|A_{h}^{3} D_{h} v| + |A_{h}^{2} D_{h}^{2} v| + |A_{h} D_{h}^{3} v|+ |D_{h}^{2} v| )
	- s^{2} \mathcal{O}_{\mu}((sh)^{2}) (|A_{h} D_{h} v|  + |A_{h}^{2} v|)
	\\
	& \quad 
 + 2 A_{h} (A_{h} q_{5}   D_{h}^{3} v ) 
	.
\end{align}

From \cref{eqQ,eqPhiE,eqdprpr,eqprpr,eqAu2leq,eqAIgbp}, thanks to Cauchy-Schwarz inequality, for $ \lambda h (\delta T^{2})^{-1}\leq 1$, $ h \leq 1 $ and $ \lambda \geq T^{2} $, we have 
\begin{align}
	\notag
	\label{eqI52}
	I_{52} 
	& =
	4 \mathbb{E} \int_{Q} r A_{h} D_{h}^{3}  \rho  A_{h}^{3} D_{h} v \Phi(v) d t
	\\ \notag
	& \geq 
	-\mathbb{E} \int_{Q} s^{5} \mathcal{O}_{\mu}((sh)^{2}) |v|^{2} d t
	- \mathbb{E} \int_{Q^{*}} [ 
		(s^{3} + s^{5}) \mathcal{O}_{\mu}((sh)^{2}) 
		+ s^{4} \mathcal{O}_{\mu}(1) 
		+ s^{5} \mu^{5} \varphi^{5}	 \mathcal{O}(1)
	] |D_{h} v|^{2} d t 
	\\ \notag
	& \quad 
	- \mathbb{E} \int_{\overline{Q}} [ s^{3}  \mathcal{O}_{\mu}((sh)^{2}) + s^{2} \mathcal{O}_{\mu}(1) ] |D_{h}^{2} v|^{2} d t 
	\\ \notag
	& \quad 
	- \mathbb{E} \int_{Q^{*}} [   \mathcal{O}_{\mu}(1) + s  \mathcal{O}_{\mu}((sh)^{2})+ s \mu \varphi \mathcal{O}(1)] |D_{h}^{3} v |^{2} d t
	\\
	& \geq 
	- \mathcal{A}_{1}  
	,
\end{align}
and
\begin{align}
	\notag
	\label{eqI53}
	I_{53}
	& = 
	4 \mathbb{E} \int_{Q} r A_{h}^{3} D_{h}  \rho  A_{h} D_{h}^{3}  v \Phi(v) d t 
	\\ \notag
	& \geq 
	-\mathbb{E} \int_{Q} s^{3} \mathcal{O}_{\mu}((sh)^{2}) |v|^{2} d t
	- \mathbb{E} \int_{Q^{*}} [ s \mathcal{O}_{\mu}((sh)^{2}) + s^{3} \mathcal{O}_{\mu}((sh)^{2}) ] |D_{h} v|^{2} d t 
	\\ \notag
	& \quad 
	- \mathbb{E} \int_{\overline{Q}}  s \mathcal{O}_{\mu}((sh)^{2})  |D_{h}^{2} v|^{2} d t 
	- \mathbb{E} \int_{Q^{*}} [  \mathcal{O}_{\mu}(1) + s  \mathcal{O}_{\mu}((sh)^{2})+ s \mu \varphi \mathcal{O}(1)] |D_{h}^{3} v |^{2} d t
	\\
	& \geq 
	- \mathcal{A}_{1}  
	.
\end{align}
Further, 
\begin{align}
	\notag
	\label{eqI54}
	I_{54}
	& = 
	2 \mathbb{E} \int_{Q}  A_{h} D_{h} (r A_{h} D_{h}^{3} \rho)v \Phi(v) d t  
	\\ \notag
	& \geq 
	-\mathbb{E} \int_{Q}  (s^{5} + s^{6} )  \mathcal{O}_{\mu} (1)  |v|^{2} d t
	- \mathbb{E} \int_{Q^{*}}  ( \mathcal{O}_{\mu}(1) +s ^{4} \mathcal{O}_{\mu}(1) ) |D_{h} v|^{2} d t 
	 \\ \notag
	& \quad 
		- \mathbb{E} \int_{\overline{Q}}   \mathcal{O}_{\mu}(1) ] |D_{h}^{2} v|^{2} d t 
	- \mathbb{E} \int_{Q^{*}}   \mathcal{O}_{\mu}(1)   |D_{h}^{3} v |^{2} d t
	\\
	& \geq 
	- \mathcal{A}_{1} 
	.
\end{align}

Combining \cref{eqI12,eqI13,eqI14,eqI22,eqI23,eqI24,eqI32,eqI33,eqI34,eqI42,eqI43,eqI44,eqI52,eqI53,eqI54}, for $ \lambda h (\delta T^{2})^{-1}\leq 1$, $ h \leq 1 $ and $ \lambda \geq T^{2} $, we have 
\begin{align}
	\notag
	\label{eqIij}
	\sum_{i=1}^{5}\sum_{j=2}^{4} I_{i j} 
	& \geq 
	8\mathbb{E} \int_{Q} s^{7} \mu^{8} \varphi^{7} |\partial_{x} \psi|^{8} | v |^{2} d t 
	+ 18 \mathbb{E} \int_{Q^{*}} s^{5} \mu^{6} \varphi^{5} |\partial_{x} \psi|^{6} |D_{h} v |^{2} d t 
	\\ \notag
	& \quad 
	+ 60 \mathbb{E} \int_{\overline{Q}} s^{3} \mu^{4} \varphi^{3} |\partial_{x} \psi|^{4} |D_{h}^{2} v |^{2} d t 
	+ 2 \mathbb{E} \int_{Q^{*}} s \mu^{2} \varphi  |\partial_{x} \psi|^{2} |D_{h}^{3} v |^{2} d t 
	\\ \notag
	& \quad 
	+ 12 \mathbb{E} \int_{\partial Q} s^{5} \mu^{5} \varphi^{5} |\partial_{x} \psi|^{5} t_{r} (|D_{h} v|^{2})  d t 
	+ 10 \mathbb{E} \int_{\partial Q} s^{3} \mu^{3} \varphi^{3} |\partial_{x} \psi|^{3}  |D_{h}^{2} v|^{2}  d t 
	\\
	& \quad 
	+ 2 \mathbb{E} \int_{\partial Q} s \mu \varphi |\partial_{x} \psi| t_{r} (|D_{h}^{3} v|^{2})  d t 
	- \mathcal{A}_{1} - \mathcal{B}
	.
\end{align}

\emph{Step 4.} In this step, we return to  the  variable $ w $.

From \cref{eqI,eqtCareleman,eqdpr2pr,eqAu2leq,eqAIgbp,eqPhiE}, for $ \lambda h (\delta T^{2})^{-1}\leq 1$ and $ \lambda \geq T + T^{2} $, we have 
\begin{align}
	\notag
	\label{eqI3}
	\mathbb{E} \int_{Q} |I_{3} |^{2} d t 
	& \leq 
	\mathbb{E} \int_{Q}  (s^{4} + s^{6}) \mathcal{O}_{\mu}(1) |v|^{2} d t 
	+ \mathbb{E} \int_{Q}  (s^{4} \mathcal{O}_{\mu}(1)+\mathcal{O}_{\mu}(1) ) |D_{h} v|^{2} d t  
	\\\notag
	& \quad 	
	+ \mathbb{E} \int_{Q}   \mathcal{O}_{\mu}(1) |D_{h}^{2} v|^{2} d t+ \mathbb{E} \int_{Q}   \mathcal{O}_{\mu}(1) |D_{h}^{3} v|^{2} d t 
	\\
	& \leq 
	\mathcal{A}_{1}
	.
\end{align}

Combining \cref{eqIi1,eqIij,eqC5,eqDefIij,eqI3,eqPsi}, for all $ \mu \geq 1 $, there exists a positive constant $ \varepsilon_{1} \leq 1 $ and a sufficiently large $ \lambda_{1} \geq 1  $, such that for $ \lambda \geq \lambda_{1} (T + T^{2})  $, $ h\leq h'_{1} $,  and $ \lambda h (\delta T^{2})^{-1}\leq \varepsilon_{1} $, we have 
\begin{align}
	\notag
	\label{eqRwGe1}
	& C \mathbb{E} \int_{Q} r^{2} | f|^{2} d t  
	+ C(\mu) \mathbb{E} \int_{Q} s^{4} r^{2} |g|^{2} d t 
	+ C \mathbb{E}  \int_{0}^{T} \int_{\mathcal{M}\cap G_{2}} s^{7} \mu^{8} \varphi^{7}  | v |^{2} d t 
	\\ \notag
	&  
	+   C \Big[ \mathbb{E}  \int_{0}^{T} \int_{\mathcal{M}^{*}\cap G_{2}} s^{5} \mu^{6} \varphi^{5}  |D_{h} v |^{2} d t 
	+    \mathbb{E}  \int_{0}^{T} \int_{\overline{\mathcal{M}}\cap G_{2}} s^{3} \mu^{4} \varphi^{3}  |D_{h}^{2} v |^{2} d t 
	+   \mathbb{E} \int_{0}^{T}  \int_{\mathcal{M}^{*}\cap G_{2}} s \mu^{2} \varphi   |D_{h}^{3} v |^{2} d t 
	\Big]
	\\ \notag
	&  \geq  
	\mathbb{E} \int_{Q} s^{7} \mu^{8} \varphi^{7}  |  v |^{2} d t 
	+   \mathbb{E} \int_{Q^{*}} s^{5} \mu^{6} \varphi^{5}  |D_{h} v |^{2} d t 
	+  \mathbb{E} \int_{\overline{Q}} s^{3} \mu^{4} \varphi^{3}  |D_{h}^{2} v |^{2} d t 
	+  \mathbb{E} \int_{Q^{*}} s \mu^{2} \varphi   |D_{h}^{3} v |^{2} d t 
	\\  
	& \quad 
	+ \mathbb{E} \int_{\overline{Q}} r^{2} |D_{h}^{2} g |^{2} d t 
	- \mathcal{A}_{1}  
	- \mathcal{A}_{3}  
	.
\end{align}

Noting $ w = \rho v $, combining \cref{eqDuv,eqdpr2pr,eqpr2pr}, for  $ \lambda h (\delta T^{2})^{-1}\leq 1$, we obtain 
\begin{align}
	\notag
	\label{eqDwe1}
	& \mathbb{E} \int_{Q^{*}} s^{5} \mu^{6} \varphi^{5} r^{2}  |D_{h} w |^{2} d t 
	 =
	\mathbb{E} \int_{Q^{*}} s^{5} \mu^{6} \varphi^{5} r^{2}  |D_{h} (\rho v) |^{2} d t 
	\\ \notag
	& \leq 
	2 \mathbb{E} \int_{Q^{*}} s^{5} \mu^{6} \varphi^{5} r^{2}  ( |D_{h} \rho A_{h} v|^{2} + |A_{h} \rho  D_{h} v|^{2}) d t 
	\\ \notag
	& \leq 
	C \mathbb{E} \int_{Q^{*}} s^{7} \mu^{8} \varphi^{7}  |A_{h} v|^{2} d t 
	+ 2 \mathbb{E} \int_{Q^{*}} s^{5} \mu^{6} \varphi^{5}  |D_{h} v|^{2} d t 
	+ \mathbb{E} \int_{Q}  [s^{6} \mathcal{O}_{\mu}(1) + s^{7} \mathcal{O}_{\mu}((sh)^{2})] | v|^{2} d t 
	\\
	& \quad 
	+ \mathbb{E} \int_{Q^{*}} s^{5} \mathcal{O}_{\mu}((sh)^{2}) | D_{h} v|^{2} d t 
	.
\end{align}
From \cref{eqAu2leq,eqAIgbp,eqAlpha,propADf,eqvBoundary}, we have 
\begin{align}
	\notag
	\label{eqDwe2}
	\mathbb{E} \int_{Q^{*}} s^{7} \mu^{8} \varphi^{7}  |A_{h} v|^{2} d t 
	& \leq
	\mathbb{E} \int_{Q^{*}} s^{7} \mu^{8} \varphi^{7}  A_{h} ( v^{2} )d t 
	\\ \notag
	& = 
	\mathbb{E} \int_{Q} s^{7} \mu^{8} A_{h} ( \varphi^{7})   v^{2} d t 
	\\
	& \leq 
	\mathbb{E} \int_{Q} s^{7} \mu^{8} \varphi^{7}  |  v|^{2} d t 
	+ \mathbb{E} \int_{Q} s^{5} \mathcal{O}_{\mu}((sh)^{2}) |  v|^{2} d t 
	.
\end{align}
Combining \cref{eqDwe1,eqDwe2}, for $ \lambda h (\delta T^{2})^{-1}\leq 1$ and $ \lambda \geq T^{2} $, it holds that 
\begin{align}
	\label{eqDw}
	\mathbb{E} \int_{Q^{*}} s^{5} \mu^{6} \varphi^{5} r^{2}  |D_{h} w |^{2} d t 
	\leq 
	2 \mathbb{E} \int_{Q} s^{7} \mu^{8} \varphi^{7}  |  v|^{2} d t 
	+ 2 \mathbb{E} \int_{Q^{*}} s^{5} \mu^{6} \varphi^{5}  |D_{h} v|^{2} d t 
	+ \mathcal{A}_{1}
	.
\end{align}
Since
\begin{align*}
	r^{2} |D_{h}^{2} w|^{2} 
	=  
	r^{2} |D_{h}^{2} (\rho v)|^{2} 
	\leq C r^{2} (
		|D_{h}^{2} \rho A_{h}^{2} v|^{2}
		+ |A_{h} D_{h}  \rho A_{h} D_{h} v|^{2}
		+ |A_{h}^{2} \rho  D_{h}^{2} v|^{2}
	)
	,
\end{align*}
similar to \cref{eqDw}, for $ \lambda h (\delta T^{2})^{-1}\leq 1$ and $ \lambda \geq T^{2} $, we have 
\begin{align}
	\notag
	\label{eqRwGe2}
	&\mathbb{E} \int_{\overline{Q}} s^{3} \mu^{4} \varphi^{3} r^{2}  |D_{h}^{2} w |^{2} d t 
	\\ \notag 
	& \leq 
	C \Big [
		\mathbb{E} \int_{Q} s^{7} \mu^{8} \varphi^{7}  |  v |^{2} d t 
		+   \mathbb{E} \int_{Q^{*}} s^{5} \mu^{6} \varphi^{5}  |D_{h} v |^{2} d t 
		+  \mathbb{E} \int_{\overline{Q}} s^{3} \mu^{4} \varphi^{3}  |D_{h}^{2} v |^{2} d t 
	\Big]
	\\ \notag
	& \quad 
	+  \mathbb{E} \int_{\overline{Q}} s^{3} \mathcal{O}_{\mu}((sh)^{2})  |D_{h}^{2} v |^{2} d t 
	+   \mathbb{E} \int_{Q^{*}} [s^{5} \mathcal{O}_{\mu}((sh)^{2}) + s^{ 4}\mathcal{O}_{\mu}(1) ) |D_{h} v |^{2} d t 
	\\ \notag 
	& \quad 
	+   \mathbb{E} \int_{Q} [s^{7} \mathcal{O}_{\mu}((sh)^{2}) +  s^{6}\mathcal{O}_{\mu}(1) ) | v |^{2} d t 
	\\ 
	& \leq 
	C \Big [
		\mathbb{E} \int_{Q} s^{7} \mu^{8} \varphi^{7}  |  v |^{2} d t 
		+   \mathbb{E} \int_{Q^{*}} s^{5} \mu^{6} \varphi^{5}  |D_{h} v |^{2} d t 
		+  \mathbb{E} \int_{\overline{Q}} s^{3} \mu^{4} \varphi^{3}  |D_{h}^{2} v |^{2} d t 
	\Big]
	+ \mathcal{A}_{1}
	.
\end{align}
Noting that
\begin{align*}
	r^{2} |D_{h}^{3} w|^{2} 
	=  
	r^{2} |D_{h}^{3} (\rho v)|^{2} 
	\leq C r^{2} (
		|D_{h}^{3} \rho A_{h}^{3} v|^{2}
		+ |A_{h} D_{h}^{2}  \rho A_{h}^{2} D_{h} v|^{2}
		+ |A_{h}^{2} D_{h}  \rho A_{h} D_{h}^{2} v|^{2}
		+ |A_{h}^{3} \rho  D_{h}^{3} v|^{2}
	)
	,
\end{align*}
similar to \cref{eqDw}, for $ \lambda h (\delta T^{2})^{-1}\leq 1$ and $ \lambda \geq T^{2} $, we have 
\begin{align}
	\notag
	\label{eqRwGe3}
	&\mathbb{E} \int_{Q^{*}} s \mu^{2} \varphi r^{2}  |D_{h}^{3} w |^{2} d t 
	\\ \notag 
	& \leq 
	C \Big [
		\mathbb{E} \int_{Q} s^{7} \mu^{8} \varphi^{7}  |  v |^{2} d t 
		+   \mathbb{E} \int_{Q^{*}} s^{5} \mu^{6} \varphi^{5}  |D_{h} v |^{2} d t 
		+  \mathbb{E} \int_{\overline{Q}} s^{3} \mu^{4} \varphi^{3}  |D_{h}^{2} v |^{2} d t 
		+  \mathbb{E} \int_{Q^{*}} s  \mu^{2} \varphi   |D_{h}^{3} v |^{2} d t 
	\Big]
	\\ \notag
	& \quad 
	+  \mathbb{E} \int_{Q^{*}} s  \mathcal{O}_{\mu}((sh)^{2})  |D_{h}^{3} v |^{2} d t 
	+  \mathbb{E} \int_{\overline{Q}} [s^{3} \mathcal{O}_{\mu}((sh)^{2}) + s^{2} \mathcal{O}_{\mu}(1)]  |D_{h}^{2} v |^{2} d t 
	\\ \notag
	& \quad 
	+   \mathbb{E} \int_{Q^{*}} [s^{5} \mathcal{O}_{\mu}((sh)^{2}) + s^{ 4}\mathcal{O}_{\mu}(1) ) |D_{h} v |^{2} d t 
	+   \mathbb{E} \int_{Q} [s^{7} \mathcal{O}_{\mu}((sh)^{2}) +  s^{6}\mathcal{O}_{\mu}(1) ) | v |^{2} d t 
	\\ \notag 
	& \leq 
	C \Big [
		\mathbb{E} \int_{Q} s^{7} \mu^{8} \varphi^{7}  |  v |^{2} d t 
		+   \mathbb{E} \int_{Q^{*}} s^{5} \mu^{6} \varphi^{5}  |D_{h} v |^{2} d t 
		+  \mathbb{E} \int_{\overline{Q}} s^{3} \mu^{4} \varphi^{3}  |D_{h}^{2} v |^{2} d t 
		+  \mathbb{E} \int_{Q^{*}} s  \mu^{2} \varphi   |D_{h}^{3} v |^{2} d t 
		\Big]
	\\
	& \quad 
	+ \mathcal{A}_{1}
	.
\end{align}

Hence, combining  \cref{eqDw,eqRwGe1,eqRwGe2,eqRwGe3}, for $ \mu \geq 1 $, $ \lambda \geq \lambda_{1} (T + T^{2})  $, $ h\leq h_{1} $,  and $ \lambda h (\delta T^{2})^{-1}\leq \varepsilon_{1} $, we have 
\begin{align}
	\notag
	\label{eqRwLe1}
	& C \mathbb{E} \int_{Q} r^{2} | f|^{2} d t  
	+ C(\mu) \mathbb{E} \int_{Q} s^{4} r^{2} |g|^{2} d t 
	+ C \mathbb{E}  \int_{0}^{T} \int_{\mathcal{M}\cap G_{2}} s^{7} \mu^{8} \varphi^{7}  | v |^{2} d t 
	\\ \notag
	&  
	+   C \Big[ \mathbb{E}  \int_{0}^{T} \int_{\mathcal{M}^{*}\cap G_{2}} s^{5} \mu^{6} \varphi^{5}  |D_{h} v |^{2} d t 
	+    \mathbb{E}  \int_{0}^{T} \int_{\overline{\mathcal{M}}\cap G_{2}} s^{3} \mu^{4} \varphi^{3}  |D_{h}^{2} v |^{2} d t 
	+   \mathbb{E} \int_{0}^{T}  \int_{\mathcal{M}^{*}\cap G_{2}} s \mu^{2} \varphi   |D_{h}^{3} v |^{2} d t 
	\Big]
	\\ \notag
	&  \geq  
	\mathbb{E} \int_{Q} s^{7} \mu^{8} \varphi^{7}  r^2|  w |^{2} d t 
	+   \mathbb{E} \int_{Q^{*}} s^{5} \mu^{6} \varphi^{5}  r^2|D_{h} w |^{2} d t 
	+  \mathbb{E} \int_{\overline{Q}} s^{3} \mu^{4} \varphi^{3}  r^2|D_{h}^{2} w |^{2} d t 
	\\  
	& \quad 
	+  \mathbb{E} \int_{Q^{*}} s \mu^{2} \varphi  r^2 |D_{h}^{3} w |^{2} d t + \mathbb{E} \int_{\overline{Q}} r^{2} |D_{h}^{2} g |^{2} d t 
	- \mathcal{A}_{1}  
	- \mathcal{A}_{3}  
	.
\end{align}

From \cref{eqDuv,eqdprpr}, for $ \lambda h (\delta T^{2})^{-1}\leq 1$, we have
\begin{align*}
	 r D_{h} w =  r D_{h} \rho A_{h} v  + r A_{h} \rho D_{h} v 
	 = r D_{h} \rho A_{h} v  + D_{h} v + \mathcal{O}_{\mu} ((sh)^{2}) D_{h} v 
	 ,
\end{align*}
which implies
\begin{align}
	\notag
	\label{eqRwLe2}
	& \mathbb{E} \int_{0}^{T} \int_{M^{*} \cap G_{2}} s^{5} \mu^{6} \varphi^{5} |D_{h} v|^{2} d t 
	\\ \notag
	& \leq 
	C \mathbb{E} \int_{0}^{T} \int_{M^{*} \cap G_{2}} s^{5} \mu^{6} \varphi^{5} (r D_{h} \rho)^{2}|A_{h} v|^{2} d t 
	+ C  \mathbb{E} \int_{0}^{T} \int_{M^{*} \cap G_{2}} s^{5} \mu^{6} \varphi^{5} r^{2} |D_{h} w|^{2} d t 
	\\ \notag
	& \quad 
	+ \mathbb{E} \int_{Q^{*}} s^{5} \mathcal{O}_{\mu}((sh)^{2}) |D_{h} v|^{2} d t
	\\ \notag
	& \leq 
	C \mathbb{E} \int_{0}^{T} \int_{M  \cap G_{1}} s^{7} \mu^{8} \varphi^{7} r^{2} | w|^{2} d t 
	+ C  \mathbb{E} \int_{0}^{T} \int_{M^{*} \cap G_{1}} s^{5} \mu^{6} \varphi^{5} r^{2} |D_{h} w|^{2} d t 
	\\ \notag
	& \quad 
	+ \mathbb{E} \int_{Q^{*}} s^{5} \mathcal{O}_{\mu}((sh)^{2}) |D_{h} v|^{2} d t
	+ \mathbb{E} \int_{Q } [s^{7} \mathcal{O}_{\mu}((sh)^{2}) + s^{6} \mathcal{O}_{\mu}(1) ]| v|^{2} d t
	\\
	& \leq 
	C \mathbb{E} \int_{0}^{T} \int_{M  \cap G_{1}} s^{7} \mu^{8} \varphi^{7} r^{2} | w|^{2} d t 
	+ C  \mathbb{E} \int_{0}^{T} \int_{M^{*} \cap G_{1}} s^{5} \mu^{6} \varphi^{5} r^{2} |D_{h} w|^{2} d t  
	+ \mathcal{A}_{1} 
	, 
\end{align}
where we use $ |A_{h} v|^{2} \leq C (|\tau_{+}v|^{2} + |\tau_{-} v|^{2}) $.

Similarly, for $ \lambda h (\delta T^{2})^{-1}\leq 1$, we have
\begin{align*}
	D_{h}^{2} v 
	=
	r D_{h}^{2} w 
	- r D_{h}^{2} \rho A_{h}^{2} v 
	- 2 r A_{h} D_{h} \rho A_{h} D_{h} v 
	- \mathcal{O}_{\mu} ((sh)^{2}) D_{h}^{2} v 
	,
\end{align*}
which implies
\begin{align}
	\notag
	\label{eqRwLe3}
	& \mathbb{E} \int_{0}^{T} \int_{\overline{\mathcal{M}} \cap G_{2}} s^{3} \mu^{4} \varphi^{3} |D_{h}^{2} v|^{2} d t 
	\\ \notag
	& \leq 
	C \mathbb{E} \int_{0}^{T} \int_{M  \cap G_{1}} s^{7} \mu^{8} \varphi^{7} r^{2} | w|^{2} d t 
	+ C  \mathbb{E} \int_{0}^{T} \int_{M^{*} \cap G_{1}} s^{5} \mu^{6} \varphi^{5} r^{2} |D_{h} w|^{2} d t 
	\\
	& \quad 
	+ C \mathbb{E} \int_{0}^{T} \int_{\overline{\mathcal{M}}\cap G_{1}} s^{3} \mu^{4} \varphi^{2} r^{2} |D_{h}^{2} w|^{2} d t 
	+ \mathcal{A}_{1}
	.
\end{align}
For $ \lambda h (\delta T^{2})^{-1}\leq 1$, we also obtain 
\begin{align}
	\notag
	\label{eqRwLe4}
	& \mathbb{E} \int_{0}^{T} \int_{\mathcal{M}^{*} \cap G_{2}} s  \mu^{2} \varphi  |D_{h}^{3} v|^{2} d t 
	\\ \notag
	& \leq 
	C \mathbb{E} \int_{0}^{T} \int_{M  \cap G_{1}} s^{7} \mu^{8} \varphi^{7} r^{2} | w|^{2} d t 
	+ C  \mathbb{E} \int_{0}^{T} \int_{M^{*} \cap G_{1}} s^{5} \mu^{6} \varphi^{5} r^{2} |D_{h} w|^{2} d t 
	\\
	& \quad 
	+ C \mathbb{E} \int_{0}^{T} \int_{\overline{\mathcal{M}}\cap G_{1}} s^{3} \mu^{4} \varphi^{2} r^{2} |D_{h}^{2} w|^{2} d t 
	+ C \mathbb{E} \int_{0}^{T} \int_{ M^{*} \cap G_{1}} s  \mu^{2} \varphi r^{2} |D_{h}^{3} w|^{2} d t 
	+ \mathcal{A}_{1}
	.
\end{align}

Noting that 
\begin{align}
	\label{eqDleq}
	|D_{h} v|^{2} \leq C h^{-2} (|\tau_{-} v|^{2} + |\tau_{+} v|^{2})
	,
\end{align}
and 
\begin{align*}
	|D_{h}^{2} v|^{2} \leq C h^{-4} (|\tau_{-}^{2} v|^{2} + |\tau_{+}^{2} v|^{2} + |v|^{2})
	,
\end{align*}
for all $ \mu \geq 1 $ and $ \lambda h (\delta T^{2})^{-1}\leq 1$, we have 
\begin{align}
	\label{eqA3}
	\mathcal{A}_{3} \leq C(\mu) h^{-4}  \Big( 
		\mathbb{E} \int_{\mathcal{M}} r^{2} |w|^{2} \Big|_{t=0}
		+ \mathbb{E} \int_{\mathcal{M}} r^{2} |w|^{2} \Big|_{t=T}
	\Big)
\end{align}

Combining \cref{eqRwLe1,eqRwLe2,eqRwLe3,eqRwLe4,eqA3}, there exists a positive constant $ \mu_{1} \geq 1 $ such that for all $ \mu \geq \mu_{1} $, one can find positive constants $ \varepsilon_{2} \leq \varepsilon_{1} $, $ h_{2} \leq h'_{1} $ and  $ \lambda_{2} \geq \lambda_{1}  $, such that for $ \lambda \geq \lambda_{2} (T + T^{2})  $, $ h\leq h_{2} $,  and $ \lambda h (\delta T^{2})^{-1}\leq \varepsilon_{2} $, we have 
\begin{align}
	\notag
	\label{eqFinale1}
	& C \mathbb{E} \int_{Q} r^{2} | f|^{2} d t  
	+ C(\mu) \mathbb{E} \int_{Q} s^{4} r^{2} |g|^{2} d t 
	+ C \mathbb{E}  \int_{0}^{T} \int_{\mathcal{M}\cap G_{1}} s^{7} \mu^{8} \varphi^{7}r^{2}  | w |^{2} d t 
	\\ \notag
	&  
	+   C   \mathbb{E}  \int_{0}^{T} \int_{\mathcal{M}^{*}\cap G_{1}} s^{5} \mu^{6} \varphi^{5} r^{2} |D_{h} w |^{2} d t 
	+    C  \mathbb{E}  \int_{0}^{T} \int_{\overline{\mathcal{M}}\cap G_{1}} s^{3} \mu^{4} \varphi^{3} r^{2} |D_{h}^{2} w |^{2} d t 
	\\ \notag
	& 
	+   C \mathbb{E} \int_{0}^{T}  \int_{\mathcal{M}^{*}\cap G_{1}} s \mu^{2} \varphi r^{2}  |D_{h}^{3} w |^{2} d t 
	+ C(\mu) h^{-4}  \Big( 
		\mathbb{E} \int_{\mathcal{M}} r^{2} |w|^{2} \Big|_{t=0}
		+ \mathbb{E} \int_{\mathcal{M}} r^{2} |w|^{2} \Big|_{t=T}
	\Big)
	\\ \notag
	&  \geq  
	\mathbb{E} \int_{Q} s^{7} \mu^{8} \varphi^{7}  |  w |^{2} d t 
	+   \mathbb{E} \int_{Q^{*}} s^{5} \mu^{6} \varphi^{5}  |D_{h} w |^{2} d t 
	+  \mathbb{E} \int_{\overline{Q}} s^{3} \mu^{4} \varphi^{3}  |D_{h}^{2} w |^{2} d t 
	+  \mathbb{E} \int_{Q^{*}} s \mu^{2} \varphi   |D_{h}^{3} w |^{2} d t 
	\\
	& \quad 
	+ \mathbb{E} \int_{\overline{Q}} r^{2} |D_{h}^{2} g |^{2} d t 
	.
\end{align}

\emph{Step 5.} In this step, we estimate the local terms on $ D_{h} w $, $ D_{h}^{2} w $ and $ D_{h}^{3} w $.

Choose a function $ \xi \in C_{0}^{\infty}(G_{0};[0,1]) $ such that $ \xi =1 $ in $ G_{1} $.
Thanks to It\^o's formula and \cref{eqW}, we have 
\begin{align*}
	d(s^{3} \varphi^{3} \xi^{4} r^{2} w^{2}) 
	& = 
	\varphi^{3} \xi^{4} \partial_{t}(s^{3} r^{2}) w^{2} d t 
	+ 2 s^{3} \varphi^{3} \xi^{4} r^{2} w d w 
	+ s^{3} \varphi^{3} \xi^{4} r^{2} (d w)^{2}
	.
\end{align*}
Hence, we get 
\begin{align} 
	\label{eqLocalD2we1}
	\notag
	& \mathbb{E} \int_{\mathcal{M} \cap G_{0}} s^{3} \varphi^{3} \xi^{4} r^{2} w^{2} \Big|_{t = 0}
	- \mathbb{E} \int_{\mathcal{M} \cap G_{0}} s^{3} \varphi^{3} \xi^{4} r^{2} w^{2} \Big|_{t = T}
	\\  \notag
	& =
	\mathbb{E} \int_{0}^{T} \int_{\mathcal{M} \cap G_{0}}\varphi^{3} \xi^{4} \partial_{t}(s^{3} r^{2}) w^{2} d t 
	+ 2 \mathbb{E} \int_{0}^{T} \int_{\mathcal{M} \cap G_{0}} s^{3} \varphi^{3} \xi^{4} r^{2} w (D_{h}^{4} w - f) d t   
	\\
	& \quad 
	+  \mathbb{E} \int_{0}^{T} \int_{\mathcal{M} \cap G_{0}} s^{3} \varphi^{3} \xi^{4} r^{2} |g|^{2}
	.
\end{align}

For all $ \mu \geq 1 $ and $ \lambda h (\delta T^{2})^{-1}\leq 1$, we have
\begin{align}
	\label{eqLocalD2we2}
	\mathbb{E} \int_{\mathcal{M} \cap G_{0}} s^{3} \varphi^{3} \xi^{4} r^{2} w^{2} \Big|_{t = 0}
	- \mathbb{E} \int_{\mathcal{M} \cap G_{0}} s^{3} \varphi^{3} \xi^{4} r^{2} w^{2} \Big|_{t = T}
	\leq C h^{-3} \Big( 
		\mathbb{E} \int_{\mathcal{M}} r^{2} |w|^{2} \Big|_{t=0}
		+ \mathbb{E} \int_{\mathcal{M}} r^{2} |w|^{2} \Big|_{t=T}
	\Big)
	.
\end{align}
For $ \lambda \geq T + T^{2} $, we obtain 
\begin{align}
	\label{eqLocalD2we3}
	\mathbb{E} \int_{0}^{T} \int_{\mathcal{M} \cap G_{0}}\varphi^{3} \xi^{4} \partial_{t}(s^{3} r^{2}) w^{2} d t 
	\leq 
	\mathbb{E} \int_{Q} s^{6} \mathcal{O}_{\mu} (1) r^{2} w^{2} d t 
	.
\end{align}
Thanks to Cauchy-Schwarz inequality, it holds that 
\begin{align}
	\label{eqLocalD2we4}
	\mathbb{E} \int_{0}^{T} \int_{\mathcal{M} \cap G_{0}} s^{3} \varphi^{3} \xi^{4} r^{2} w f d t  
	\leq 
	\mathbb{E} \int_{Q} s^{6} \mathcal{O}_{\mu} (1)r^{2} w^{2} d t 
	+ \mathbb{E} \int_{Q} r^{2} f ^{2} d t 
	.
\end{align}
From \cref{eqD2Ibp,eqDuv}, we have 
\begin{align}
	& \notag
	\label{eqLocalD2we5}
	2 \mathbb{E} \int_{0}^{T} \int_{\mathcal{M} \cap G_{0}} s^{3} \varphi^{3} \xi^{4} r^{2} w  D_{h}^{4} w   d t   
	\\\notag
	& = 
	2 \mathbb{E} \int_{0}^{T} \int_{\mathcal{M} \cap G_{0}} s^{3} D_{h}^{2} (\varphi^{3} \xi^{4} r^{2} w ) D_{h}^{2} w   d t   
	\\
	& = 
	2 \mathbb{E} \int_{0}^{T} \int_{\mathcal{M} \cap G_{0}} s^{3} [
		 D_{h}^{2} (\varphi^{3} \xi^{4} r^{2} ) A_{h}^{2} w 
		 + 2A_{h} D_{h} (\varphi^{3} \xi^{4} r^{2} ) A_{h} D_{h} w 
		 + A_{h}^{2}(\varphi^{3} \xi^{4} r^{2} ) D_{h}^{2} w 
	]D_{h}^{2} w   d t   
	.
\end{align}
Thanks to \cref{propADf,eqAu2leq,eqAIgbp} and Cauchy-Schwarz inequality, for all $ \varepsilon > 0 $ and $ \lambda h (\delta T^{2})^{-1}\leq 1$, it holds that 
\begin{align}
	\label{eqLocalD2we6}
	\notag
	& 
	2 \mathbb{E} \int_{0}^{T} \int_{\mathcal{M} \cap G_{0}} s^{3}  D_{h}^{2} (\varphi^{3} \xi^{4} r^{2} ) A_{h}^{2} w D_{h}^{2} w   d t   
	\\ \notag
	& \leq 
	\varepsilon \mathbb{E} \int_{0}^{T} \int_{\mathcal{M} \cap G_{0}} s^{3}    \varphi^{3} \xi^{4} r^{2}   |D_{h}^{2} w|^{2}   d t   
	+ C(\varepsilon) \mathbb{E} \int_{0}^{T} \int_{\mathcal{M} \cap G_{0}} s^{7} \mu^{4}    \varphi^{7} \xi^{4} r^{2}   |  w|^{2}   d t   
	\\
	& \quad 
	+ \mathbb{E} \int_{Q} [s^{2} \mathcal{O}_{\mu}(1) +  s^{3} \mathcal{O}_{\mu}((sh)^{2})] r^{2}  |D_{h}^{2} w|^{2}   d t   
	+ \mathbb{E} \int_{Q} [s^{6} \mathcal{O}_{\mu}(1) +  s^{7} \mathcal{O}_{\mu}((sh)^{2})]  r^{2} |w|^{2}   d t 
	,  
\end{align}
and
\begin{align}
	\label{eqLocalD2we7}
	\notag
	& 
	4 \mathbb{E} \int_{0}^{T} \int_{\mathcal{M} \cap G_{0}} s^{3} A_{h}  D_{h}  (\varphi^{3} \xi^{4} r^{2} ) A_{h} D_{h} w D_{h}^{2} w   d t   
	\\ \notag
	& \leq 
	\varepsilon \mathbb{E} \int_{0}^{T} \int_{\mathcal{M} \cap G_{0}} s^{3}    \varphi^{3} \xi^{4} r^{2}   |D_{h}^{2} w|^{2}   d t   
	+ C(\varepsilon) \mathbb{E} \int_{0}^{T} \int_{\mathcal{M}^{*} \cap G_{0}} s^{5} \mu^{2}    \varphi^{5} \xi^{4} r^{2}   | D_{h} w|^{2}   d t   
	\\
	& \quad 
	+ \mathbb{E} \int_{Q} [s^{2} \mathcal{O}_{\mu}(1) +  s^{3} \mathcal{O}_{\mu}((sh)^{2})] r^{2}  |D_{h}^{2} w|^{2}   d t   
	+ \mathbb{E} \int_{Q^{*}} [s^{4} \mathcal{O}_{\mu}(1) +  s^{5} \mathcal{O}_{\mu}((sh)^{2})] r^{2}  |D_{h} w|^{2}   d t 
	.
\end{align}
Similarly, for all $ \varepsilon > 0 $ and $ \lambda h (\delta T^{2})^{-1}\leq 1$, we get 
\begin{align}
	\label{eqLocalD2we8}
	& \notag
	2 \mathbb{E} \int_{0}^{T} \int_{\mathcal{M} \cap G_{0}} s^{3} A_{h}^{2}  (\varphi^{3} \xi^{4} r^{2} )  |D_{h}^{2} w|^{2}   d t   
	\\
	& \geq 
	2 \mathbb{E} \int_{0}^{T} \int_{\mathcal{M} \cap G_{0}} s^{3}  \varphi^{3} \xi^{4} r^{2}   |D_{h}^{2} w|^{2}   d t   
	- \mathbb{E} \int_{0}^{T} \int_{Q} s^{3} \mathcal{O}_{\mu}((sh)^{2}) r^{2} |D_{h}^{2} w|^{2}   d t   
	.
\end{align}
Combining \cref{eqLocalD2we1,eqLocalD2we2,eqLocalD2we3,eqLocalD2we4,eqLocalD2we5,eqLocalD2we6,eqLocalD2we7,eqLocalD2we8}, letting $ \varepsilon $ be sufficiently small, for $ \mu \geq 1 $, $ \lambda \geq T + T^{2} $, $ h \leq 1  $  and $ \lambda h (\delta T^{2})^{-1}\leq 1$, we obtain 
\begin{align}
	\notag
	\label{eqLocalD2w}
	& 
	\mathbb{E} \int_{0}^{T} \int_{\overline{\mathcal{M}} \cap G_{1}} s^{3} \mu^{4} \varphi^{3}   r^{2}   |D_{h}^{2} w|^{2}   d t   
	\\ \notag
	& \leq 
	\mathbb{E} \int_{0}^{T} \int_{\overline{\mathcal{M}} \cap G_{0}} s^{3} \mu^{4}  \varphi^{3} \xi^{4} r^{2}   |D_{h}^{2} w|^{2}   d t   
	\\ \notag
	& \leq 
	C \mathbb{E} \int_{0}^{T} \int_{\mathcal{M} \cap G_{0}} s^{7} \mu^{8}    \varphi^{7}  r^{2}   |  w|^{2}   d t  
	+ C \mathbb{E}  \int_{0}^{T} \int_{\mathcal{M} \cap G_{0}} s^{5} \mu^{6} \varphi^{5} \xi^{4} r^{2} |D_{h} w|^{2} d t  
	\\ \notag
	& \quad 
	+ \mathbb{E} \int_{Q} [s^{2} \mathcal{O}_{\mu}(1) +  s^{3} \mathcal{O}_{\mu}((sh)^{2})]  r^{2} |D_{h}^{2} w|^{2}   d t   
	+ \mathbb{E} \int_{Q} [s^{6} \mathcal{O}_{\mu}(1) +  s^{7} \mathcal{O}_{\mu}((sh)^{2})]  r^{2} |w|^{2}   d t 
	\\ \notag
	& \quad 
	+ \mathbb{E} \int_{Q^{*}} [s^{4} \mathcal{O}_{\mu}(1) +  s^{5} \mathcal{O}_{\mu}((sh)^{2})] r^{2}  |D_{h} w|^{2}   d t 
	+ C(\mu) \mathbb{E} \int_{Q} s^{4} r^{2} |g|^{2}
	\\
	& \quad 
	+ C(\mu) h^{-4} \Big( 
		\mathbb{E} \int_{\mathcal{M}} r^{2} |w|^{2} \Big|_{t=0}
		+ \mathbb{E} \int_{\mathcal{M}} r^{2} |w|^{2} \Big|_{t=T}
	\Big)
	.
\end{align}

From \cref{eqDIgbp,eqDuv,propADf}, thanks to Cauchy-Schwarz inequality, for all $ \varepsilon > 0 $ and $ \lambda h (\delta T^{2})^{-1}\leq 1$, we have 
\begin{align}
	\notag
	\label{eqLocalDwe1}
	&
	\mathbb{E}  \int_{0}^{T} \int_{\mathcal{M}^{*} \cap G_{0}} s^{5} \mu^{6} \varphi^{5} \xi^{4} r^{2} |D_{h} w|^{2} d t  
	\\ \notag
	& =
	- \mathbb{E}  \int_{0}^{T} \int_{\mathcal{M} \cap G_{0}} s^{5} D_{h}(\mu^{6} \varphi^{5} \xi^{4} r^{2} D_{h} w )  w d t  
	\\ \notag
	& =
	- \mathbb{E}  \int_{0}^{T} \int_{\mathcal{M} \cap G_{0}} s^{5} D_{h}(\mu^{6} \varphi^{5} \xi^{4} r^{2} ) A_{h} D_{h} w  w d t  
	- \mathbb{E}  \int_{0}^{T} \int_{\mathcal{M} \cap G_{0}} s^{5} A_{h}(\mu^{6} \varphi^{5} \xi^{4} r^{2} )   D_{h}^{2} w  w d t  
	\\ \notag
	& \leq 
	\varepsilon \mathbb{E}  \int_{0}^{T} \int_{\mathcal{M}^{*} \cap G_{0}} s^{5} \mu^{6} \varphi^{5} \xi^{4} r^{2} |D_{h} w|^{2} d t  
	+ \varepsilon \mathbb{E} \int_{0}^{T} \int_{ \mathcal{M}  \cap G_{0}} s^{3} \mu^{4}  \varphi^{3} \xi^{4} r^{2}   |D_{h}^{2} w|^{2}   d t  
	\\ \notag
	& \quad 
	+ C(\varepsilon) \mathbb{E} \int_{0}^{T} \int_{\mathcal{M} \cap G_{0}} s^{7} \mu^{8}    \varphi^{7}  r^{2}   |  w|^{2}   d t  
	+ \mathbb{E} \int_{Q^{*}} [s^{4} \mathcal{O}_{\mu}(1) +  s^{5} \mathcal{O}_{\mu}((sh)^{2})]  r^{2} |D_{h} w|^{2}   d t 
	\\
	& \quad 
	+ \mathbb{E} \int_{Q} [s^{6} \mathcal{O}_{\mu}(1) +  s^{7} \mathcal{O}_{\mu}((sh)^{2})]  r^{2} |w|^{2}   d t 
	.
\end{align}
From  \cref{eqLocalD2w,eqLocalDwe1}, letting $ \varepsilon $ be sufficiently small, for $ \mu \geq 1 $, $ \lambda \geq T + T^{2} $, $ h \leq 1  $  and $ \lambda h (\delta T^{2})^{-1}\leq 1$, we obtain 
\begin{align}
	\notag
	\label{eqLocalDw}
	&
	\mathbb{E}  \int_{0}^{T} \int_{\mathcal{M}^{*} \cap G_{1}} s^{5} \mu^{6} \varphi^{5}  r^{2} |D_{h} w|^{2} d t  
	\\ \notag
	& \leq
	\mathbb{E}  \int_{0}^{T} \int_{\mathcal{M}^{*} \cap G_{0}} s^{5} \mu^{6} \varphi^{5} \xi^{4} r^{2} |D_{h} w|^{2} d t  
	\\ \notag
	& \leq 
	C \mathbb{E} \int_{0}^{T} \int_{\mathcal{M} \cap G_{0}} s^{7} \mu^{8}    \varphi^{7}  r^{2}   |  w|^{2}   d t  
	+ \mathbb{E} \int_{Q^{*}} [s^{4} \mathcal{O}_{\mu}(1) +  s^{5} \mathcal{O}_{\mu}((sh)^{2})]  r^{2} |D_{h} w|^{2}   d t 
	\\ \notag
	& \quad 
	+ \mathbb{E} \int_{Q} [s^{6} \mathcal{O}_{\mu}(1) +  s^{7} \mathcal{O}_{\mu}((sh)^{2})]  r^{2} |w|^{2}   d t 
	+ \mathbb{E} \int_{Q} [s^{2} \mathcal{O}_{\mu}(1) +  s^{3} \mathcal{O}_{\mu}((sh)^{2})]  r^{2} |D_{h}^{2} w|^{2}   d t   
	\\
	& \quad 
	+ C(\mu) \mathbb{E} \int_{Q} s^{4} r^{2} |g|^{2}
	+ C(\mu) h^{-4} \Big( 
		\mathbb{E} \int_{\mathcal{M}} r^{2} |w|^{2} \Big|_{t=0}
		+ \mathbb{E} \int_{\mathcal{M}} r^{2} |w|^{2} \Big|_{t=T}
	\Big)
	.
\end{align}

Thanks to It\^o's formula and \cref{eqW}, we have 
\begin{align*}
	d(s  \varphi  \xi^{4} r^{2} |D_{h} w|^{2}) 
	& = 
	\varphi  \xi^{4} \partial_{t}(s  r^{2}) |D_{h}w|^{2} d t 
	+ 2 s \varphi \xi^{4} r^{2} D_{h }w d (D_{h} w) 
	+ s \varphi  \xi^{4} r^{2} |d (D_{h} w)|^{2}
	.
\end{align*}
Hence, we get 
\begin{align} 
	\label{eqLocalD3we1}
	\notag
	& \mathbb{E} \int_{\mathcal{M}^{*} \cap G_{0}} s  \varphi  \xi^{4} r^{2} |D_{h} w|^{2} \Big|_{t = 0}
	- \mathbb{E} \int_{\mathcal{M}^{*} \cap G_{0}} s  \varphi  \xi^{4} r^{2} |D_{h} w|^{2} \Big|_{t = T}
	\\  \notag
	& =
	\mathbb{E} \int_{0}^{T} \int_{\mathcal{M}^{*} \cap G_{0}}\varphi  \xi^{4} \partial_{t}(s  r^{2}) |D_{h}w|^{2} d t 
	+ 2 \mathbb{E} \int_{0}^{T} \int_{\mathcal{M}^{*} \cap G_{0}} s \varphi \xi^{4} r^{2} D_{h }w  (D_{h}^{5} w - D_{h} f) d t   
	\\
	& \quad 
	+  \mathbb{E} \int_{0}^{T} \int_{\mathcal{M}^{*} \cap G_{0}}  s \varphi  \xi^{4} r^{2}  |D_{h} g|^{2}
	.
\end{align}
From \cref{eqDleq}, for all $ \mu \geq 1 $ and $ \lambda h (\delta T^{2})^{-1}\leq 1$, we have
\begin{align}
	\label{eqLocalD3we2}
	\notag
	& \mathbb{E} \int_{\mathcal{M}^{*} \cap G_{0}} s  \varphi  \xi^{4} r^{2} |D_{h} w|^{2} \Big|_{t = 0}
	- \mathbb{E} \int_{\mathcal{M}^{*} \cap G_{0}} s  \varphi  \xi^{4} r^{2} |D_{h} w|^{2} \Big|_{t = T}
	\\
	& \leq C h^{-3} \Big( 
		\mathbb{E} \int_{\mathcal{M}} r^{2} |w|^{2} \Big|_{t=0}
		+ \mathbb{E} \int_{\mathcal{M}} r^{2} |w|^{2} \Big|_{t=T}
	\Big)
	.
\end{align}
For $ \lambda \geq T + T^{2} $, we obtain 
\begin{align}
	\label{eqLocalD3we3}
	\mathbb{E} \int_{0}^{T} \int_{\mathcal{M}^{*} \cap G_{0}}\varphi  \xi^{4} \partial_{t}(s  r^{2}) |D_{h}w|^{2} d t 
	\leq 
	\mathbb{E} \int_{Q} s^{4} \mathcal{O}_{\mu} (1) r^{2} |D_{h} w|^{2} d t 
	.
\end{align}
Similar to \cref{eqIi1e3}, for $ \varepsilon > 0 $  and $ \lambda h (\delta T^{2})^{-1}\leq 1$, we have 
\begin{align}
	\label{eqLocalD3we4}
	\mathbb{E} \int_{0}^{T} \int_{\mathcal{M}^{*} \cap G_{0}}  s \varphi  \xi^{4} r^{2}  |D_{h} g|^{2}
	\leq 
	\varepsilon \mathbb{E} \int_{\overline{Q}} r^{2} |D_{h}^{2} g |^{2} d t 
	+ C(\varepsilon,\mu)  \mathbb{E} \int_{Q} s^{3} r^{2} |g|^{2} d t 
	.
\end{align}
From \cref{eqDIgbp,eqDuv,propADf,eqAu2leq,eqAIgbp}, thanks to Cauchy-Schwarz inequality, for $ \lambda h (\delta T^{2})^{-1}\leq 1$, it holds that 
\begin{align}
	\notag
	\label{eqLocalD3we5}
	& -2\mathbb{E} \int_{0}^{T} \int_{\mathcal{M}^{*} \cap G_{0}} s  \varphi  \xi^{4} r^{2} D_{h} w D_{h} f d t  
	\\ \notag
	& = 
	\mathbb{E} \int_{0}^{T} \int_{\mathcal{M} \cap G_{0}} D_{h} (s  \varphi  \xi^{4} r^{2} D_{h} w) f d t  
	\\
	& \leq 
	\mathbb{E} \int_{Q} s^{2} \mathcal{O}_{\mu} (1)r^{2} |D_{h}^{2} w|^{2} d t 
	+ \mathbb{E} \int_{Q^{*}} s^{4} \mathcal{O}_{\mu} (1)r^{2} |D_{h} w|^{2} d t  
	+ \mathbb{E} \int_{Q} r^{2} f ^{2} d t 
	.
\end{align}
From \cref{eqD2Ibp,eqDuv}, we have 
\begin{align}
	\notag
	\label{eqLocalD3we6}
	& 
	2 \mathbb{E} \int_{0}^{T} \int_{\mathcal{M}^{*} \cap G_{0}} s  \varphi  \xi^{4} r^{2} D_{h} w D_{h}^{5} w d t  
	\\ \notag
	& =
	2 \mathbb{E} \int_{0}^{T} \int_{\mathcal{M}^{*} \cap G_{0}} D_{h}^{2}(s  \varphi  \xi^{4} r^{2}) A_{h}^{2} D_{h}  w  D_{h}^{3} w d t  
	+ 4 \mathbb{E} \int_{0}^{T} \int_{\mathcal{M}^{*} \cap G_{0}}A_{h} D_{h} (s  \varphi  \xi^{4} r^{2}) A_{h}  D_{h}^{2}  w  D_{h}^{3} w d t  
	\\
	& \quad 
	+ 2 \mathbb{E} \int_{0}^{T} \int_{\mathcal{M}^{*} \cap G_{0}} A_{h}^{2}   (s  \varphi  \xi^{4} r^{2})  |D_{h}^{3} w|^{2} d t 
	. 
\end{align}
Similar to \cref{eqLocalD2we6}, for  all $ \varepsilon > 0 $ and $ \lambda h (\delta T^{2})^{-1}\leq 1$, it holds that 
\begin{align}
	\notag
	\label{eqLocalD3we7}
	& 2 \mathbb{E} \int_{0}^{T} \int_{\mathcal{M}^{*} \cap G_{0}} D_{h}^{2}(s  \varphi  \xi^{4} r^{2}) A_{h}^{2} D_{h}  w  D_{h}^{3} w d t  
	+ 4 \mathbb{E} \int_{0}^{T} \int_{\mathcal{M}^{*} \cap G_{0}}A_{h} D_{h} (s  \varphi  \xi^{4} r^{2}) A_{h}  D_{h}^{2}  w  D_{h}^{3} w d t  
	\\ \notag
	& \leq 
	\varepsilon \mathbb{E} \int_{\mathcal{M}^{*} \cap G_{0}} s \varphi \xi^{4} r^{2} |D_{h}^{3} w |^{2} d t 
	+ C(\varepsilon) \mathbb{E} \int_{0}^{T} \int_{\mathcal{M}^{*} \cap G_{0}} s^{5} \mu^{4}    \varphi^{5} \xi^{4} r^{2}   | D_{h} w|^{2}   d t   
	\\ \notag
	& \quad 
	+ C(\varepsilon) \mathbb{E} \int_{0}^{T} \int_{\mathcal{M}  \cap G_{0}} s^{3} \mu^{2}    \varphi^{3} \xi^{4} r^{2}   | D_{h}^{2} w|^{2}   d t   
	+ \mathbb{E} \int_{Q^{*}} [  \mathcal{O}_{\mu}(1) +  s  \mathcal{O}_{\mu}((sh)^{2})]  r^{2} |D_{h}^{3} w|^{2}   d t   
	\\
	& \quad 
	+ \mathbb{E} \int_{Q} [s^{2} \mathcal{O}_{\mu}(1) +  s^{3} \mathcal{O}_{\mu}((sh)^{2})]   r^{2}|D_{h}^{2} w|^{2}   d t   
	+ \mathbb{E} \int_{Q^{*}} [s^{4} \mathcal{O}_{\mu}(1) +  s^{5} \mathcal{O}_{\mu}((sh)^{2})]   r^{2}|D_{h} w|^{2}   d t 
	,
\end{align}
and 
\begin{align}
	\notag
	\label{eqLocalD3we8}
	& 
	2 \mathbb{E} \int_{0}^{T} \int_{\mathcal{M}^{*} \cap G_{0}} A_{h}^{2}   (s  \varphi  \xi^{4} r^{2})  |D_{h}^{3} w|^{2} d t 
	\\
	& \geq
	2\mathbb{E} \int_{0}^{T} \int_{\mathcal{M}^{*} \cap G_{0}}  s  \varphi  \xi^{4} r^{2}  |D_{h}^{3} w|^{2} d t 
	- \mathbb{E} \int_{Q^{*}}   s  \mathcal{O}_{\mu}((sh)^{2}) r^{2}  |D_{h}^{3} w|^{2}   d t   
	.
\end{align}
Combining  \cref{eqLocalD3we1,eqLocalD3we2,eqLocalD3we3,eqLocalD3we4,eqLocalD3we5,eqLocalD3we6,eqLocalD3we7,eqLocalD3we8,eqLocalD2w,eqLocalDwe1,eqLocalDw}, letting $ \varepsilon $ be sufficiently small, for $ \mu \geq 1 $, $ \lambda \geq T + T^{2} $, $ h \leq 1  $  and $ \lambda h (\delta T^{2})^{-1}\leq 1$, we obtain 
\begin{align}
	\notag
	\label{eqLocalD3w}
	& 
	\mathbb{E} \int_{0}^{T} \int_{\mathcal{M}^{*} \cap G_{1}}  s  \mu^{2} \varphi   r^{2}  |D_{h}^{3} w|^{2} d t 
	\\ \notag
	& \leq 
	\mathbb{E} \int_{0}^{T} \int_{\mathcal{M}^{*} \cap G_{0}}  s  \mu^{2} \varphi  \xi^{4} r^{2}  |D_{h}^{3} w|^{2} d t 
	\\ \notag
	& \leq 
	C \mathbb{E} \int_{0}^{T} \int_{\mathcal{M} \cap G_{0}} s^{7} \mu^{8}    \varphi^{7}  r^{2}   |  w|^{2}   d t  
	+ \mathbb{E} \int_{Q^{*}} [s^{4} \mathcal{O}_{\mu}(1) +  s^{5} \mathcal{O}_{\mu}((sh)^{2})]  r^{2} |D_{h} w|^{2}   d t 
	\\ \notag
	& \quad 
	+ \mathbb{E} \int_{Q} [s^{6} \mathcal{O}_{\mu}(1) +  s^{7} \mathcal{O}_{\mu}((sh)^{2})]  r^{2} |w|^{2}   d t 
	+ \mathbb{E} \int_{Q} [s^{2} \mathcal{O}_{\mu}(1) +  s^{3} \mathcal{O}_{\mu}((sh)^{2})]  r^{2} |D_{h}^{2} w|^{2}   d t   
	\\ \notag
	& \quad 
	+ \mathbb{E} \int_{Q^{*}} [ \mathcal{O}_{\mu}(1) +  s  \mathcal{O}_{\mu}((sh)^{2})] r^{2}  |D_{h}^{3} w|^{2}   d t   
	+ C(\mu) \mathbb{E} \int_{Q} s^{4} r^{2} |g|^{2}
	\\
	& \quad 
	+ C(\mu) h^{-4} \Big( 
		\mathbb{E} \int_{\mathcal{M}} r^{2} |w|^{2} \Big|_{t=0}
		+ \mathbb{E} \int_{\mathcal{M}} r^{2} |w|^{2} \Big|_{t=T}
	\Big)
	.
\end{align}

Combining \cref{eqFinale1,eqLocalD2w,eqLocalDwe1,eqLocalDw,eqLocalD3w}, there exists a positive constant $ \mu_{0} \geq \mu_{1} $ such that for all $ \mu \geq \mu_{2} $, one can find positive constants $ \varepsilon_{0} \leq \varepsilon_{2} $, $ h_{0} \leq h_{2} $ and  $ \lambda_{0} \geq \lambda_{2}  $, such that for $ \lambda \geq \lambda_{0} (T + T^{2})  $, $ 0 < h\leq h_{0} $,  and $ \lambda h (\delta T^{2})^{-1}\leq \varepsilon_{0} $, we obtain \cref{eqFinalCarleman}.
\end{proof}

\section{Proof of the observability inequality}
\label{sec4}

This section is devoted to the proof of \cref{thmObservabilityEstiames}.

\begin{proof}[Proof of \cref{thmObservabilityEstiames}]
From \cref{thmCarleman}, choosing $ \mu = \mu_{0} $, for $ \lambda \geq \lambda_{0} (T+T^{2}) $, $ 0 < h \leq h_{0} $ and $\lambda h\left(\delta T^2\right)^{-1} \leq \varepsilon_0$, we have 
\begin{align*}
	\mathbb{E}\int_Q s^7 e^{2 s \phi}|z|^2 d t
	& \leq 
	C \Big(
		\mathbb{E} \int_{0}^{T}\int_{G_{0} \cap \mathcal{M}} s^{7} e^{2 s \phi} |z|^{2} d t
		+ \mathbb{E} \int_{Q} e^{2 s \phi} |a_{1} z + a_{2} Z|^{2} d t
		+ \mathbb{E} \int_{Q} s^{4} e^{2 s \phi} |Z|^{2} d t
	\\
	& \quad \quad \quad  
	+ h^{-4} \mathbb{E} \int_{\mathcal{M}} e^{2 s \phi}|z|^2 \Big|_{t=0}
	+ h^{-4} \mathbb{E} \int_{\mathcal{M}} e^{2 s \phi}|z|^2 \Big |_{t=T} \Big),
\end{align*}
From this, there exists $ \lambda_{1} \geq \lambda_{0} $ such that for a fixed
\begin{align}
	\label{eqFixedL}
	\lambda = \lambda_{1} (T + T^{2} + \mathcal{H}^{2/7} T^{2})
\end{align}
we obtain 
\begin{align}
	\notag
	\label{eqObe1}
	\mathbb{E}\int_Q \theta^7 e^{2 s \phi}|z|^2 d t
	& \leq 
	C \Big(
		\mathbb{E} \int_{0}^{T}\int_{G_{0} \cap \mathcal{M}} \theta^{7} e^{2 s \phi} |z|^{2} d t
		+ \mathbb{E} \int_{Q} \theta^{7} e^{2 s \phi} |Z|^{2} d t
		+ h^{-4} \mathbb{E} \int_{\mathcal{M}} e^{2 s \phi}|z|^2 \Big|_{t=0}
	\\
	& \quad \quad \quad  
	+ h^{-4} \mathbb{E} \int_{\mathcal{M}} e^{2 s \phi}|z|^2 \Big |_{t=T} \Big)\
	.
\end{align}
Noting that \cref{eqAlpha}, we have $ \theta(t) \geq \theta(T/2) \geq T^{-2} $ for $ t \in [0,T] $ and $ \theta(t) \leq \theta(T/4) \leq \frac{16}{3 T^{2}} $ for $ t\in [T/4, 3T/4] $. 
Then
\begin{align}
	\label{eqObe2}
	\mathbb{E}\int_Q \theta^7 e^{2 s \phi}|z|^2 d t
	 \geq
	\mathbb{E}\int_{\frac{T}{4}}^{\frac{3 T}{4}}\int_{\mathcal{M}} \theta^7 e^{2 s \phi}|z|^2 d t
	\geq 
	C e^{-C\lambda T^{-2}}\mathbb{E}\int_{\frac{T}{4}}^{\frac{3 T}{4}}  \int_{\mathcal{M}}|z|^2 d t
	.
\end{align}

From \cref{eqDiscreteEqAdjoint,eqD2Ibp}, thanks to Cauchy-Schwarz inequality, for any $ 0 \leq t_{1} \leq t_{2} \leq T $, we have 
\begin{align*}
	 \mathbb{E} \int_{\mathcal{M}} |z (t_{2})|^{2}
	- \mathbb{E} \int_{\mathcal{M}} |z (t_{1})|^{2}
	& =
	\mathbb{E} \int_{t_{1}}^{t_{2}} \int_{\mathcal{M}} d (z^{2})
	= 
	\mathbb{E} \int_{t_{1}}^{t_{2}} \int_{\mathcal{M}} [2 z d z + (d z)^{2}]
	\\
	& 
	= \mathbb{E} \int_{t_{1}}^{t_{2}} \int_{\mathcal{M}} [2 z (D_{h}^{4}z - a_{1} z - a_{2} Z) + Z^{2}] d t 
	\\
	& 
	= \mathbb{E} \int_{t_{1}}^{t_{2}} \int_{\mathcal{M}} [2 |D_{h}^{2} z|^{2} - 2 a_{1} z^{2} - a_{2} z Z  + Z^{2}] d t 
	\\
	& 
	\geq - C (1+ \mathcal{H}^{2}) \mathbb{E} \int_{t_{1}}^{t_{2}} \int_{\mathcal{M}} |z|^{2} d t 
	.
\end{align*}
Thus, in terms of Gronwall inequality, for $ 0 \leq t_{1} \leq t_{2} \leq T $, we have 
\begin{align}
	\label{eqObe3}
	\mathbb{E}   \int_{\mathcal{M}} |z(t_{1})|^{2} d t  
	\leq
	C e^{CT(1+ \mathcal{H}^{2})}\mathbb{E}   \int_{\mathcal{M}} |z(t_{2})|^{2} d t  
	.
\end{align}
From it, since $ \theta(T) = \theta(0) \geq \frac{1}{2} (\delta T^{2})^{-1} $, we obtain 
\begin{align}\notag
	\label{eqObe4}
	& 
	h^{-4} \mathbb{E} \int_{\mathcal{M}} e^{2 s \phi}|z|^2 \Big|_{t=0}
	+ h^{-4} \mathbb{E} \int_{\mathcal{M}} e^{2 s \phi}|z|^2 \Big |_{t=T}  
	\\ \notag
	& \leq 
	C h^{-4} e^{- C\lambda (\delta T^{2})^{-1}} \Big[ 
		\mathbb{E} \int_{\mathcal{M}} |z|^{2} \Big|_{t=0}
		+ \mathbb{E} \int_{\mathcal{M}} |z|^{2} \Big|_{t=T}
	\Big]
	\\ 
	& \leq 
	C h^{-4} e^{- C\lambda (\delta T^{2})^{-1} + C T (1+\mathcal{H}^{2})}  
	\mathbb{E} \int_{\mathcal{M}} |z|^{2} \Big|_{t=T}
	.
\end{align}
We have 
\begin{align}
	\notag
	\label{eqObe5}
	& \mathbb{E} \int_{0}^{T}\int_{G_{0} \cap \mathcal{M}} \theta^{7} e^{2 s \phi} |z|^{2} d t
	+ \mathbb{E} \int_{Q} \theta^{7} e^{2 s \phi} |Z|^{2} d t
	\\
	& \leq 
	C e^{- C \lambda T^{-2}} \Big(
		\mathbb{E} \int_{0}^{T}\int_{G_{0} \cap \mathcal{M}} |z|^{2} d t
	+ \mathbb{E} \int_{Q} |Z|^{2} d t
	\Big)
	.
\end{align}
Combining  \cref{eqObe1,eqObe2,eqObe3,eqObe4,eqObe5}, for $ 0 < h \leq h_{0} $ and $\lambda h\left(\delta T^2\right)^{-1} \leq \varepsilon_0$, it holds that 
\begin{align}
	\label{eqObe6}
	\mathbb{E} \int_{\mathcal{M}} |z(0)|^{2} 
	\leq 
	C e^{C \lambda T^{-2} + C T (1+\mathcal{H}^{2})}
	\Big(
		\mathbb{E} \int_{0}^{T}\int_{G_{0} \cap \mathcal{M}} |z|^{2} d t
		+ \mathbb{E} \int_{Q} |Z|^{2} d t
		+h^{-4} e^{- C\lambda (\delta T^{2})^{-1} } \int_{\mathcal{M}} |z|^{2} \Big|_{t=T}
	\Big)
\end{align}
Since we required $\lambda h\left(\delta T^2\right)^{-1} \leq \varepsilon_0$ and \cref{eqFixedL}, define 
\begin{align*}
	h_{1} = \frac{\varepsilon_{0} T^{2} }{ \lambda_{1} (T + T^{2} + \mathcal{H}^{2/7} T^{2}) }
	.
\end{align*} 
Then, for every $ h \leq \min\{ h_{0}, h_{1} \} $, from \cref{eqFixedL,eqObe6}, we obtain 
\begin{align*}
	\mathbb{E} \int_{\mathcal{M}} |z(0)|^{2} 
	\leq 
	e^{C(1+T^{-1}+\mathcal{H}^{2/7}+T+\mathcal{H}^{2}T)}
	\Big(
		\mathbb{E} \int_{0}^{T}\int_{G_{0} \cap \mathcal{M}} |z|^{2} d t
		+ \mathbb{E} \int_{Q} |Z|^{2} d t
		+ e^{- \frac{C}{h} } \int_{\mathcal{M}} |z|^{2} \Big|_{t=T}
	\Big)
	.
\end{align*}
\end{proof}

\section{Proof of the controllability result}
\label{sec5}


\begin{proof}[Proof of \cref{thmCon}]

For any $ y_{0} \in L^{2}_{h} (\mathcal{M}) $, define 
\begin{align*}
	J_{\varepsilon} (z_{T}) = 
	\frac{1}{2} \mathbb{E} \int_{Q} |Z|^{2} d t
	+ \frac{1}{2} \mathbb{E} \int_{0}^{T} \int_{\mathcal{M} \cap G_{0}} |z|^{2} d t
	+ \frac{\varepsilon}{2} \mathbb{E} |z_{T}|_{L^{2}_{h}(\mathcal{M})}^{2}
	- \langle y_{0}, z(0) \rangle_{L_h^2(\mathcal{M})}
	,
\end{align*}
where $ \varepsilon = e^{- \frac{C}{h}} $ and $ (z,Z) $ solves \cref{eqDiscreteEqAdjoint} with final datum $ z_{T} $.
It is clearly that $ J_{\varepsilon} (\cdot) $ is continuous and convex. 
From \cref{thmObservabilityEstiames} and Young's inequality, for $ h \leq \min\{h_{0}, h_{1}\} $, we obtain 
\begin{align*}
	& J_{\varepsilon}(z_{T}) 
	\\
	& \geq 
	\frac{1}{2} \mathbb{E} \int_{Q} |Z|^{2} d t
	+ \frac{1}{2} \mathbb{E} \int_{0}^{T} \int_{\mathcal{M} \cap G_{0}} |z|^{2} d t
	+ \frac{\varepsilon}{2} \mathbb{E}|z_{T}|_{L^{2}_{h}(\mathcal{M})}^{2}
	- \frac{1}{4 C_{obs}} \mathbb{E}|z(0)|_{L^{2}_{h}(\mathcal{M})}^{2}
	- C_{obs}  | y_{0}|_{L^{2}_{h}(\mathcal{M})}^{2}
	\\
	& \geq 
	\frac{\varepsilon}{4} \mathbb{E} |z_{T}|_{L^{2}_{h}(\mathcal{M})}^{2}
	- C_{obs} | y_{0}|_{L^{2}_{h}(\mathcal{M})}^{2}
	,
\end{align*}
which implies that $ J_{\varepsilon} (\cdot) $ is coercive. 
Hence, $ J_{\varepsilon} (\cdot) $ admits a unique minimizer $ z^{*}_{T} $.
Denote $ (z^{*}, Z^{*}) $ the solution to \cref{eqDiscreteEqAdjoint} with final datum $ z_{T}^{*} $. 
Then we obtain 
\begin{align}
	\label{eqEuler}
	\mathbb{E} \int_{Q} ZZ^{*} d t 
	+ \mathbb{E} \int_{0}^{T} \int_{\mathcal{M} \cap G_{0}} z z^{*} d t 
	+ \varepsilon \mathbb{E} \langle z_{T}, z_{T}^{*} \rangle _{L^{2}_{h}(\mathcal{M})}
	-   \langle  y_{0}, z(0) \rangle _{L^{2}_{h}(\mathcal{M})}
	=0
	,
\end{align}
for any $ z_{T} $ with the associated solution $ (z,T) $ of  \cref{eqDiscreteEqAdjoint}.
Choose $ u = - z^{*} $ and $ v = - Z^{*} $, by the duality of \cref{eqDiscreteEqOpe} with \cref{eqDiscreteEqAdjoint}, we obtain 
\begin{align*}
	\mathbb{E} \langle  y(T), z_{T} \rangle _{L^{2}_{h}(\mathcal{M})}
	-   \langle  y_{0}, z(0) \rangle _{L^{2}_{h}(\mathcal{M})}
	=
	\mathbb{E} \int_{0}^{T} \int_{\mathcal{M} \cap G_{0}} u z  d t
	+ \mathbb{E} \int_{Q} v Z d t 
	.
\end{align*}
Then $ y(T) = \varepsilon z_{T}^{*} $.
Noting that \cref{eqEuler}, we have 
\begin{align}
	\label{eqCone1}
	\mathbb{E} \int_{Q} |Z^{*}|^{2} d t 
	+ \mathbb{E} \int_{0}^{T} \int_{\mathcal{M} \cap G_{0}}  |z^{*}|^{2} d t 
	+ \varepsilon \mathbb{E} |z_{T}^{*}|_{L^{2}_{h}(\mathcal{M})}^{2}
	\leq 
	C_{obs}  | y_{0}|_{L^{2}_{h}(\mathcal{M})}^{2}
	+ \frac{1}{4 C_{obs}} | z^{*}(0)|_{L^{2}_{h}(\mathcal{M})}^{2}
	.
\end{align}
From \cref{thmObservabilityEstiames}, for $ h \leq \min\{h_{0}, h_{1}\} $, we obtain 
\begin{align}
	\label{eqCone2}
	| z^{*}(0)|_{L^{2}_{h}(\mathcal{M})}^{2}
	\leq 
	C_{obs} \Bigl(
		\mathbb{E} \int_{Q} |Z^{*}|^{2} d t 
		+ \mathbb{E} \int_{0}^{T} \int_{\mathcal{M} \cap G_{0}}  |z^{*}|^{2} d t 
		+ \varepsilon \mathbb{E} |z_{T}^{*}|_{L^{2}_{h}(\mathcal{M})}^{2}
	\Bigr)
\end{align}
Combining \cref{eqCone1,eqCone2}, for $ h \leq \min\{h_{0}, h_{1}\} $, we have 
\begin{align*}
	\mathbb{E} \int_{\mathcal{M}}|y(T)|^2 
	= \varepsilon^{2} \mathbb{E} |z_{T}^{*}|_{L^{2}_{h}(\mathcal{M})}^{2}
	\leq C_{obs} e^{-\frac{C}{h}}  \int_{\mathcal{M}}\left|y_0\right|^2
	,
\end{align*}
and 
\begin{align*}
	\mathbb{E} \int_Q|v|^2 d t
	+\mathbb{E} \int_0^T \int_{\omega \cap \mathcal{M}}|u|^2 d t 
	= \mathbb{E} \int_{Q} |Z^{*}|^{2} d t 
	+ \mathbb{E} \int_{0}^{T} \int_{\mathcal{M} \cap G_{0}}  |z^{*}|^{2} d t 
	\leq C_{obs}   \int_{\mathcal{M}}\left|y_0\right|^2
	,
\end{align*}
which implies that $ (u, v) $ satisfy our need.	
\end{proof}

\appendix

\end{document}